\documentclass[10pt,a4paper]{amsart}

\usepackage{fullpage}
\usepackage[utf8]{inputenc}
\usepackage[T2A,T1]{fontenc}
\usepackage{textcomp}
\usepackage{yfonts}

\usepackage{amsthm}
\usepackage{mathtools}
\usepackage{bm}
\usepackage{microtype}

\usepackage{amsmath}
\usepackage{amsfonts}
\usepackage{amssymb}
\usepackage{amscd}
\usepackage{stmaryrd}
\usepackage{mathrsfs}
\usepackage{accents}
\usepackage{nicefrac}

\DeclareMathAlphabet{\EuScript}{U}{eus}{m}{n}
\usepackage{comment}
\usepackage{enumerate}
\usepackage{fvextra}
\usepackage{xparse}

\usepackage[dvipsnames]{xcolor}
\usepackage{tikz}
\usetikzlibrary{arrows.meta,calc}
\usepackage[all,tips]{xy}
\usepackage{soul}
\usepackage[normalem]{ulem}
\usepackage{tcolorbox}
\tcbuselibrary{breakable,skins}

\makeatletter
\@ifclassloaded{svjour3}{}{\RequirePackage{amsthm}}
\makeatother

\usepackage{hyperref}

\makeatletter
\@ifclassloaded{svjour3}{}{%
\newtheorem{theorem}{Theorem}[section]
\newtheorem{lemma}[theorem]{Lemma}
\newtheorem{corollary}[theorem]{Corollary}
\newtheorem{proposition}[theorem]{Proposition}
\theoremstyle{definition}

\theoremstyle{remark}
\newtheorem{remark}[theorem]{Remark}
\numberwithin{equation}{section}
}
\makeatother

\newcommand{\bbB}{{\mathbb B}}
\newcommand{\bbC}{{\mathbb C}}

\newcommand{\bbN}{{\mathbb N}}

\newcommand{\bbR}{{\mathbb R}}

\newcommand{\bbZ}{{\mathbb Z}}

\newcommand{\bfF}{{\mathbf F}}

\newcommand{\bfH}{{\mathbf H}}

\newcommand{\bfL}{{\mathbf L}}

\newcommand{\bfP}{{\mathbf P}}

\newcommand{\bfW}{{\mathbf W}}

\newcommand{\calC}{{\mathcal C}}
\newcommand{\calD}{{\mathcal D}}

\newcommand{\calF}{{\mathcal F}}
\newcommand{\calG}{{\mathcal G}}

\newcommand{\calO}{{\mathcal O}}

\newcommand{\calS}{{\mathcal S}}

\newcommand{\calU}{{\mathcal U}}

\newcommand{\calX}{{\mathcal X}}

\newcommand{\calZ}{{\mathcal Z}}

\newcommand{\euA}{{\EuScript A}}

\newcommand{\euE}{{\EuScript E}}
\newcommand{\euF}{{\EuScript F}}

\newcommand{\euH}{{\EuScript H}}

\newcommand{\euK}{{\EuScript K}}

\newcommand{\euM}{{\EuScript M}}

\newcommand{\euP}{{\EuScript P}}
\newcommand{\euQ}{{\EuScript Q}}

\newcommand{\euS}{{\EuScript S}}

\newcommand{\euY}{{\EuScript Y}}

\newcommand{\sfE}{{\mathsf E}}

\newcommand{\sfN}{{\mathsf N}}

\makeatletter
\newcommand\suchthat{%
 \@ifstar
  {\mathrel{}\middle|\mathrel{}}
  {\mid}%
}
\makeatother

\newcommand{\closure}{\overline}

\renewcommand{\restriction}{|}

\newcommand{\Ball}[2]{\bbB(#1,#2)}

\newcommand{\smallball}{\bbB}

\newcommand{\diameter}[1]{\operatorname{diam}(#1)}

\newcommand{\dist}{\operatorname{dist}}

\newcommand{\dif}{{\mathrm d}}

\newcommand{\tsum}{{\textstyle{\sum}}}

\newcommand{\Id}{\operatorname{Id}}

\DeclareSymbolFont{bbold}{U}{bbold}{m}{n}
\DeclareSymbolFontAlphabet{\mathbbold}{bbold}
\newcommand{\bboldone}{\mathbbold{1}}

\newcommand{\boldPsi}{\boldsymbol{\Psi}}

\newcommand{\supp}{\operatorname{supp}}
\newcommand{\interior}[1]{\operatorname{int}{#1}}

\newcommand{\contract}{\mathbin{\lrcorner}}

\newcommand{\vecnorm}[1]{|#1|}

\DeclareMathOperator*{\esssup}{ess\,sup}

\newcommand{\cartan}{{\mathrm d}}
\newcommand{\cartanx}{{\mathrm d}x}

\newcommand{\Domain}{{\Omega}}
\newcommand{\subdomain}{{\omega}}
\newcommand{\Cover}{\calU}

\newcommand{\extDomain}{{\sfE\Omega}}
\newcommand{\collaredDomain}{{\widetilde\Omega}}
\newcommand{\Collar}{{\sfN\Omega}}

\newcommand{\Alt}[1]{\Lambda^{#1}}

\newcommand{\Proj}{\pi}
\newcommand{\Projection}{\Pi}
\newcommand{\Poincare}{\euP}

\newcommand{\Defect}{\euK}

\newcommand{\Extension}{\euE}
\newcommand{\Erweiterung}{\euA}
\newcommand{\kommErweiterung}{\euA_{\cartan}}
\newcommand{\Harmonic}{{\euH}}
\newcommand{\Mollifier}{{\euM}}
\newcommand{\Fourier}{{\euF}}

\newcommand{\weight}{\vartheta}

\newcommand{\Sobdiff}{\mathcal{W}}

\newcommand{\Wint}{W}
\newcommand{\Wres}{W_{\mathrm{res}}}

\newcommand{\grad}{\operatorname{grad}}
\newcommand{\curl}{\operatorname{curl}}

\newcommand{\divergence}{\operatorname{div}}

\newcommand{\VEC}{{V}}
\newcommand{\bounded}{\mathrm{b}}

\newcommand{\rng}{\mathrm{rng}}

\begin{document}

\title
[Universal extension de~Rham complex]
{Universal extension operators for de~Rham complexes on bounded Lipschitz domains}

\author{Martin W. Licht}
\address{TU Dresden}
\email{martin\_werner.licht@tu-dresden.de}



\begin{abstract}
    We construct universal extension operators for differential forms on bounded Lipschitz domains. These are continuous between Sobolev--Slobodeckij spaces, Besov--Triebel--Lizorkin spaces, and tempered distribution spaces of differential forms. 
    The extension operators preserve the regularity class of the coefficients and the regularity class of any weak exterior derivative.
    Regularized Poincar\'e operators play a key role in the construction.
    The first class of operators extends the differential form onto the entire ambient Euclidean space. 
    The second class of operators extends to a collar neighborhood of the domain and commutes with the exterior derivative.
    We apply these results to construct commuting smoothing operators for differential forms in numerous graph norms. 
    An exposition of the main results in the language of vector calculus is included. 
\end{abstract}

\maketitle

\setcounter{tocdepth}{4}
\tableofcontents

\section{Introduction}

Vector and tensor fields involve phenomena that are not present in the study of scalar fields. 
The de~Rham complex of differential forms is the most important example:
the exterior derivative gives rise to a natural regularity theory of differential forms that is not fully captured by the regularity theory of scalar fields.
In light of that, generalizing mathematical concepts for scalar fields to the de~Rham complex requires particular attention to remain compatible with the exterior derivative.

We address extension operators for Sobolev de~Rham complexes.
The question of how to extend scalar fields on domains onto the enclosing space has been addressed for a broad range of domains and function spaces.
But an extension operator for differential forms that is compatible with the exterior derivative is not a straightforward application of the scalar-valued result, in particular if the extension operator must commute with the exterior derivative. 
The present manuscript provides such an extension on bounded Lipschitz domains.

First, we develop extension operators for differential forms on bounded Lipschitz domains onto the entire ambient Euclidean space. 
The extensions are \emph{universal} in the sense that the same operator is continuous on multiple regularity scales of harmonic analysis: the tempered distributional differential forms, the scales of Besov and Triebel--Lizorkin spaces, and in particular the Bessel potential spaces and the Sobolev--Slobodeckij spaces. 
Importantly, the extensions are compatible with the exterior derivative in the sense that they preserve the regularity classes not only of the differential forms but also of their exterior derivatives. 

Second, we contribute commuting extension operators.
These satisfy the same continuity properties as the universal extension operators, but they extend onto a collar neighborhood of the Lipschitz domain instead of the whole space.
The operators in this second class of extensions commute with the exterior derivative and induce isomorphisms on cohomology. 
In particular, they are morphisms between Banach complexes of Sobolev--Slobodeckij spaces of differential forms.
Their properties are summarized in the following main result.

\begin{theorem}\label{introduction:theorem:commuting-universal-extension}
    Let $\Domain \subseteq \bbR^{n}$ be a bounded Lipschitz domain. 
    Then there exists a bounded Lipschitz domain $\extDomain$ with $\Domain \Subset \extDomain$ and the following properties.
    We have continuous linear operators 
    \begin{align}\label{introduction:math:commuting-universal-extension}
        \kommErweiterung^{k} : 
        \calS'(\Domain,\Alt{k})
        \to 
        \calS'(\extDomain,\Alt{k})
        ,
        \qquad 
        k \in \bbZ
        ,
    \end{align}
    such that for all $u \in \calS'(\Domain,\Alt{k})$ we have 
    \begin{align*}
        ( \kommErweiterung^{k} u) \restriction_{\Domain} = u,
        \qquad 
        \kommErweiterung^{k+1} \cartan^{k} u = \cartan^{k} \kommErweiterung^{k} u
        .
    \end{align*}
    Moreover, we have continuous operators 
    \begin{align}\label{introduction:math:commuting-universal-extension:smooth}
        &\kommErweiterung^{k} : C^{\infty}(\closure\Domain,\Alt{k}) \to C^{\infty}(\closure\extDomain,\Alt{k})
        ,
        \\
        &\kommErweiterung^{k} : \Wint^{s,p}(\Domain,\Alt{k}) \to \Wint^{s,p}(\extDomain,\Alt{k}), 
        \qquad 
        s \geq 0, \quad 1 < p < \infty
        ,
        \\
        &\kommErweiterung^{k} : B^{s}_{p,q}(\Domain,\Alt{k}) \to B^{s}_{p,q}(\extDomain,\Alt{k}),
        \qquad
        s \in \bbR, \quad 1 < p < \infty, \quad 1 < q < \infty
        ,
        \\
        &\kommErweiterung^{k} : F^{s}_{p,q}(\Domain,\Alt{k}) \to F^{s}_{p,q}(\extDomain,\Alt{k}),
        \qquad
        s \in \bbR, \quad 1 < p < \infty, \quad 1 < q < \infty
        .
    \end{align}
\end{theorem}

The distinguishing functionality is that closed differential forms that are smooth up to the boundary admit smooth closed extensions within the collar neighborhood. 
As an auxiliary result that may be of independent interest, we develop the cohomology theory of smooth differential forms on Lipschitz domains with a partial flatness condition.

Our extension operators combine two distinct components.
On the one hand, we use classical universal extension operators for scalar fields, applied componentwise to differential forms. 
On the other hand, regular decompositions of differential forms via universal potential operators play a central role.
We review the general background on these concepts.

\subsection{Universal extensions}

It is a classical topic of analysis how functions initially defined on a domain can be continued to the whole ambient space. 
The continuation of scalar fields has received the bulk of attention. 
Whitney~\cite{whitney1934analytic} discussed the continuation of $C^{m}$ functions on closed subsets of $\bbR^{n}$. 
Seeley~\cite{seeley1964extension}, generalizing earlier work by Lichtenstein~\cite{lichtenstein1929elementare}, developed a continuous linear extension operator for functions defined on the half-space that preserves smoothness of every order simultaneously. 
Important work on the extension of Sobolev spaces on Lipschitz domains up to a fixed order of smoothness is due to Calder\'on~\cite{calderon1961lebesgue}. 
Universal extension operators were developed over the following years,  
on Lipschitz domains or more general domain classes, and for different notions of function smoothness;
the term \emph{universal extension operator} seems to have been popularized in harmonic analysis.
For example, we mention work by Stein~\cite{stein1970singular}, Jones~\cite{jones1981quasiconformal}, Christ~\cite{christ1984extension}, Rychkov~\cite{rychkov1999restrictions}, Rogers~\cite{rogers2006degree}, and by Shi and Yao~\cite{shi2024estimates}.
We focus on Rychkov's main result: his extension is defined on tempered distributions and is continuous on the Besov and Triebel--Lizorkin classes of functions.

De~Rham complexes composed of rough spaces of functions, beyond the smooth spaces of classical analysis and differential geometry, 
have been researched extensively~\cite{federer,derham1984differentiable,goldshtein1982differential,scott1995Lp,iwaniec1999nonlinear,mitrea2008traces,goldshtein2011hodge}. 
Extension operators for Sobolev de~Rham complexes on Lipschitz domains must be compatible with the exterior derivative. 
Reflection-based extensions onto tubular neighborhoods of Lipschitz domains commute with the exterior derivative and preserve Lebesgue spaces of differential forms~\cite{schoeberl2008posteriori,AFW1,licht2019smoothed,licht2019mixed,hiptmair2020review},
but they generally do not preserve any smoothness, classical or otherwise, of higher order. 

The main reference for a universal extension on differential forms is due to Hiptmair, Li, and Zou~\cite{hiptmair2012universal}. 
As a direct generalization of Stein's construction from scalar fields to differential forms,
their operator preserves the regularity of the differential forms and their exterior derivatives. 
However, this universal extension operator does not commute with the exterior derivative. 
We also point out that their exposition focuses on Hilbert spaces for reasons of brevity, and only integer-order smoothness indices are covered.
These are widely used throughout the analysis of linear partial differential equations; 
our goal is to cover a wide range of function spaces, including Banach spaces of differential forms, which is pertinent to the analysis of nonlinear partial differential equations.
Their work already points out the connection of extension operators and regular decompositions.

\subsection{Regular decompositions}

Regular decompositions are a feature of the de~Rham complex that does not have an analogue in the scalar theory. 
Regular decompositions are closely related to regularizing potential operators. 
An illustrative classical result is the following decomposition into regular fields on any convex bounded Lipschitz domain $\Domain \subseteq \bbR^{3}$:
if $\bfF \in \bfL^{2}(\Domain)$ is a square-integrable vector field such that $\curl \bfF \in \bfL^{2}(\Domain)$ is square-integrable,
then there exist a scalar field $\phi \in H^{1}(\Domain)$ and a vector field $\boldPsi \in \bfH^{1}(\Domain)$,
each with first-order Sobolev regularity, such that 
\begin{align*}
    \bfF = \grad \phi + \boldPsi, \qquad \curl \bfF = \curl \boldPsi.
\end{align*}
There exist regularizing potential operators 
\begin{gather*}
    P_{\grad} : \bfL^{2}(\Domain) \to H^{1}(\Domain), \qquad \bfP_{\curl} : \bfL^{2}(\Domain) \to \bfH^{1}(\Domain)
\end{gather*}
such that 
\begin{align*}
    \phi = P_{\grad}( \bfF ), \qquad \boldPsi = \bfP_{\curl}( \curl \bfF ).
\end{align*}
On general bounded Lipschitz domains, not necessarily convex, such a decomposition still holds with a three-term sum instead of a two-term sum:
one must complement $\grad \phi$ and $\bfP_{\curl}( \curl \bfF )$ by a third summand that accounts for the cohomology of the domain.
These \emph{cohomology fields} are not gradients and have vanishing curl, but they are smooth up to the boundary of the domain.
This decomposes the field into the gradient of a more regular function, a more regular cohomology vector field, and a more regular vector field that is reconstructed as a potential of the curl of the original field.
We find regular decompositions in numerous forms throughout the literature%
~\cite{birman1987theory, iwaniec1993integral, amrouche1994decomposition, amrouche1998vector, mitrea2004sharp, hiptmair2007nodal, bourgain2007new, mitrea2008traces,   gopalakrishnan2012partial, costabel2013potential, demlow2014posteriori, bauer2016maxwell, gopalakrishnan2018auxiliary, hiptmair2019discrete, baldi2019l1, hiptmair2020review, pauly2022hilbert, zbMATH07859130, goldshtein2026some}. %

A significant milestone is the theory of averaged Poincar\'e and Bogovskii operators on bounded Lipschitz domains, 
developed by Iwaniec and Lutoborski \cite{iwaniec1993integral}, by Costabel and McIntosh~\cite{costabel2010}, and by Mitrea, Mitrea, and Monniaux~\cite{mitrea2008poisson}, among others. 
These potential operators are pseudodifferential operators of order $-1$ with symbols in the local H\"ormander class $S^{-1}_{1,0}(\bbR^{n})$. 
This fact already warrants a gain of smoothness and continuity across a scale of locally convex function spaces simultaneously. 
In particular, this enables regular decompositions for a large class of de~Rham complexes. 

The averaged Poincar\'e operators can be interpreted as generalized inverses. 
In the case of Hilbert complexes, we distinguish them from the operator pseudoinverses of the exterior derivative between Hilbert spaces of differential forms, which on exact forms always provide the preimage with minimum norm. 
The averaged Poincar\'e operators, by contrast, have no such norm-minimizing property. 
More generally, the averaged Poincar\'e operators are generalized inverses for Banach complexes when general Sobolev--Slobodeckij spaces are considered.

In analogy with the universal extension operators, 
Costabel and McIntosh's potential operators are henceforth referred to as \emph{universal potential operators}.
They define a \emph{universal regularized decomposition}, a generalized Hodge decomposition of differential forms valid simultaneously for different function spaces.
This decomposition holds not only for Hilbert complexes of differential forms, but also for Banach complexes and most generally for certain distributional de~Rham complexes of locally convex spaces. 

\subsection{Construction and application}

We now point out a few features of our constructions. 
The first operator is a comparatively simple combination of the universal regular decomposition and the universal extension operators. 
Componentwise application of the universal extension operator to the regular decomposition terms then produces the desired extension of differential forms onto the ambient space. 
This operation is compatible with the exterior derivative: 
the extension onto the enclosing space inherits the regularity of the original differential form, and likewise its exterior derivative inherits the regularity of the original exterior derivative.

We need to appreciate the cohomology of the de~Rham complex to develop our commuting extension operator. 
In a differential complex, the range of one differential operator is included in the kernel of the next differential operator. 
For example, this structure is present in the three-dimensional $L^{2}$ de~Rham complex 
\begin{align*}
    \begin{CD}
        0 \to H^{1}(\Domain) @>\grad>> \bfH(\curl,\Domain) @>\curl>> \bfH(\divergence,\Domain) @>\divergence>> L^{2}(\Domain) \to 0
    \end{CD}
\end{align*}
on a bounded Lipschitz domain $\Domain$.
Generally, the range is a proper subspace of the kernel. 
The questions arise whether we can identify the dimension of the quotient space and whether the range is topologically complemented in the kernel. 

The quotient spaces are known as the cohomology spaces of the de~Rham complex.
We show that they are finite-dimensional and that their dimensions are the Betti numbers of the domain. 
In fact, more is true: the ranges have a finite-dimensional complement spanned by a basis of smooth differential forms. 
In particular, the closed differential forms allow a decomposition into an exact differential form and this finite-dimensional complement.
This holds for the common Hilbert complexes of differential forms, and more generally for the Banach complexes and Fr\'echet complexes. 

The universal regular decomposition involves a \emph{universal harmonic complement}, which is a finite-dimensional topological complement to the ranges of the differential operators inside the kernel of the next differential, simultaneously for all regularity classes of the source.
If we pick an arbitrary smooth extension for each of these finite-dimensional harmonic complements, 
then this leads to an extension operator on the full space. 
However, the resulting extension operators may not preserve the cocycle property of these harmonic forms, which obstructs the commutativity with the exterior derivative. 
For example, the non-zero constant functions span the kernel of the gradient operator but admit no extension onto $\bbR^{n}$ that is both gradient-free (i.e., constant) and square-integrable. Another example is that the extension of harmonic complements in $\bfH(\curl,\Domain)$ may not be globally curl-free.

This is resolved by our second class of extension operators, which modifies the machinery previously described.
If we consider universal extensions onto a small collar neighborhood $\extDomain$ of the bounded Lipschitz domain $\Domain$, but not onto the whole ambient space, 
then a prudent choice of the extension of harmonic forms ensures that the extension operators commute with the exterior derivative.
This gives commuting diagrams such as 
\begin{align*}
    \begin{CD}
        H^{s}(\Domain) @>\grad>> \bfH^{s}(\curl,\Domain) @>\curl>> \bfH^{s}(\divergence,\Domain) @>\divergence>> H^{s}(\Domain)
        \\
        @V{\kommErweiterung^{0}}VV 
        @V{\kommErweiterung^{1}}VV 
        @V{\kommErweiterung^{2}}VV 
        @V{\kommErweiterung^{3}}VV
        \\
        H^{s}(\extDomain) @>\grad>> \bfH^{s}(\curl,\extDomain) @>\curl>> \bfH^{s}(\divergence,\extDomain) @>\divergence>> H^{s}(\extDomain).
    \end{CD}
\end{align*}
Here, we enlarge the original domain by a small neighborhood such that the topology remains the same and the Betti numbers are preserved. 
This observation motivates the search for a commuting extension operator from the bounded Lipschitz domain onto the enlarged domain that is then an isomorphism on cohomology.

The decisive insight is that smooth closed forms, such as in the universal harmonic complement, admit a smooth closed extension onto the collar neighborhood. This follows from the exactness of the smooth de~Rham complex of the outer collar neighborhood $\extDomain \setminus \closure\Domain$ when imposing flat boundary conditions along the interior Lipschitz interface $\partial\Domain$. We develop the following auxiliary result through the apparatus of sheaf theory~\cite{cap1998garben,theeigenbook,mitrea2008poisson,goldshtein2011hodge,deitmar2025kohomologie}.

\begin{lemma}\label{introduction:lemma:smooth-extension}
    Let $\Domain \subseteq \bbR^{n}$ be a bounded Lipschitz domain. 
    Then there exists a bounded Lipschitz domain $\extDomain$ with $\Domain \Subset \extDomain$ 
    such that for any $u \in C^{\infty}(\closure{\Domain},\Alt{k})$ with $\cartan^{k} u = 0$
    there exists $v \in C^{\infty}(\closure{\extDomain},\Alt{k})$ with $\cartan^{k} v = 0$
    and $v \restriction_{\Domain} = u$. 
\end{lemma}

De~Rham complexes on the $L^{2}$ scale have received the largest share of attention; Costabel and McIntosh focus on these operators for simplicity of exposition.
As we leverage pseudodifferential calculus, we cover a much wider class of function spaces: for the range of Lebesgue exponents $1 < p < \infty$, we recover the full Sobolev--Slobodeckij regularity scale on Lipschitz domains, including spaces with fractional smoothness and the negative-order spaces.
In fact, we utilize de~Rham complexes of locally convex function spaces, such as tempered distributions or bounded smooth forms.
This allows us to separate the algebraic properties of the operators, such as homotopy formulas, from their continuity in particular function spaces. 
Within the pseudodifferential calculus, we can extend the discussion to the Besov and Triebel--Lizorkin classes with little additional overhead. 
To accommodate a variety of different readers, our presentation is modular: 
continuity statements are given for the common Sobolev--Slobodeckij spaces and the Besov--Triebel--Lizorkin spaces separately.

Extension operators are convenient because they make techniques for spaces on Euclidean space available to function spaces on domains.
For example, a direct application of these results is a commuting smooth approximation operator on Sobolev de~Rham complexes:
a bounded family of smoothing operators
from a Sobolev--Slobodeckij de~Rham complex with additional smoothness assumptions into the smooth de~Rham complex on a bounded Lipschitz domain
that converges pointwise to the identity.
These give commuting diagrams such as 
\begin{align*}
    \begin{CD}
        W^{s,p}(\Domain) @>\grad>> \bfW^{s,p}(\curl,\Domain) @>\curl>> \bfW^{s,p}(\divergence,\Domain) @>\divergence>> W^{s,p}(\Domain)
        \\
        @V{\Mollifier^{0}_{\epsilon,\grad}}VV @V{\Mollifier^{1}_{\epsilon,\curl}}VV @V{\Mollifier^{2}_{\epsilon,\divergence}}VV @V{\Mollifier^{3}_{\epsilon}}VV
        \\
        C^{\infty}(\closure\Domain) @>\grad>> C^{\infty}(\closure\Domain,\bbR^{3}) @>\curl>> C^{\infty}(\closure\Domain,\bbR^{3}) @>\divergence>> C^{\infty}(\closure\Domain)
    \end{CD}
\end{align*}
Here, we are concerned with non-standard smoothness spaces such as 
\begin{align}\label{math:Wspcurl}
    \bfW^{s,p}(\curl,\Domain) 
    := 
    \{ \mathbf{u} \in \bfW^{s,p}(\Domain) 
    \suchthat 
    \curl \mathbf{u} \in \bfW^{s,p}(\Domain) \}
    ,
\end{align}
where the curl of a vector field provides additional smoothness. 
While the mollification of a componentwise extension of these vector fields is a standard smoothing argument, 
this does not warrant convergence in~\eqref{math:Wspcurl} because the componentwise extension generally does not have a curl with the improved smoothness required by $\bfW^{s,p}(\curl,\Domain)$. 
But if we use the compatible extension operator, then smoothing produces an approximation in the correct norm. 
In particular, 
\begin{align}\label{math:smooth-forms-are-dense-in-Wspcurl}
    C^{\infty}(\closure\Domain,\bbR^{3})
    \subseteq
    \bfW^{s,p}(\curl,\Domain) 
\end{align}
is a dense subspace.

Universal extension operators for de~Rham complexes allow us to define operations on bounded Lipschitz domains in terms of operations on the ambient space. 
Even though we use universal scalar extension and universal potential operators, the construction principle can be implemented with different scalar extension operators as a foundation. Our methodology can be instantiated with any scalar extension operator and regular decomposition under the condition that the regularity parameters interact accordingly for the parameter regime of interest. 
The class of Lipschitz domains is very general while still ensuring convenient mathematical structure. However, this class of domains is by no means the only class relevant for mathematical modeling.

\subsection{Context}

The work in the present manuscript was undertaken to generalize an approximation result that has been used in the context of finite element quasi-interpolation operators~\cite{licht2023averaging}. 

The connection between extension operators and regular decompositions has been recognized for some time.
Indeed, Hiptmair, Li, and Zou~\cite{hiptmair2012universal} use the generalized Stein extension operators to construct regular decompositions (and Poincar\'e inequalities) from the regular decompositions on the ambient space. 
The possibility of using regular decompositions in the construction of extension operators seems to have been anticipated in their work~\cite[Remark~3.2]{hiptmair2012universal}. The present manuscript seems to be the first to elaborate on that idea.
We also point out recent work by Nilsson and Pitassi~\cite{nilsson2026uniformly} that has addressed the extension problem for $L^{2}$ de~Rham complexes to obtain uniform Poincar\'e inequalities relevant for unfitted finite element methods. Their extension operator is defined on the orthogonal complement of the harmonic forms. By contrast, our construction intrinsically relies on Poincar\'e inequalities because it starts from universal regular decompositions, and the harmonic forms are extended to smooth closed forms. 
We therefore draw the overall picture that extension operators and regular decompositions for de~Rham complexes on bounded Lipschitz domains are conceptually equivalent.

Sobolev--Slobodeckij spaces are of concrete interest in the numerical analysis of partial differential equations, where smoothness, even if only partial, is reflected in the provable convergence rates of finite element methods. This especially pertains to numerical methods in vector calculus on domains, where solutions to partial differential equations feature non-integer smoothness that is typically rougher than the smoothness of scalar partial differential equations.
While Sobolev--Slobodeckij de~Rham complexes are fundamental for partial differential equations, discrete de~Rham complexes now play a central role in numerical analysis. 
Alain Bossavit~\cite{bossavit1992new,rapetti2009whitney} identified Hassler Whitney's eponymous forms as finite element differential forms, and pioneering contributions by Ralf Hiptmair~\cite{hiptmair2002finite} and Douglas Arnold, Richard Falk, and Ragnar Winther~\cite{AFW1} have established finite element de~Rham complexes as paradigmatic in mixed finite element methods. 

While Hilbert complexes have been the subject of increased research in recent years, there is a long tradition of studying Banach complexes and distributional de~Rham complexes, which are differential complexes of locally convex spaces~\cite{palamodov1972homological,khelemskij1989homology,kuzminov1998homological,kuzminov2003compactness}. 
Even though these are generally not situated in the theory of abelian categories, the cohomology theory of these differential complexes is well understood.

\subsection{Organization of the manuscript}

The manuscript is organized as follows.
Section~\ref{section:function-spaces-on-domains} reviews the function spaces under consideration.
Section~\ref{section:extension-operators-on-lipschitz-domains} reviews the first main component of our proof, namely universal extension operators.
Section~\ref{section:differential-forms} reviews notions of exterior calculus and spaces of differential forms.
Section~\ref{section:regular-decomposition} reviews the second main component, namely the universal regular decomposition. 
Section~\ref{section:extension-of-forms} first constructs the universal extension operator compatible with the exterior derivative and then provides the commuting extension operator.
Finally, Section~\ref{section:vector-calculus} presents the main results in the language of vector calculus and for Sobolev--Slobodeckij spaces.

\subsection*{Notation}

For any set $X \subseteq \bbR^{n}$, we write $\bboldone_{X} : \bbR^{n} \to \{0,1\}$ for the indicator function of $X$. 
We use the notation $A(n):=\bbN_{0}^{n}$ for the set of multi-indices over $\{1,\dots,n\}$.
For $\alpha = (\alpha_1, \dots, \alpha_n) \in A(n)$, we set
\begin{gather*}
    |\alpha| := \alpha_1 + \dots + \alpha_n,
    \qquad
    \partial^\alpha
    :=
    \partial_1^{\alpha_1} \cdots \partial_n^{\alpha_n},
    \qquad 
    x^{\alpha} := x_1^{\alpha_1} \cdots x_n^{\alpha_n}.
\end{gather*}

\section{Function spaces on domains}\label{section:function-spaces-on-domains}

This section introduces the function spaces and the corresponding notation that are used throughout the remainder of the manuscript. 
All definitions hold for functions on domains taking values in a finite-dimensional real inner product space. 
Of course, the corresponding vector-valued function spaces can be defined equally as tuples of functions, one for each component, with the respective regularity classes. 

Throughout this section, $\Domain \subseteq \bbR^{n}$ is a domain, i.e., a connected open set.

\subsection{Spaces of smooth functions and distributions}

We consider spaces of functions taking values in a finite-dimensional inner product space $\VEC$ with norm $\vecnorm{\cdot}$. We identify $\VEC$ with its dual space $\VEC^{*}$ via the inner product.

We write $C^{\infty}(\Domain,\VEC)$ for the space of smooth functions on $\Domain$.
Furthermore, 
we write $C^{\infty}_{\bounded}(\Domain,\VEC)$ for the space of smooth functions on $\Domain$ whose derivatives of all orders are bounded,
we write $C^{\infty}(\closure\Domain,\VEC)$ for the space of smooth functions on $\Domain$ that are restrictions of functions in $C^{\infty}(\bbR^{n},\VEC)$, 
and we write $C_{c}^{\infty}(\Domain,\VEC)$ for the space of smooth functions with compact support in $\Domain$. 

While we mostly consider these spaces as vector spaces, they can also be equipped with the structure of locally convex spaces.
The space $C^{\infty}(\Domain,\VEC)$ is a locally convex space by the family of seminorms
that induce the topology of uniform convergence of all derivatives on compact sets:
\begin{align*}
    p_{K,\alpha}(u) := \max_{x \in K} \vecnorm{ \partial^{\alpha} u(x) }, \qquad K \Subset \Domain, \quad \alpha \in A(n)
    .
\end{align*}
In fact, this topology is generated by a countable family of seminorms:
if $\emptyset = K_{0} \subseteq K_{1} \subseteq K_{2} \subseteq \dots$ is an exhaustion of $\Domain$ by compact sets, i.e., 
\begin{align*}
    \bigcup_{i \in \bbN } K_{i} = \Domain,
    \qquad 
    \forall i \in \bbN : K_{i-1} \subset \interior{K_{i}}
    ,
\end{align*}
then the corresponding countable family of seminorms induces the topology of $C^{\infty}(\Domain,\VEC)$. 
The space $C^{\infty}(\closure\Domain,\VEC)$ of restrictions of smooth functions 
\begin{align*}
    C^{\infty}(\closure\Domain,\VEC) := \left\{ u \restriction_{\Domain} \suchthat u \in C^{\infty}(\bbR^{n},\VEC) \right\}
\end{align*}
carries a natural quotient topology with a corresponding countable family of seminorms.
The space $C^{\infty}_{\bounded}(\Domain,\VEC)$ is locally convex with the countable family of seminorms that describe the topology of (global) uniform convergence of all derivatives:
\begin{align*}
    p_{\alpha}(u) := \sup_{x \in \Domain} \vecnorm{ \partial^{\alpha} u(x) }, \qquad \alpha \in A(n).
\end{align*}

Finally, the space $C^{\infty}_{c}(\Domain,\VEC)$ carries the usual locally convex topology; we follow the convention of writing $\calD(\Domain,\VEC)$ for that locally convex vector space. An explicit family of seminorms is constructed as follows: if $\emptyset = K_{0} \subseteq K_{1} \subseteq K_{2} \subseteq \dots \subseteq \Domain$ is some compact exhaustion of the domain $\Domain$, 
then 
\begin{align*}
    p_{m,c}(u) 
    := 
    \sup_{ i \in \bbN } \sup_{ x \in K_{i} \setminus K_{i-1} } 
    \max_{ \substack{ \alpha \in A(n) \\ |\alpha| \leq m_i } } 
    c_{i} \vecnorm{ \partial^{\alpha} u(x) }
    ,
    \qquad 
    \{m_{i}\}_{i\in\bbN}, \{c_{i}\}_{i\in\bbN} \in \bbN^{\bbN}
\end{align*}
is a family of seminorms 
that 
induces the topology of $\calD(\Domain,\VEC)$. This topology is characterized as an inductive limit. A sequence converges in this topology if its terms are supported within a common compact subset of $\Domain$ and all derivatives converge uniformly.

We review a few spaces of distributions. 
We let $\calD'(\Domain,\VEC)$ denote the dual space of $\calD(\Domain,\VEC)$, known as the space of \emph{distributions} on $\Domain$ with values in $\VEC$.
We follow the convention that the dual space carries the topology of uniform convergence on bounded sets. 
Recall that a subset of a locally convex space is called bounded if every continuous seminorm has a finite upper bound on that subset.
The topology of $\calD'(\Domain,\VEC)$ is induced by the seminorms
\begin{align*}
    p_{B}(u) := \sup_{ \phi \in B } \lvert \langle u, \phi \rangle \rvert, \qquad B \subseteq \calD(\Domain,\VEC) \text{ bounded}.
\end{align*}

The following smooth function space and its distributional dual are particularly important in harmonic analysis.
Recall the Schwartz space on $\bbR^{n}$, which is the vector space of fast-decaying functions
\begin{gather*}
    \calS(\bbR^{n},\VEC) 
    := 
    \left\{ 
        u \in C^{\infty}_{\bounded}(\bbR^{n},\VEC) 
        \suchthat* 
        \forall \alpha, \beta \in A(n) : x^{\alpha} \partial^{\beta} u \in C^{\infty}_{\bounded}(\bbR^{n},\VEC) 
    \right\} 
\end{gather*}
equipped with the locally convex topology generated by the seminorms
\begin{gather*}
    \|u\|_{\alpha,\beta}
    :=
    \sup_{x \in \bbR^{n}} \vecnorm{x^{\alpha} \partial^{\beta} u(x)}
    ,
    \qquad 
    \alpha,\beta \in A(n)
    .
\end{gather*}
The Schwartz space has a topological dual space $\calS'(\bbR^{n},\VEC)$, known as the space of \emph{tempered distributions}, which is again a locally convex vector space.
The tempered distributions, too, carry the topology of uniform convergence on bounded sets. 

For any domain $\Domain \subseteq \bbR^{n}$, we define 
\begin{align*}
    \calS'(\Domain,\VEC) 
    := 
    \left\{ u \restriction_{\Domain} \suchthat* u \in \calS'(\bbR^{n},\VEC) \right\} 
    \subseteq 
    \calD'(\Domain,\VEC)
    .
\end{align*}
This is a locally convex vector space equipped with the quotient topology. This space embeds into $\calD'(\Domain,\VEC)$. 
We will write $\langle u, \vartheta \rangle$ for the duality pairing of a distribution $u \in \calS'(\Domain,\VEC)$ and a test function $\vartheta \in \calD(\Domain,\VEC)$.

With respect to these topologies, the aforementioned spaces satisfy the continuous inclusions 
\begin{gather*}
    \calD(\bbR^{n},\VEC) \subseteq \calS(\bbR^{n},\VEC) \subseteq C^{\infty}(\bbR^{n},\VEC), 
    \qquad 
    \calS'(\Domain,\VEC) \subseteq \calD'(\Domain,\VEC), 
    \\
    \calD(\Domain,\VEC) \subseteq C^{\infty}(\closure\Domain,\VEC), 
    \qquad 
    C^{\infty}(\closure\Domain,\VEC) \subseteq C^{\infty}(\Domain,\VEC), 
    \\ 
    \calD(\Domain,\VEC) \subseteq C^{\infty}_{\bounded}(\Domain,\VEC), 
    \qquad 
    C^{\infty}_{\bounded}(\Domain,\VEC) \subseteq C^{\infty}(\Domain,\VEC).
\end{gather*}

\subsection{Lebesgue spaces and Sobolev spaces}

For $1 \leq p < \infty$, we write $L^{p}(\Domain,\VEC)$ for the Banach space of measurable functions $u : \Domain \to \VEC$, 
up to equality almost everywhere, for which
\begin{gather}\label{function-spaces-on-domains:lp-norm}
    \| u \|_{L^{p}(\Domain,\VEC)}
    :=
    \left(
        \int_{\Domain} \vecnorm{u(x)}^{p} \,\dif x
    \right)^{1/p}
    < \infty
    .
\end{gather}
For $p = \infty$, we write $L^{\infty}(\Domain,\VEC)$ for the Banach space of essentially bounded measurable functions, modulo equivalence almost everywhere, equipped with the norm
\begin{gather}\label{function-spaces-on-domains:linfty-norm}
    \|u\|_{L^{\infty}(\Domain,\VEC)}
    :=
    \esssup_{x \in \Domain} \vecnorm{u(x)}.
\end{gather}
We first introduce Sobolev spaces on domains in terms of weak derivatives and review the Sobolev--Slobodeckij spaces in terms of the Gagliardo--Slobodeckij seminorms. 
For any $m \in \bbN_{0}$ and $1 \leq p \leq \infty$, we define the \emph{Sobolev space} $\Wint^{m,p}(\Domain,\VEC)$ by setting 
\begin{gather*}
    \Wint^{m,p}(\Domain,\VEC)
    :=
    \left\{
        u \in L^{p}(\Domain,\VEC)
        \suchthat
        \forall \alpha \in A(n),\, |\alpha| \leq m : \partial^\alpha u \in L^{p}(\Domain,\VEC)
    \right\}.
\end{gather*}
With this choice of notation, $\Wint^{0,p}(\Domain,\VEC) = L^{p}(\Domain,\VEC)$. 

For $1 \leq p < \infty$, we equip $\Wint^{m,p}(\Domain,\VEC)$ with the norm
\begin{gather*}
    \|u\|_{\Wint^{m,p}(\Domain,\VEC)}
    :=
    \left(
        \tsum_{|\alpha| \leq m}
        \|\partial^\alpha u\|_{L^{p}(\Domain,\VEC)}^{p}
    \right)^{1/p}.
\end{gather*}
For $p = \infty$, we use the norm 
\begin{gather*}
    \|u\|_{\Wint^{m,\infty}(\Domain,\VEC)}
    :=
    \sum_{|\alpha|\leq m}
    \|\partial^\alpha u\|_{L^{\infty}(\Domain,\VEC)}.
\end{gather*}
With these norms, the spaces $\Wint^{m,p}(\Domain,\VEC)$ are Banach spaces.

Let $0 < \theta < 1$ and $1 \leq p \leq \infty$.
In the case $p < \infty$, the Slobodeckij seminorm of any $u \in L^{p}(\Domain,\VEC)$ is defined as 
\begin{gather*}
    |u|_{\Wint^{\theta,p}(\Domain,\VEC)}
    :=
    \left(
        \int_{\Domain} \int_{\Domain}
        \frac{ \vecnorm{u(x)-u(y)}^{p} }{ \vecnorm{x-y}^{n+\theta p} }
        \,\dif x\,\dif y
    \right)^{1/p}.
\end{gather*}
In the case $p = \infty$, the Slobodeckij seminorm of any $u \in L^{\infty}(\Domain,\VEC)$ is defined as 
\begin{gather*}
    |u|_{\Wint^{\theta,\infty}(\Domain,\VEC)}
    :=
    \esssup_{ \substack{ x, y \in \Domain \\ x \neq y } } 
    \frac{ \vecnorm{u(x)-u(y)} }{ \vecnorm{x-y}^{\theta} }
    .
\end{gather*}
The fractional Sobolev--Slobodeckij space
\begin{gather*}
    \Wint^{\theta,p}(\Domain,\VEC)
    :=
    \left\{
        u \in L^{p}(\Domain,\VEC)
        \suchthat
        |u|_{\Wint^{\theta,p}(\Domain,\VEC)} < \infty
    \right\}
\end{gather*}
is a Banach space with norm
\begin{gather*}
    \|u\|_{\Wint^{\theta,p}(\Domain,\VEC)}
    :=
    \|u\|_{L^{p}(\Domain,\VEC)}
    +
    |u|_{\Wint^{\theta,p}(\Domain,\VEC)}.
\end{gather*}

More generally, for $m \in \bbN_{0}$, $0 < \theta < 1$, and $1 \leq p \leq \infty$,
the Sobolev--Slobodeckij space of order $m+\theta$ is defined by 
\begin{gather*}
    \Wint^{m+\theta,p}(\Domain,\VEC)
    :=
    \left\{
        u \in \Wint^{m,p}(\Domain,\VEC)
        \suchthat %
        \forall \alpha \in A(n) \text{ with } |\alpha|=m :
        \partial^\alpha u \in \Wint^{\theta,p}(\Domain,\VEC)
    \right\}.
\end{gather*}
This is a Banach space with the norm 
\begin{gather*}
    \|u\|_{\Wint^{m+\theta,p}(\Domain,\VEC)}
    :=
    \|u\|_{\Wint^{m,p}(\Domain,\VEC)}
    +
    \sum_{ \substack{ \alpha \in A(n) \\ |\alpha| = m } }
    |\partial^\alpha u|_{\Wint^{\theta,p}(\Domain,\VEC)}.
\end{gather*}

The Sobolev--Slobodeckij spaces on domains $\Wint^{s,p}(\Domain,\VEC)$ as defined above are to be distinguished from the Sobolev--Slobodeckij spaces on domains defined in terms of restrictions.
The restriction-based definition of Sobolev--Slobodeckij spaces gives us access to the machinery of harmonic analysis and pseudodifferential calculus, but its equivalence to the intrinsic definition depends on the domain.
We focus on the reflexive spaces.

The Sobolev--Slobodeckij spaces on $\bbR^{n}$ with negative smoothness are defined as a duality construction. 
If $s \in (0,\infty)$ and $1 < p,q < \infty$ with $\frac{1}{p} + \frac{1}{q} = 1$,
then we define the dual Banach space 
\begin{align*}
    \Wint^{-s,p}(\bbR^{n},\VEC) := \Wint^{s,q}(\bbR^{n},\VEC)'.
\end{align*}
This Banach space is equipped with the usual operator norm. 

For $s \in \bbR$ and $1 < p < \infty$,
the restriction-based Sobolev--Slobodeckij space $\Wres^{s,p}(\Domain,\VEC)$ on any open set $\Domain$ is defined as the space of restrictions of distributions in $\Wint^{s,p}(\bbR^{n},\VEC)$ to $\Domain$.
The space $\Wres^{s,p}(\Domain,\VEC)$ is equipped with the quotient norm
\begin{align*}
    \|u\|_{\Wres^{s,p}(\Domain,\VEC)}
    :=
    \inf\left\{
        \| \widetilde{u} \|_{\Wint^{s,p}(\bbR^{n},\VEC)}
        \suchthat*
        \widetilde{u} \in \Wint^{s,p}(\bbR^{n},\VEC), \; \widetilde{u} \restriction_{\Domain} = u
    \right\}.
\end{align*}
It is evident from the construction that 
\begin{align*}
    \Wres^{s,p}(\bbR^{n},\VEC) = \Wint^{s,p}(\bbR^{n},\VEC)
    .
\end{align*}
For non-negative smoothness, an analogous equality of Banach spaces, with equivalent norms, holds on sufficiently regular domains.
We will prove it for bounded Lipschitz domains.

Partial derivatives are defined either in the weak sense or in the sense of duality:
if $s \in \bbR$, $m \in \bbN_{0}$, $1 < p < \infty$, and $\alpha \in A(n)$ with $|\alpha| \leq m$, then we have a bounded linear operator
\begin{align*}
    \partial^{\alpha} : \Wres^{s,p}(\Domain,\VEC) \to \Wres^{s-m,p}(\Domain,\VEC)
    .
\end{align*}
Moreover, we will use 
\begin{align*}
    \Wres^{s,p}(\Domain,\VEC) \subseteq \calD'(\Domain,\VEC), \qquad s \in \bbR, \quad 1 < p < \infty.
\end{align*}

\begin{remark}
    We have first defined the intrinsic Sobolev--Slobodeckij spaces on domains in terms of weak derivatives and the Gagliardo--Slobodeckij seminorms. 
    This is the definition most common in the theory of partial differential equations, which is our main application. 
    The definition in terms of restrictions of Sobolev--Slobodeckij spaces on $\bbR^{n}$ is common in harmonic analysis. 
    On general open sets, these two definitions of Sobolev--Slobodeckij spaces are generally not equivalent, even at positive or integer orders.
    However, this dissonance resolves if the domain is a bounded Lipschitz domain:
    then the definitions agree for non-negative smoothness parameters.
    This observation is standard (e.g.,~\cite{grisvard2011elliptic}) and will be formally reviewed later in Section~\ref{section:extension-operators-on-lipschitz-domains}.
    
    The situation for the Sobolev--Slobodeckij spaces with negative smoothness parameter is more delicate. 
    We have defined the spaces with negative smoothness as the restriction (in the sense of distributions) of the dual spaces. 
    These need not coincide with the duals of the corresponding positive-order restriction spaces, even on bounded Lipschitz domains.
    This issue is not as pressing for the purpose of this manuscript, though, 
    because we use negative-order Sobolev--Slobodeckij spaces only in the sense of restrictions to apply the methods of harmonic analysis.
    We point out that, by definition, for any $s > 0$ and $1 < p,q < \infty$ with $\frac{1}{p} + \frac{1}{q} = 1$,
    the space $\Wres^{-s,p}(\Domain,\VEC)$ coincides with the dual space of the closure of $C^{\infty}_{c}(\Domain,\VEC) \subseteq C^{\infty}_{c}(\bbR^{n},\VEC)$ within the Banach space $\Wint^{s,q}(\bbR^{n},\VEC)$. 
\end{remark}

\subsection{Besov and Triebel--Lizorkin spaces}

Sobolev--Slobodeckij spaces on Euclidean space and their restrictions to domains have a natural generalization in terms of Besov and Triebel--Lizorkin spaces.
The Besov and Triebel--Lizorkin spaces on $\bbR^{n}$ are quasi-Banach spaces that continuously embed into the space of tempered distributions.
We first define them on $\bbR^{n}$ and then on arbitrary domains in terms of restriction.
The following outline is based on the exposition in~\cite{mitrea2008poisson}.

We need the notion of a Littlewood--Paley decomposition of unity, which describes a partition of unity in the Fourier domain.
We begin with a smooth non-negative radially symmetric function $\zeta_{0} \in C^{\infty}_{c}(\bbR^{n},\bbR)$ that satisfies
\begin{gather*}
    \forall \xi \in \Ball{0}{1} : \zeta_{0}(\xi) = 1,
    \qquad
    \forall \xi \in \Ball{0}{2} : \zeta_{0}(\xi) > 0,
    \qquad
    \forall \xi \notin \Ball{0}{2} : \zeta_{0}(\xi) = 0,
\end{gather*}
and that we assume to be strictly decreasing in $|\xi|$ on the annulus $\Ball{0}{2} \setminus \Ball{0}{1}$.

For $j \in \bbN$, we define $\zeta_{j} \in C^{\infty}_{c}(\bbR^{n},\bbR)$ through
\begin{align*}
    \zeta_{j}(\xi) := \zeta_{0}( 2^{-j} \xi ) - \zeta_{0}( 2^{1-j} \xi )
    .
\end{align*}
This entails the support properties
\begin{gather*}
    \forall \xi \in \Ball{0}{2^{j+1}} \setminus \closure{\Ball{0}{2^{j-1}}} 
    : 
    \zeta_{j}(\xi) > 0,
    \qquad
    \forall \xi \notin \Ball{0}{2^{j+1}} \setminus \closure{\Ball{0}{2^{j-1}}} 
    : 
    \zeta_{j}(\xi) = 0,
\end{gather*}
and the derivative bounds 
\begin{align*}
    \forall \alpha \in A(n) 
    :
    \sup_{j \in \bbN_{0}} 
    \sup_{\xi \in \bbR^{n}} 
    2^{j|\alpha|}
    | \partial^{\alpha} \zeta_{j}(\xi) |
    \leq 
    ( 1 + 2^{|\alpha|} ) 
    \sup_{\xi \in \bbR^{n}} 
    | \partial^{\alpha} \zeta_{0}(\xi) |
    .
\end{align*}
These functions form a smooth partition of unity
\begin{align*}
    \sum_{j=0}^{\infty} \zeta_{j}(\xi) = 1, \qquad \xi \in \bbR^{n}
    .
\end{align*}
We denote the Fourier transform and its inverse by $\Fourier$ and $\Fourier^{-1}$, respectively.
They act on Schwartz functions with values in the complexification $\VEC_{\bbC} := \VEC \otimes_{\bbR} \bbC$ as continuous operators:
\begin{align*}
    \Fourier      &: \calS(\bbR^{n},\VEC_{\bbC}) \to \calS(\bbR^{n},\VEC_{\bbC}),
    \qquad 
    u \mapsto \int_{\bbR^{n}} u(x) \, e^{-2\pi i x \cdot \xi} \,\dif x,
    \\
    \Fourier^{-1} &: \calS(\bbR^{n},\VEC_{\bbC}) \to \calS(\bbR^{n},\VEC_{\bbC}),
    \qquad
    u \mapsto \int_{\bbR^{n}} u(x) \, e^{2\pi i x \cdot \xi} \,\dif x.
\end{align*}
Duality gives continuous operators $\Fourier : \calS'(\bbR^{n},\VEC_{\bbC}) \to \calS'(\bbR^{n},\VEC_{\bbC})$ and $\Fourier^{-1} : \calS'(\bbR^{n},\VEC_{\bbC}) \to \calS'(\bbR^{n},\VEC_{\bbC})$.

The Littlewood--Paley decomposition of unity is now constituted in terms of Fourier multiplier operators.
For $j \in \bbN_{0}$, we define the continuous linear Fourier multiplier operators
\begin{align*}
    \Delta_{j} : \calS'(\bbR^{n},\VEC) \to \calS'(\bbR^{n},\VEC),
    \quad 
    u \mapsto \Fourier^{-1}\left( \zeta_{j} \Fourier u \right).
\end{align*}
Since the cutoffs $\zeta_j$ are real and even, $\Delta_j$ preserves real $\VEC$-valued distributions.
These operators satisfy the decomposition property
\begin{align*}
    u = \sum_{j=0}^{\infty} \Delta_{j} u, \qquad u \in \calS'(\bbR^{n},\VEC).
\end{align*}
Here, the series converges in the topology of $\calS'(\bbR^{n},\VEC)$.

We can now define the Besov spaces. 
Let $s \in \bbR$, $0 < p \leq \infty$, and $0 < q \leq \infty$.
The Besov space $B^{s}_{p,q}(\bbR^{n},\VEC)$ is the quasi-Banach space of tempered distributions $u \in \calS'(\bbR^{n},\VEC)$ such that 
\begin{align}\label{math:besov-quasinorm}
    \| u \|_{B^{s}_{p,q}(\bbR^{n},\VEC)}
    :=
    \left\|
        \left(
            2^{js} \| \Delta_{j} u \|_{L^{p}(\bbR^{n},\VEC)}
        \right)_{j \in \bbN_{0}}
    \right\|_{\ell^{q}(\bbN_{0},\bbR)}
    < 
    \infty.
\end{align}
The usual modification applies if $q = \infty$.
Here, we use the quasi-norm of the generalized Lebesgue spaces
\begin{gather*}
    \| u \|_{L^{p}(\bbR^{n},\VEC)}
    :=
    \left(
        \int_{\bbR^{n}} \vecnorm{u(x)}^{p} \,\dif x
    \right)^{1/p}
    < \infty
    ,
\end{gather*}
for any $0 < p < \infty$, generalizing the $L^{p}$-norm as defined in~\eqref{function-spaces-on-domains:lp-norm} for $1 \leq p < \infty$.
Analogously, we define the quasi-norm of the generalized sequence spaces
\begin{gather*}
    \| (a_{j})_{j \in \bbN_{0}} \|_{\ell^{q}(\bbN_{0},\bbR)}
    :=
    \left(
        \sum_{j=0}^{\infty} |a_{j}|^{q}
    \right)^{1/q}
    < \infty
    ,
\end{gather*}
for any $0 < q < \infty$, which generalizes the usual $\ell^{q}$-norm for $1 \leq q < \infty$.

The Besov space $B^{s}_{p,q}(\Domain,\VEC)$ on any domain is defined as the space of restrictions of tempered distributions in $B^{s}_{p,q}(\bbR^{n},\VEC)$ to the domain $\Domain$. 
The space $B^{s}_{p,q}(\Domain,\VEC)$ is equipped with the quotient quasi-norm
\begin{align}\label{math:besov-quasinorm-on-domain}
    \| u \|_{B^{s}_{p,q}(\Domain,\VEC)}
    :=
    \inf\left\{
        \| \widetilde{u} \|_{B^{s}_{p,q}(\bbR^{n},\VEC)}
        \suchthat*
        \widetilde{u} \in B^{s}_{p,q}(\bbR^{n},\VEC), \; \widetilde{u} \restriction_{\Domain} = u
    \right\}
    .
\end{align}
The Triebel--Lizorkin spaces are defined similarly.
Let $s \in \bbR$, $0 < p < \infty$, and $0 < q \leq \infty$.
The Triebel--Lizorkin space $F^{s}_{p,q}(\bbR^{n},\VEC)$ is the quasi-Banach space of tempered distributions $u \in \calS'(\bbR^{n},\VEC)$ such that 
\begin{align}\label{math:triebel-lizorkin-quasinorm}
    \| u \|_{F^{s}_{p,q}(\bbR^{n},\VEC)}
    :=
    \left\|     
        \left(
            2^{js} \Delta_{j} u
        \right)_{j \in \bbN_{0}}
    \right\|_{L^{p}(\bbR^{n},\ell^{q}(\bbN_{0},\VEC))}
    <
    \infty.
\end{align}
Again, the usual modification applies if $q = \infty$.

The Triebel--Lizorkin space $F^{s}_{p,q}(\Domain,\VEC)$ is defined as the space of restrictions of tempered distributions in $F^{s}_{p,q}(\bbR^{n},\VEC)$ to the domain $\Domain$. 
The space $F^{s}_{p,q}(\Domain,\VEC)$ is equipped with the quotient quasi-norm
\begin{align}\label{math:triebel-lizorkin-quasinorm-on-domain}
    \| u \|_{F^{s}_{p,q}(\Domain,\VEC)}
    :=
    \inf\left\{     
        \| \widetilde{u} \|_{F^{s}_{p,q}(\bbR^{n},\VEC)}
        \suchthat*
        \widetilde{u} \in F^{s}_{p,q}(\bbR^{n},\VEC), \; \widetilde{u} \restriction_{\Domain} = u
    \right\}.
\end{align}

\begin{remark}\label{remark:nature-of-besov-triebel-lizorkin}
    The Besov and Triebel--Lizorkin spaces are quasi-Banach spaces.
    If $p,q \geq 1$, then they are Banach spaces and the quasi-norms~\eqref{math:besov-quasinorm} and~\eqref{math:triebel-lizorkin-quasinorm} are norms.
    However, if $p < 1$ or $q < 1$, then the spaces are quasi-Banach spaces that are not locally convex\footnote{Any locally convex quasi-Banach space has an equivalent norm that makes it a Banach space.}, but they continuously embed into the locally convex space of tempered distributions.
    
    Different choices of the smooth partition of unity $( \zeta_j )_{j \in \bbN_{0}}$ yield different Littlewood--Paley decompositions, but the different resulting quasi-norms on the Besov and Triebel--Lizorkin spaces are equivalent, both on $\bbR^{n}$ and on domains.
\end{remark}

\begin{remark}\label{remark:bessel-sobolev-besov-triebel-lizorkin}
    We summarize a few identifications of function spaces and refer to the literature~\cite[Section~2.1.2]{runstsickel} for further details.
    The Besov and Triebel--Lizorkin spaces generalize the Sobolev--Slobodeckij spaces in the sense that
    \begin{align*}
        \Wres^{s,p}(\bbR^{n},\VEC) &= F^s_{p,2}(\bbR^{n},\VEC),
        \qquad 
        s \in \bbZ, \quad 1 < p < \infty,
        \\
        \Wres^{s,p}(\bbR^{n},\VEC) &= B^s_{p,p}(\bbR^{n},\VEC) = F^s_{p,p}(\bbR^{n},\VEC),
        \qquad 
        s \in \bbR \setminus \bbZ, \quad 1 < p < \infty,
    \end{align*}
    with equivalent norms. 
    Moreover, for every $s \in \bbR$,
    \begin{align*}
        \Wint^{s,2}(\bbR^{n},\VEC)
        = B^s_{2,2}(\bbR^{n},\VEC)
        = F^s_{2,2}(\bbR^{n},\VEC).
    \end{align*}
    By restriction, the corresponding equivalences hold on every domain with $\Wres^{s,p}(\Domain,\VEC)$ in place of $\Wint^{s,p}(\bbR^{n},\VEC)$.

    Write $\nu(x) := \sqrt{ 1 + |x|^2 }$.
    For $s \in \bbR$ and $1 < p < \infty$, the Bessel potential space is defined by
    \begin{gather*}
        H^{s,p}(\bbR^{n},\VEC)
        := \left\{
            u \in \calS'(\bbR^{n},\VEC)
            \suchthat
            \Fourier^{-1}\left( \nu^s \Fourier u\right)
            \in L^p(\bbR^{n},\VEC)
        \right\}.
    \end{gather*}
    Equipped with the $L^p$-norm of $\Fourier^{-1}(\nu^s\Fourier u)$, it coincides with $F^s_{p,2}(\bbR^{n},\VEC)$ with equivalent norms.
    Another important special case is that $B^{s}_{\infty,\infty}(\bbR^{n},\VEC)$ coincides with the H\"older--Zygmund space $\calC^{s}(\bbR^{n},\VEC)$ for $s > 0$.
\end{remark}

\section{Universal extension operators on Lipschitz domains}\label{section:extension-operators-on-lipschitz-domains}

We introduce Lipschitz domains and review another central concept of this manuscript: universal extension operators on Sobolev spaces on bounded Lipschitz domains. 

Any domain $\Domain \subseteq \bbR^{n}$ is called a \emph{Lipschitz domain} if each $x \in \partial\Domain$ admits an open neighborhood $U \subseteq \bbR^{n}$ of $x$,
numbers $r > 0$ and $h > 0$, 
an isometry 
\begin{gather}\label{math:lipschitz-chart}
    \varphi : U \to (-r,r)^{n-1} \times (-h,h),
\end{gather}
and a Lipschitz-continuous function $\gamma : [-r,r]^{n-1} \to \bbR$ such that 
\begin{subequations}
\label{math:stronglylipschitz}
\begin{align}
    \label{math:stronglylipschitz:origin}
    \varphi(x) &= 0,
    \\
    \label{math:stronglylipschitz:condition1}
    \varphi\left( \Domain \cap U \right)
    &=
    \left\{\; y = (y',y_{n}) \suchthat y \in (-r,r)^{n-1} \times (-h,h), \; y_{n} < \gamma(y') \;\right\}, 
    \\
    \label{math:stronglylipschitz:condition2}
    \varphi\left( \partial\Domain \cap U \right)
    &=
    \left\{\; y = (y',y_{n}) \suchthat y \in (-r,r)^{n-1} \times (-h,h), \; y_{n} = \gamma(y') \;\right\}, 
    \\
    \label{math:stronglylipschitz:condition3}
    \varphi\left( {\closure\Domain}^{c} \cap U \right)
    &=
    \left\{\; y = (y',y_{n}) \suchthat y \in (-r,r)^{n-1} \times (-h,h), \; y_{n} > \gamma(y') \;\right\}. 
\end{align}
\end{subequations}
We call the mapping $\varphi : U \to (-r,r)^{n-1} \times (-h,h)$ in~\eqref{math:lipschitz-chart} a Lipschitz-graph coordinate chart,
and we call $U$ a \emph{Lipschitz-graph chart domain}. 

\begin{remark}
    One can show that each $x \in \partial\Domain$ admits a Lipschitz-graph chart domain $U \subseteq \bbR^{n}$, $x \in U$, 
    such that the intersection $U \cap \Domain$ is itself a Lipschitz domain. 
\end{remark}

The class of bounded Lipschitz domains features the equivalence of different definitions of function spaces.
For example, the intrinsic definition of Sobolev--Slobodeckij spaces on bounded Lipschitz domains,
using weak derivatives and Gagliardo seminorms, 
is equivalent to the definition of Sobolev--Slobodeckij spaces in terms of restrictions from $\bbR^{n}$. 
This result is standard, being an easy consequence of the following extension theorem, 
and we briefly review it for completeness.

\begin{proposition}\label{proposition:equivalence-of-sobolev-definitions}
    Let $\Domain \subseteq \bbR^{n}$ be a bounded Lipschitz domain.
    Let $s \geq 0$ and $1 < p < \infty$. 
    Then there exists a bounded linear operator 
    \begin{align*}
        \Extension_{\Domain,s,p} : \Wint^{s,p}(\Domain,\VEC) \to \Wint^{s,p}(\bbR^{n},\VEC)
    \end{align*}
    such that $\Extension_{\Domain,s,p} u \restriction_{\Domain} = u$ for all $u \in \Wint^{s,p}(\Domain,\VEC)$.
    
    In particular, the Banach spaces $\Wint^{s,p}(\Domain,\VEC)$ and $\Wres^{s,p}(\Domain,\VEC)$ coincide, with equivalent norms.
\end{proposition}
\begin{proof}
    For $s = 0$, extension by zero has the required properties.
    For $s > 0$, we refer to the literature; see~\cite[Theorem~1.4.3.1]{grisvard2011elliptic},~\cite{ern2021finite}.
    The equivalence of the intrinsic and restriction-based definitions of Sobolev--Slobodeckij spaces on bounded Lipschitz domains is now evident.
\end{proof}

A universal extension operator for Sobolev spaces (of integer order) on bounded Lipschitz domains is due to Rogers~\cite{rogers2006degree}. 
One easily sees that the extension operator is continuous between the corresponding spaces of smooth functions equipped with their Fr\'echet topologies.

\begin{proposition}\label{proposition:rogers-operator}
    Let $\Domain \subseteq \bbR^{n}$ be a bounded Lipschitz domain.
    Then there exists a bounded linear operator 
    \begin{align}\label{math:rogers-operator:sobolev}
        \Extension_{\Domain,W} : \Wint^{m,p}(\Domain,\VEC) \to \Wint^{m,p}(\bbR^{n},\VEC),
        \qquad 
        m \in \bbN_{0}, \quad 1 \leq p \leq \infty,
    \end{align}
    such that $\Extension_{\Domain,W} u \restriction_{\Domain} = u$ for all $u \in \Wint^{m,p}(\Domain,\VEC)$.

    In particular, we have a continuous linear operator 
    \begin{align}\label{math:rogers-operator:smooth}
        \Extension_{\Domain,W} : C^{\infty}_{\bounded}(\Domain,\VEC) \to C^{\infty}_{\bounded}(\bbR^{n},\VEC),
    \end{align}
    and the locally convex vector spaces $C^{\infty}_{\bounded}(\Domain,\VEC)$ and $C^{\infty}(\closure\Domain,\VEC)$ are isomorphic.
\end{proposition}
\begin{proof}
    We refer to the literature for an extension operator that satisfies~\eqref{math:rogers-operator:sobolev} for all $m \in \bbN_{0}$ and $1 \leq p \leq \infty$; see~\cite[Theorem~8]{rogers2006degree} if $n \geq 2$ and~\cite{seeley1964extension} if $n = 1$. %
    We briefly show that it satisfies the continuity condition~\eqref{math:rogers-operator:smooth}.
    
    Suppose that $u \in C^{\infty}_{\bounded}(\Domain,\VEC)$.
    Then $u \in \Wint^{m,1}(\Domain,\VEC)$ for all $m \in \bbN_{0}$, with
    \begin{align*}
        \| u \|_{\Wint^{m,1}(\Domain,\VEC)} 
        \leq 
        C \sup_{ x \in \Domain } \sup_{ \substack{ \alpha \in A(n) \\ |\alpha| \leq m } } \vecnorm{ \partial^{\alpha} u(x) }
        .
    \end{align*}
    We write $\widetilde{u} := \Extension_{\Domain,W} u$ and observe that $\widetilde{u} \in \Wint^{m,1}(\bbR^{n},\VEC)$ for all $m \in \bbN_{0}$,
    with
    \begin{align*}
        \| \widetilde{u} \|_{\Wint^{m,1}(\bbR^{n},\VEC)} 
        \leq 
        C \| u \|_{\Wint^{m,1}(\Domain,\VEC)} 
        .
    \end{align*}
    By the Sobolev embedding theorem, we have $\widetilde{u} \in C^{\infty}_{\bounded}(\bbR^{n},\VEC)$,
    i.e., $\widetilde{u}$ is a smooth function with bounded derivatives of all orders. We have 
    \begin{align*}
        \sup_{ x \in \bbR^{n} } \sup_{ \substack{ \alpha \in A(n) \\ |\alpha| \leq m } } \vecnorm{ \partial^{\alpha} \widetilde{u}(x) } 
        \leq 
        C \| \widetilde{u} \|_{\Wint^{m+n,1}(\bbR^{n},\VEC)} 
        .
    \end{align*}
    This also shows~\eqref{math:rogers-operator:smooth}, which means that the operator $\Extension_{\Domain,W} : C^{\infty}_{\bounded}(\Domain,\VEC) \to C^{\infty}_{\bounded}(\bbR^{n},\VEC)$ is continuous with respect to the topology of uniform convergence of all derivatives.

    Recall that $u \in C^{\infty}(\closure\Domain,\VEC)$ implies that there exists $\widetilde{u} \in C^{\infty}(\bbR^{n},\VEC)$ with $\widetilde{u} \restriction_{\Domain} = u$, and thus $u \in C^{\infty}_{\bounded}(\Domain,\VEC)$. 
    Hence, there exists a continuous embedding $\imath : C^{\infty}(\closure\Domain,\VEC) \to C^{\infty}_{\bounded}(\Domain,\VEC)$ of locally convex vector spaces.
    Applying $\Extension_{\Domain,W} : C^{\infty}_{\bounded}(\Domain,\VEC) \to C^{\infty}(\bbR^{n},\VEC)$ and taking the restriction to $\Domain$ gives another continuous embedding 
    $\jmath : C^{\infty}_{\bounded}(\Domain,\VEC) \to C^{\infty}(\closure\Domain,\VEC)$ that is inverse to the first one.
    We conclude $C^{\infty}_{\bounded}(\Domain,\VEC)$ and $C^{\infty}(\closure\Domain,\VEC)$ are isomorphic.
\end{proof}

We now reference the universal extension theorem for vector-valued distributions on bounded Lipschitz domains.
This provides an extension operator that is simultaneously continuous across the Besov and Triebel--Lizorkin scales of function spaces. This includes its continuity across the restriction-based Sobolev--Slobodeckij scale.

\begin{theorem}\label{theorem:extension-componentwise}
    Let $\Domain \subseteq \bbR^{n}$ be a bounded Lipschitz domain and let $\VEC$ be a finite-dimensional normed space.
    There exists a linear map 
    \begin{align*}
        \Extension_{\Domain} : \calS'(\Domain,\VEC) \to \calS'(\bbR^{n},\VEC),
    \end{align*}
    such that for every $u \in \calS'(\Domain,\VEC)$ we have 
    \begin{gather}\label{theorem:extension-componentwise:extension-property}
        \Extension_{\Domain} u \restriction_{\Domain} = u
        ,
    \end{gather}
    and the following operators are continuous:
    \begin{subequations}
    \begin{align}
        &
        \Extension_{\Domain} : C^{\infty}_{\bounded}(\Domain,\VEC) \to C^{\infty}_{\bounded}(\bbR^{n},\VEC), \label{math:extension-componentwise:smooth}
        \\
        &
        \Extension_{\Domain} : \Wres^{s,p}(\Domain,\VEC) \to \Wres^{s,p}(\bbR^{n},\VEC), \qquad s \in \bbR, \quad p \in (1,\infty), \label{math:extension-componentwise:sobolev}
        \\
        &
        \Extension_{\Domain} : B^{s}_{p,q}(\Domain,\VEC) \to B^{s}_{p,q}(\bbR^{n},\VEC), \qquad s \in \bbR, \quad p,q \in (0,\infty], \label{math:extension-componentwise:besov}
        \\
        &
        \Extension_{\Domain} : F^{s}_{p,q}(\Domain,\VEC) \to F^{s}_{p,q}(\bbR^{n},\VEC), \qquad s \in \bbR, \quad p \in (0,\infty), \quad q \in (0,\infty]. \label{math:extension-componentwise:triebellizorkin}
    \end{align}
    \end{subequations}
\end{theorem}
\begin{proof}
    Rychkov's extension operator~\cite{rychkov1999restrictions} is applied componentwise. 
    This yields a continuous operator $\Extension_{\Domain} : \calS'(\Domain,\VEC) \to \calS'(\bbR^{n},\VEC)$
    that satisfies the extension property~\eqref{theorem:extension-componentwise:extension-property}
    and the continuity conditions 
   ~\eqref{math:extension-componentwise:sobolev},
   ~\eqref{math:extension-componentwise:besov}, and
   ~\eqref{math:extension-componentwise:triebellizorkin}.
    Since~\eqref{math:extension-componentwise:smooth} is not explicitly stated in~\cite{rychkov1999restrictions}, we describe a short proof.
    Combining~\eqref{math:extension-componentwise:sobolev} with Proposition~\ref{proposition:equivalence-of-sobolev-definitions}, we obtain bounded operators
    \begin{align*}
        \Extension_{\Domain} : \Wint^{s,p}(\Domain,\VEC) \to \Wint^{s,p}(\bbR^{n},\VEC)
    \end{align*}
    for all $s \geq 0$ and $p \in (1,\infty)$.
    Adapting the argument in the proof of Proposition~\ref{proposition:rogers-operator}, using $p=2$, now yields the continuity condition~\eqref{math:extension-componentwise:smooth}.
\end{proof}

\begin{remark}
    Rychkov's universal extension operator~\cite{rychkov1999restrictions} covers Besov spaces and Triebel--Lizorkin spaces. 
    Up to endpoints of the integrability exponents, this includes Sobolev--Slobodeckij spaces. 
    If one is only interested in Sobolev spaces (with non-negative integer smoothness indices), 
    then the extension operator by Rogers is one of many alternatives~\cite{rogers2006degree}. 
\end{remark}

\section{Differential forms}\label{section:differential-forms}

\emph{Exterior calculus} reformulates the familiar three-dimensional vector calculus and generalizes it to arbitrary spatial dimension. 
The differential operators of vector calculus, which are the gradient, the curl, and the divergence, can be conceived as special instances of the exterior derivative. 
Exterior calculus is a subset of tensor calculus, and its objects are particular classes of tensor fields, namely alternating (skew-symmetric) tensor fields. 
These are fields over a domain that take values in skew-symmetric multilinear forms. 
Such fields are better known as differential forms, and exterior calculus also comes under the name of \emph{calculus of differential forms}.

We refer to the standard literature~\cite{cartan1967formes,amann2009analysisiii,agricola2010vektoranalysis,lee2012smooth,edwards2014advanced} for the necessary background on multilinear algebra and tensor fields.

\subsection{Exterior algebra}

Let $\Alt{k}(\bbR^{n})$ be the vector space of alternating $k$-linear forms on $\bbR^{n}$. 
For example, $\Alt{1}(\bbR^{n})$ is the dual space of $\bbR^{n}$, and $\Alt{0}(\bbR^{n}) = \bbR$. 
We follow the common convention to define $\Alt{k}(\bbR^{n}) := 0$ unless $0 \leq k \leq n$.

Given two alternating forms $u_1 \in \Alt{k}(\bbR^{n})$ and $u_2 \in \Alt{l}(\bbR^{n})$,
we write $u_1 \wedge u_2 \in \Alt{k+l}(\bbR^{n})$ for their \emph{exterior product}.
The exterior product is bilinear, and it is graded-commutative and associative:
\begin{gather*}
    u_1 \wedge u_2 = (-1)^{kl} u_2 \wedge u_1, 
    \qquad 
    u_1 \in \Alt{k}(\bbR^{n}), \quad u_2 \in \Alt{l}(\bbR^{n})
    ,
    \\
    ( u_1 \wedge u_2 ) \wedge u_3 = u_1 \wedge ( u_2 \wedge u_3 ),
    \qquad 
    u_1 \in \Alt{k}(\bbR^{n}), \quad u_2 \in \Alt{l}(\bbR^{n}), \quad u_3 \in \Alt{m}(\bbR^{n})
    .
\end{gather*}
Given a vector $a \in \bbR^{n}$ and an alternating $k$-linear form $u \in \Alt{k}(\bbR^{n})$ with $k > 0$,
we write $a \contract u$ for their \emph{interior product}, also known as \emph{contraction}, which is defined by 
\begin{gather*}
    (a \contract u)( v_1,\dots,v_{k-1}) := u( a, v_1,\dots,v_{k-1} ), \qquad v_1,\dots,v_{k-1} \in \bbR^{n}.
\end{gather*}
We define the interior product as zero if $k \leq 0$. 
The interior product is anti-commutative:
\begin{gather*}
    a_1 \contract ( a_2 \contract u )
    =
    - a_2 \contract ( a_1 \contract u )
    ,
    \qquad 
    a_1, a_2 \in \bbR^{n}, \quad u \in \Alt{k}(\bbR^{n})
    .
\end{gather*}
Recall that the standard basis of $\bbR^{n}$ gives rise to a basis of the dual space $\Alt{1}(\bbR^{n})$.
That basis is written $\cartanx^{1}, \dots, \cartanx^{n} \in \Alt{1}(\bbR^{n})$. 
We have a corresponding basis of $\Alt{k}(\bbR^{n})$ given by 
\begin{gather*}
    \left\{ 
        \cartanx^{i_1} \wedge \cartanx^{i_2} \wedge \dots \wedge \cartanx^{i_k} 
        \suchthat 
        1 \leq i_1 < i_2 < \dots < i_k \leq n 
    \right\}.
\end{gather*}
By definition, this set is empty unless $0 \leq k \leq n$.
Also, in the case $k=0$, we have $\Alt{0}(\bbR^{n}) = \bbR$ and the basis consists only of the constant function $1$. 
In particular, any $u \in \Alt{k}(\bbR^{n})$ has the form 
\begin{gather}\label{math:standard-representation}
    u 
    = 
    \sum_{1 \leq i_1 < \dots < i_k \leq n} 
    u_{i_1,\dots,i_k} \cartanx^{i_1} \wedge \cartanx^{i_2} \wedge \dots \wedge \cartanx^{i_k},
    \qquad 
    u_{i_1,\dots,i_k} \in \bbR.
\end{gather}
Here, the coefficients satisfy 
\begin{gather*}
    u_{i_1,\dots,i_k} = u( e_{i_1}, \dots, e_{i_k} ).
\end{gather*}
This definition includes $\cartanx^{\emptyset} := 1$ and $u_{\emptyset} := u$.
We define the \emph{Euclidean inner product} on $\Alt{k}(\bbR^{n})$ as 
\begin{gather*}
    \langle u, v \rangle 
    := 
    \sum_{1 \leq i_1 < \dots < i_k \leq n}
    u_{i_1,\dots,i_k} v_{i_1,\dots,i_k}
    .
\end{gather*}
The corresponding \emph{Euclidean norm} $\vecnorm{\cdot}$ is characterized by
\begin{gather*}
    \vecnorm{u}^{2}
    := 
    \sum_{1 \leq i_1 < \dots < i_k \leq n}
    u_{i_1,\dots,i_k}^{2}
    .
\end{gather*}
In what follows, we may use the \emph{musical isomorphisms}, 
which are the mutually inverse isometries 
\begin{gather*}
    \sharp : \Alt{1}(\bbR^{n}) \to \bbR^{n}, \quad \tsum_{i=1}^{n} v_i \cartanx^{i} \mapsto \tsum_{i=1}^{n} v_i e_{i},
    \\
    \flat  : \bbR^{n} \to \Alt{1}(\bbR^{n}), \quad \tsum_{i=1}^{n} a_i e_{i}        \mapsto \tsum_{i=1}^{n} a_i \cartanx^{i}.
\end{gather*}
The exterior product and the interior product satisfy a duality relationship with respect to the Euclidean inner product 
\begin{gather*}
    \langle a \contract u_1, u_2 \rangle = \langle u_1, a^{\flat} \wedge u_2 \rangle,
    \qquad 
    a \in \bbR^{n}, \quad u_1 \in \Alt{k+1}(\bbR^{n}), \quad u_2 \in \Alt{k}(\bbR^{n}).
\end{gather*}
Furthermore, the following inequalities are known:
\begin{align*}
    \vecnorm{ a \contract u }  &\leq \vecnorm{a} \cdot \vecnorm{u}, \qquad a \in \bbR^{n}, \quad u \in \Alt{k}(\bbR^{n}),
    \\
    \vecnorm{ u_1 \wedge u_2 } &\leq C_{k,l} \vecnorm{u_1} \cdot \vecnorm{u_2}, \qquad u_1 \in \Alt{k}(\bbR^{n}), \quad u_2 \in \Alt{l}(\bbR^{n}),
\end{align*}
where $C_{k,l} > 0$ is a numerical constant~\cite[Section~1.7.5]{federer}, which is known to satisfy
\begin{gather*}
    C_{k,l} \leq \sqrt{ \tbinom{k+l}{k} },
    \quad
    C_{k,n-k} = 1,
    \quad
    C_{k,0} = 1.
\end{gather*}

\subsection{Smooth differential forms}

A differential $k$-form on an $n$-dimensional open set is a function with values in the space of alternating $k$-linear forms $\Alt{k} = \Alt{k}(\bbR^{n})$.
The spaces of smooth differential forms 
\begin{gather*}
    C^{\infty}(\Domain,\Alt{k}),
    \qquad 
    C^{\infty}_{\bounded}(\Domain,\Alt{k}),
    \qquad 
    C^{\infty}(\closure\Domain,\Alt{k}),
    \qquad 
    C^{\infty}_{c}(\Domain,\Alt{k})
\end{gather*}
are of particular interest. The \emph{exterior derivative}
\begin{align}\label{math:exterior-derivative}
	\cartan^{k}: C^{\infty}(\Domain,\Alt{k}) \to C^{\infty}(\Domain,\Alt{k+1})
\end{align}
is defined as usual.
In explicit form, the exterior derivative of any $u \in C^{\infty}(\Domain,\Alt{k})$ reads 
\begin{equation}\label{math:exterior-wedge-formula}
    \cartan^{k} u
    := 
    \sum_{1 \leq j \leq n} \cartanx^{j} \wedge \partial_{j} u 
    =
    \sum_{\substack{ 1 \leq j \leq n \\ 1 \leq i_1 < \dots < i_k \leq n }}
    \partial_j u_{i_1,\dots,i_k} \cartanx^{j} \wedge \cartanx^{i_1} \wedge \dots \wedge \cartanx^{i_k}
    .
\end{equation}
It satisfies 
\begin{gather*}
    \cartan^{k+1} \cartan^{k} u = 0, \qquad u \in C^{\infty}(\Domain,\Alt{k}),
    \\
    \cartan^{k+l} ( u \wedge v ) = \cartan^{k} u \wedge v + (-1)^{k} u \wedge \cartan^{l} v, 
    \qquad 
    u \in C^{\infty}(\Domain,\Alt{k}), \quad v \in C^{\infty}(\Domain,\Alt{l}).
\end{gather*}

\subsection{Rough differential forms}

We define spaces of differential forms whose coefficients are in spaces of distributions and in non-smooth function spaces. 
The spaces $\calD'(\Domain,\Alt{k})$ and $\calS'(\Domain,\Alt{k})$ are preserved by the exterior derivative:
\begin{align*}
    \cartan^{k} : \calD'(\Domain,\Alt{k}) \to \calD'(\Domain,\Alt{k+1}), \qquad \cartan^{k} : \calS'(\Domain,\Alt{k}) \to \calS'(\Domain,\Alt{k+1}).
\end{align*}
For any $s \in \bbR$ and $p \in (1,\infty)$, 
we consider the spaces $L^{p}(\Domain,\Alt{k})$ and $\Wres^{s,p}(\Domain,\Alt{k})$.
These are Banach spaces. 
They can be characterized equivalently as those differential $k$-forms whose coefficients show the respective regularity properties. 
On the Sobolev scale, the exterior derivative lowers smoothness by one:
\begin{align*}
    \cartan^{k} : \Wres^{s,p}(\Domain,\Alt{k}) \to \Wres^{s-1,p}(\Domain,\Alt{k+1}), \qquad s \in \bbR, \quad p \in (1,\infty)
    .
\end{align*}
In addition to the Sobolev spaces, we also consider the Besov spaces $B^{s}_{p,q}(\Domain,\Alt{k})$ and the Triebel--Lizorkin spaces $F^{s}_{p,q}(\Domain,\Alt{k})$.
The exterior derivative also lowers smoothness by one on these scales:
\begin{align*}
    &
    \cartan^{k} : B^{s}_{p,q}(\Domain,\Alt{k}) \to B^{s-1}_{p,q}(\Domain,\Alt{k+1}), \qquad s \in \bbR, \quad p,q \in (0,\infty]
    ,
    \\
    &
    \cartan^{k} : F^{s}_{p,q}(\Domain,\Alt{k}) \to F^{s-1}_{p,q}(\Domain,\Alt{k+1}), \qquad s \in \bbR, \quad p \in (0,\infty), \quad q \in (0,\infty]
    .
\end{align*}

\section{Universal regular decomposition}\label{section:regular-decomposition}

The Poincar\'e homotopy operator is instrumental in proving the exactness in positive degrees of the smooth de~Rham complex on contractible domains. 
This homotopy operator defines a differential form via a line integral starting from a fixed pivot point and therefore assumes strong regularity.
If we take a smooth average over the pivot point, then the construction also works as a mapping between much rougher spaces. 
The averaged homotopy operators are pseudodifferential operators that increase the smoothness on the scales of Sobolev--Slobodeckij spaces, Besov spaces, and Triebel--Lizorkin spaces. 
The homotopy operators are interpreted as potential operators, since they produce a preimage for every exact form,
and the homotopy identities describe regular decompositions. 

The classical and averaged homotopy operators (or \emph{potential operators}) work on domains that are star-shaped with respect to a point or an open set, respectively. 
We can construct potential operators that satisfy a homotopy formula up to a defect, subject to the non-trivial topology of the domain. 
This results in a \emph{universal potential operator} and a \emph{universal regularized decomposition}.

\subsection{Poincar\'e homotopy operators}

The classical Poincar\'e homotopy operator with pivot point $a \in \bbR^{n}$ is defined for every $k \in \{1,\dots,n\}$ and any $u \in C^{\infty}(\bbR^{n},\Alt{k})$ by 
\begin{align*}
    \Poincare^{k}_{a} u(x)
    &:=
    \int_{0}^{1} t^{k-1} (x-a)\contract u\left(a+t(x-a)\right)\,\dif t
    ,
    \qquad
    x \in \bbR^{n}.
\end{align*}
We also formally define $\Poincare^{k}_{a} := 0$ for any $k \notin \{ 1, \dots, n \}$. %

This operator preserves smooth differential forms and satisfies a homotopy formula.

\begin{theorem}\label{theorem:poincare-homotopy}
    Let $a \in \bbR^{n}$. For every $k \in \{1,\dots,n\}$, we have a linear map 
    \begin{align*}
        \Poincare^{k}_{a} : C^{\infty}(\bbR^{n},\Alt{k}) \to C^{\infty}(\bbR^{n},\Alt{k-1})
        .
    \end{align*}
    For $u \in C^{\infty}(\bbR^{n},\Alt{k})$ with $k \in \bbZ$ and $x \in \bbR^{n}$, we have 
    \begin{subequations}\label{math:poincare-homotopy}
        \begin{alignat}{4}
            u(x) &= \Poincare^{1}_{a} \cartan^{0} u(x) + u(a),                           \qquad && k = 0,
            \\
            u(x) &= \cartan^{k-1} \Poincare^{k}_{a} u(x) + \Poincare^{k+1}_{a} \cartan^{k} u(x), \qquad && 1 \leq k \leq n-1,
            \\
            u(x) &= \cartan^{n-1} \Poincare^{n}_{a} u(x),                                  \qquad && k = n.
        \end{alignat}
    \end{subequations}
    We have the continuity condition 
    \begin{align}\label{math:poincare-homotopy-is-bounded}
        \sup_{ \substack{ \alpha \in A(n) \\ |\alpha| \leq m } }
        \vecnorm{ \partial^{\alpha} \Poincare^{k}_{a} u(x) }
        \leq 
        \sup_{ \substack{ 0 \leq t \leq 1 } }
        \sup_{ \substack{ \alpha \in A(n) \\ |\alpha| \leq m } }
        \left( \lvert x - a \rvert + m \right) 
        \vecnorm{ \partial^{\alpha} u( t x + (1-t) a ) }
        ,
        \qquad 
        x \in \bbR^{n}
        ,
        \quad 
        m \in \bbN_{0}
        .
    \end{align}
\end{theorem}
\begin{proof}
    Suppose that $u \in C^{\infty}(\bbR^{n},\Alt{k})$.
    Suppose that $k \geq 1$ and that there exists a ball $\smallball$ with $a, x \in \smallball$.
    Then a short induction proof and the dominated convergence theorem imply 
    that for all sequences of $m > 0$ indices $i_{1},\dots,i_{m} \in \{ 1, \dots, n \}$ we have 
    \begin{gather*}
        \partial_{i_{1}} \cdots \partial_{i_{m}}
        \int_{0}^{1} t^{k-1} u\left(a+t(x-a)\right)\,\dif t
        = 
        \int_{0}^{1} t^{k-1} t^{m}
        \partial_{i_{1}} \cdots \partial_{i_{m}}
        u\left(a+t(x-a)\right)\,\dif t
        .
    \end{gather*}
    Contraction of a smooth vector field with a smooth differential form preserves smoothness, 
    and therefore $\Poincare^{k}_{a} u \in C^{\infty}(\bbR^{n},\Alt{k-1})$. 
    Explicitly, 
    \begin{align*}
        \partial_{i_{1}} \cdots \partial_{i_{m}}
        \Poincare^{k}_{a} u(x)
        &
        =
        (x-a) 
        \contract 
        \int_{0}^{1} t^{k-1} t^{m}
        \partial_{i_{1}} \cdots \partial_{i_{m}}
        u\left(a+t(x-a)\right)\,\dif t
        \\&\quad 
        +
        \sum_{1 \leq j \leq m}
        e_{i_{j}}
        \contract 
        \int_{0}^{1} t^{k-1} t^{m-1}
        \partial_{i_{1}} \cdots \partial_{i_{j-1}} \partial_{i_{j+1}} \cdots \partial_{i_{m}}
        u\left(a+t(x-a)\right)\,\dif t
        .
    \end{align*}
    The inequality~\eqref{math:poincare-homotopy-is-bounded} is an easy consequence. 
    
    We prove the homotopy formula. 
    In the special case $k = 0$, this is the standard equation 
    \begin{align*}
        \Poincare^{1}_{a} \cartan^{0} u(x)
        = 
        \int_{0}^{1} \frac{\partial}{\partial t} u( a + t ( x - a ) ) \,\dif t
        =
        u(x) - u(a)
        .
    \end{align*}
    Suppose now $k \geq 1$. 
    On the one hand, the exterior derivative of the homotopy operator equals 
    \begin{align*}
        &
        \cartan^{k-1} \Poincare^{k}_{a} u(x)
        =
        \cartan^{k-1}\left( \int_{0}^{1} t^{k-1} (x-a)\contract u\left(a+t(x-a)\right)\,\dif t \right)
        \\&\quad 
        =
        \sum_{1 \leq j \leq n}
        \cartanx^{j} \wedge 
        \int_{0}^{1}
        t^{k-1} e_{j} \contract u\left(a+t(x-a)\right) \,\dif t
        +
        \sum_{1 \leq j \leq n}
        \cartanx^{j} \wedge 
        \int_{0}^{1}
        t^{k} 
        (x-a) \contract         
        \frac{ \partial u }{\partial x_{j}}\left(a+t(x-a)\right) 
        \,\dif t
        .
    \end{align*}
    On the other hand, 
    applying the Poincar\'e potential to the exterior derivative yields
    \begin{align*}
        &
        \Poincare^{k+1}_{a} \cartan^{k} u(x)
        =
        (x-a)\contract
        \sum_{1 \leq j \leq n}
        \int_{0}^{1} t^{k} 
        \cartanx^{j} \wedge \frac{ \partial u }{\partial x_{j}}\left(a+t(x-a)\right) \,\dif t
        \\&\quad 
        =
        \sum_{1 \leq j \leq n}
        (x-a)_{j}
        \int_{0}^{1} t^{k} \frac{ \partial u }{\partial x_{j}}\left(a+t(x-a)\right) \,\dif t
        -
        \sum_{1 \leq j \leq n}
        \cartanx^{j} \wedge
        \int_{0}^{1} t^{k} \,(x-a) \contract
        \frac{ \partial u }{\partial x_{j}}\left(a+t(x-a)\right) \,\dif t
        .
    \end{align*}
    We sum these terms and use cancellation effects:
    \begin{align*}
        &
        \cartan^{k-1} \Poincare^{k}_{a} u(x)
        +
        \Poincare^{k+1}_{a} \cartan^{k} u(x)
        \\&
        =
        \sum_{1 \leq j \leq n}
        \cartanx^{j} \wedge 
        \int_{0}^{1}
        t^{k-1} e_{j} \contract u\left(a+t(x-a)\right) \,\dif t
        +
        \sum_{1 \leq j \leq n}
        (x-a)_{j}
        \int_{0}^{1} t^{k} \frac{ \partial u }{\partial x_{j}}\left(a+t(x-a)\right) \,\dif t
        \\&
        =
        \int_{0}^{1} k \, t^{k-1} u\left(a+t(x-a)\right) \,\dif t
        +
        \sum_{1 \leq j \leq n}
        (x-a)_{j}
        \int_{0}^{1} t^{k} \frac{ \partial u }{\partial x_{j}}\left(a+t(x-a)\right) \,\dif t
        \\&
        =
        \int_{0}^{1} 
        \frac{\partial}{\partial t} \left( t^{k} u\left(a+t(x-a)\right) \right) \,\dif t
        .
    \end{align*}
    Here, we have used $\tsum_{j} \cartanx^{j} \wedge (e_{j} \contract u) = k \cdot u$ for $k$-forms. 
    The fundamental theorem of calculus now applies. 
    We conclude that the identities~\eqref{math:poincare-homotopy} hold, which completes the proof. 
\end{proof}

\subsection{Averaged local potential operators}

Suppose that $\weight \in C^{\infty}_{c}(\bbR^{n},\bbR)$ is non-negative, supported in a ball $\smallball$, and normalized to have unit integral:
\begin{align*}
    \int_{\bbR^{n}} \weight(x) \;\dif x = 1.
\end{align*}
We define the averaged homotopy operators
\begin{align}\label{math:averaged-poincare-operator-definition}
    \Poincare^{k}_{\weight} u(x)
    :=
    \int_{\bbR^{n}} \weight(a) \Poincare^{k}_{a} u(x) \,\dif a
    ,
    \qquad 
    k \in \bbZ
    .
\end{align}
This is a linear map 
\begin{align*}
    \Poincare^{k}_{\weight} : C^{\infty}(\bbR^{n},\Alt{k}) \to C^{\infty}(\bbR^{n},\Alt{k-1})
    .
\end{align*}
From the definition, it is clear that $\Poincare^{k}_{\weight} u(x)$ depends only on the values of $u$ over the convex hull of $\smallball \cup \{x\}$. 

We reformulate the averaged homotopy operator. Suppose that $1 \leq k \leq n$. By the substitution $\tau := 1/(1-t)$, one finds
\begin{align*}  
    \Poincare^{k}_{\weight} u(x)
    &
    =
    \int_{\bbR^{n}} 
    \int_{0}^{1} 
    t^{k-1}
    \weight(a) \,(x-a)\contract u\left( a + t(x-a) \right)
    \,\dif t
    \,\dif a
    \\&
    =
    \int_{\bbR^{n}} 
    \int_{1}^{\infty} 
    (\tau-1)^{k-1} \tau^{-k-1}
    \weight(a) \,(x-a)\contract u\left( x + (a-x)/\tau \right)
    \,\dif \tau
    \,\dif a
    ,
\end{align*} 
and by the substitution
$y(a) := x + (a - x)/\tau$, %
one arrives at the formulation 
\begin{align}\label{math:poincare-operator-reformulated}
    \Poincare^{k}_{\weight} u(x) 
    &= 
    \int_{\bbR^{n}} 
    \int_{1}^{\infty} (\tau-1)^{k-1} \tau^{n-k} \weight\left(x+\tau(y-x)\right)\,\dif \tau
    \,(x-y)\contract u(y)
    \,\dif y
    .
\end{align}
We define the kernel function $\calG_{\weight,k} : \bbR^{n} \times \bbR^{n} \to [0,\infty]$ by
\begin{equation}\label{math:mainkernel}
    \calG_{\weight,k}(\xi,\eta) := \int_{1}^{\infty} (\tau-1)^{n-k} \tau^{k-1} \weight\left(\eta+\tau(\xi-\eta)\right)\,\dif \tau, \qquad 1 \leq k \leq n
    .
\end{equation}
The kernel function enters the definition of the \emph{Poincar\'e potential} of any differential form $u \in C^{\infty}(\bbR^{n},\Alt{k})$:
\begin{align}\label{math:poincare-operator}
    \Poincare^{k}_{\weight} u(x) &= \int_{\bbR^{n}} \calG_{\weight,n-k+1}(y,x) \,(x-y)\contract u(y)\,\dif y, \qquad 1 \leq k \leq n.
\end{align}
In particular, the averaged homotopy operators are integral kernel operators.

We define a localized version of the potential operator on domains that are star-shaped with respect to a ball.
For any $x, y \in \bbR^{n}$, the closed line segment between $x$ and $y$ is written 
\begin{align*}
    [x,y] := \left\{ (1-t)\,x + t\,y \suchthat 0 \leq t \leq 1 \right\}
    .
\end{align*}
Suppose that $\Domain \subseteq \bbR^{n}$ is an open set. 
We say that $\Domain$ is \emph{star-shaped} with respect to a point $x \in \Domain$ 
if for all $y \in \Domain$ it holds that $[x,y] \subseteq \Domain$.
We say that $\Domain$ is \emph{star-shaped with respect to a ball} $B \subseteq \Domain$
if for all $x \in B$ and $y \in \Domain$ we have $[x,y] \subseteq \Domain$.

Domains star-shaped with respect to a ball are Lipschitz domains; see~\cite[Lemma~3.15]{lewicka2023calculus} and~\cite[Lemma~2]{melzak1959class}.

\begin{remark}\label{remark:star-shaped}
    Each star-shaped open set is connected and contractible, and thus it is simply connected. 
    Domains star-shaped with respect to a point need not be Lipschitz domains, but domains star-shaped with respect to a ball are always Lipschitz domains. 
    Cartan's homotopy relies on the set being star-shaped with respect to a point, 
    but the regularized homotopy operator requires the domain to be star-shaped with respect to a ball. 
\end{remark}

The essential properties of the Poincar\'e operator are summarized in the following statement.

\begin{theorem}\label{theorem:averaged-poincare-operator-on-domains}
    Let $\Domain \subseteq \bbR^{n}$ be bounded and star-shaped with respect to a ball $\smallball \Subset \Domain$.
    Suppose that $\weight \in C^{\infty}_{c}(\bbR^{n},\bbR)$ is non-negative, has support in $\smallball$, and has unit integral.
    We have continuous linear operators
    \begin{align}\label{math:averaged-poincare-operator-on-domains:distributions}
        \Poincare^{k}_{\weight} : \calD'(\Domain,\Alt{k}) \to \calD'(\Domain,\Alt{k-1}), \qquad k \in \bbZ, 
    \end{align}
    that restrict to the following continuous operators:
    \begin{subequations}\label{math:averaged-poincare-operator-on-domains:continuity}
    \begin{align}
        &\Poincare^{k}_{\weight} : \calS'(\Domain,\Alt{k}) \to \calS'(\Domain,\Alt{k-1}), \label{math:averaged-poincare-operator-on-domains:tempered-distributions}
        \\
        &\Poincare^{k}_{\weight} : C^{\infty}(\Domain,\Alt{k}) \to C^{\infty}(\Domain,\Alt{k-1}), \label{math:averaged-poincare-operator-on-domains:smooth}
        \\
        &\Poincare^{k}_{\weight} : C^{\infty}_{\bounded}(\Domain,\Alt{k}) \to C^{\infty}_{\bounded}(\Domain,\Alt{k-1}), \label{math:averaged-poincare-operator-on-domains:bounded-smooth}
        \\
        &\Poincare^{k}_{\weight} : \Wres^{s,p}(\Domain,\Alt{k}) \to \Wres^{s+1,p}(\Domain,\Alt{k-1}), \qquad s \in \bbR, \quad p \in (1,\infty), \label{math:averaged-poincare-operator-on-domains:sobolev}
        \\
        &\Poincare^{k}_{\weight} : B^{s}_{p,q}(\Domain,\Alt{k}) \to B^{s+1}_{p,q}(\Domain,\Alt{k-1}), \qquad s \in \bbR, \quad p,q \in (0,\infty], \label{math:averaged-poincare-operator-on-domains:besov}
        \\
        &\Poincare^{k}_{\weight} : F^{s}_{p,q}(\Domain,\Alt{k}) \to F^{s+1}_{p,q}(\Domain,\Alt{k-1}), \qquad s \in \bbR, \quad p \in (0,\infty), \quad q \in (0,\infty]. \label{math:averaged-poincare-operator-on-domains:triebellizorkin}
    \end{align}
    \end{subequations}
    In addition, if $u \in \calD'(\Domain,\Alt{k})$, then 
    \begin{subequations}\label{math:averaged-poincare-operator-on-domains:homotopy}
        \begin{alignat}{4}
            u &= \Poincare^{1}_{\weight} \cartan^{0} u + \langle u, \weight \rangle,         \qquad && k = 0,
            \\
            u &= \cartan^{k-1} \Poincare^{k}_{\weight} u + \Poincare^{k+1}_{\weight} \cartan^{k} u,  \qquad && 1 \leq k \leq n-1,
            \\
            u &= \cartan^{k-1} \Poincare^{k}_{\weight} u,                                        \qquad && k = n.
        \end{alignat}
    \end{subequations}
    If $u \in \calD'(\Domain,\Alt{k})$ and $\subdomain \subseteq \closure\Domain$ is relatively open with $u \restriction_{\subdomain \cap \Domain} \in C^{\infty}_{\bounded}(\subdomain \cap \Domain,\Alt{k})$, then for every $x_{0} \in \subdomain$
    there exists $\subdomain' \subseteq \closure\Domain$ relatively open with $x_{0} \in \subdomain'$
    and $( \Poincare^{k}_{\weight} u ) \restriction_{\subdomain' \cap \Domain} \in C^{\infty}_{\bounded}(\subdomain' \cap \Domain,\Alt{k-1})$.
    
\end{theorem}

\begin{proof}
    For any compact $K \subset \bbR^n$, write 
    \begin{align*}
        H(K) &:= \{ a + t(x-a) \suchthat a \in \supp(\weight), \ x \in K, \ t \in [0,1] \}
        .
    \end{align*}
    The set $H(K)$ is compact. Since $\Domain$ is star-shaped with respect to $\smallball$, it holds that 
    \begin{align}\label{math:convex-hull-in-domain}
        K \Subset \Domain \implies H(K) \Subset \Domain
        .
    \end{align}
    When $K \subseteq \bbR^{n}$ is compact, then $\calD_{K}(\bbR^{n},\Alt{k}) \subseteq \calD(\bbR^{n},\Alt{k})$ denotes the subspace of forms whose support is contained in $K$.

    It is clear that we have linear maps 
    \begin{align}
        \Poincare^{k}_{\weight} : C^{\infty}(\Domain,\Alt{k}) \to C^{\infty}(\Domain,\Alt{k-1})
        .
    \end{align}
    The bound~\eqref{math:poincare-homotopy-is-bounded} for the classical homotopy operator yields the following bound for the averaged homotopy operator:
    \begin{align}\label{math:averaged-poincare-homotopy-is-bounded}
        \max_{ \substack{ \alpha \in A(n) \\ |\alpha| \leq m } }
        \vecnorm{ \partial^{\alpha} \Poincare^{k}_{\weight} u(x) }
        \leq 
        \left( \diameter{\Domain} + m \right) 
        \sup_{ z \in H(\{x\}) }
        \max_{ \substack{ \alpha \in A(n) \\ |\alpha| \leq m } }
        \vecnorm{ \partial^{\alpha} u( z ) }
        ,
        \qquad 
        x \in \Domain
        ,
        \quad 
        m \in \bbN_{0}
        .
    \end{align}
    As an immediate consequence, we have continuous operators~\eqref{math:averaged-poincare-operator-on-domains:smooth} and~\eqref{math:averaged-poincare-operator-on-domains:bounded-smooth}, the first one with respect to the topology of uniform convergence of all derivatives on compact subsets of $\Domain$,
    the second with respect to the topology of uniform convergence of all derivatives on $\Domain$.
    
    The same formula defines linear maps 
    \begin{align}
        \Poincare^{k}_{\weight} : C^{\infty}(\bbR^{n},\Alt{k}) \to C^{\infty}(\bbR^{n},\Alt{k-1})
    \end{align}
    continuous with respect to the topology of uniform convergence of all derivatives on compact sets. 
    This is compatible in the sense that 
    \begin{align*}
        ( \Poincare^{k}_{\weight} u ) \restriction_{\Domain}
        =
        \Poincare^{k}_{\weight} ( u \restriction_{\Domain} )
        ,
        \qquad 
        u \in C^{\infty}(\bbR^{n},\Alt{k})
        .
    \end{align*}

    The discussion of distributional spaces requires a recourse to pseudodifferential calculus. 
    The averaged Poincar\'e operator, say, as a continuous linear map 
    $\Poincare^{k}_{\weight} : C^{\infty}_{c}(\bbR^{n},\Alt{k}) \subseteq C^{\infty}(\bbR^{n},\Alt{k}) \to C^{\infty}(\bbR^{n},\Alt{k-1})$,
    is a pseudodifferential operator on $\bbR^{n}$
    whose symbol belongs to the local H\"ormander symbol class $S^{-1}_{1,0}(\bbR^{n})$; see~\cite[Theorem~3.2]{costabel2010}.
    This means that its symbol in Fourier space is a smooth function $a_{k} \in C^{\infty}( \bbR^{n} \times \bbR^{n}, \Alt{k-1} \otimes \Alt{k} )$, taking values in the $\bbC$-linear mappings from $\Alt{k}$ into $\Alt{k-1}$,
    such that for any compact set $K \subseteq \bbR^{n}$ and multi-indices $\alpha, \beta \in A(n)$
    there exists $C_{K,\alpha,\beta} \geq 0$ satisfying  
    \begin{align*}
        &
        \forall (x,\xi) \in K \times \bbR^{n} 
        :
        \vecnorm{ \partial_{x}^{\alpha} \partial_{\xi}^{\beta} a_{k}( x, \xi ) }
        \leq 
        C_{K,\alpha,\beta} 
        \left( 1 + \lvert \xi \rvert \right)^{ - 1 - |\beta| }
        ,
        \\
        &
        \forall (x,\xi) \in K \times \bbR^{n} 
        :
        a_{k}( x, \xi ) = \overline{ a_{k}( x, - \xi ) }
        .
    \end{align*}
    Thus, using the Fourier transform and its inverse, we can write 
    \begin{align*}
        ( \Poincare^{k}_{\weight} u )(x)
        &
        = 
        \int_{\bbR^{n}} \int_{\bbR^{n}} e^{2\pi i (x-y) \cdot \xi} a_{k}(x,\xi) u(y) \,\dif y \,\dif \xi
        = 
        \int_{\bbR^{n}} e^{2\pi i x \cdot \xi} a_{k}(x,\xi) (\Fourier u)(\xi) \,\dif \xi   
    \end{align*}
    for any $u \in C^{\infty}_{c}(\bbR^{n},\Alt{k})$ and $x \in \bbR^{n}$. 
    Here, we only observe that the integral is convergent 
    because $\Fourier u \in \calS(\bbR^{n},\Alt{k}_{\bbC})$. 
    We refer to the literature (\cite[Chapter XVIII]{hoermander},~\cite[Definition~7.1]{grubb},~\cite[7.1.4]{taylor2},~\cite[Definition 2.1, Definition 2.3]{treves1},~\cite{wong}) for the definition of pseudodifferential operators and their symbols. 
    
    We let $\varrho \in C^{\infty}_{c}(\bbR^{n},\bbR)$ with $\varrho = 1$ within an open bounded neighborhood of $\closure\Domain$.
    For simplicity, we assume that $\supp(\varrho)$ is convex and thus star-shaped with respect to $\smallball$. 

    Define $\widetilde{a}_{k}(x,\xi) := \varrho(x) a_{k}(x,\xi)$. 
    Now, $\widetilde{a}_{k} \in C^{\infty}( \bbR^{n} \times \bbR^{n}, \Alt{k-1} \otimes \Alt{k} )$,
    and for all multi-indices $\alpha, \beta \in A(n)$ we have the estimate 
    \begin{align*}
        \forall
        (x,\xi) \in \bbR^{n} \times \bbR^{n}
        :
        \vecnorm{ \partial_{x}^{\alpha} \partial_{\xi}^{\beta} \widetilde{a}_{k}( x, \xi ) }
        \leq 
        C_{\alpha,\beta} 
        \left( 1 + \lvert \xi \rvert \right)^{ - 1 - |\beta| }
        .
    \end{align*}
    Thus, $\widetilde{a}_{k}$ belongs to the global H\"ormander symbol class $S^{-1}_{1,0}(\bbR^{n})$;
    see~\cite[Definition 18.1.1]{hoermander},~\cite[Definition 6.1]{wong}.
    This immediately implies that this symbol defines a pseudodifferential operator ${\euQ}^{k}$,
    and this operator is a continuous map between Schwartz spaces~\cite[Theorem~18.1.6]{hoermander},~\cite[Proposition~6.7]{wong}:
    \begin{align}\label{math:improved-operator:schwartz}
        {\euQ}^{k} : \calS(\bbR^{n},\Alt{k}) \to \calS(\bbR^{n},\Alt{k-1})
        .
    \end{align}
    By definition, 
    \begin{align}\label{math:improved-operator-coincides-with-original}
        ( {\euQ}^{k} u ) \restriction_{\Domain} = ( \Poincare^{k}_{\weight} u ) \restriction_{\Domain}, \qquad u \in \calS(\bbR^{n},\Alt{k}).
    \end{align}
    The pseudodifferential operator ${\euQ}^{k}$ defines a transpose $({\euQ}^{k})^{t}$ with the same symbol class~\cite[Theorem~18.1.7]{hoermander},~\cite[Theorem I.3.1]{treves1},
    which is again a well-defined continuous map between Schwartz spaces:
    \begin{align}\label{math:improved-operator-transpose:schwartz}
        ({\euQ}^{k})^{t} : \calS(\bbR^{n},\Alt{k-1}) \to \calS(\bbR^{n},\Alt{k})
        .
    \end{align}
    The transpose pseudodifferential operator is also uniquely defined by the pairing 
    \begin{align*}
        \left\langle 
        {\euQ}^{k} u, v
        \right\rangle 
        :=
        \left\langle 
        u, ({\euQ}^{k})^{t} v
        \right\rangle 
        ,
        \qquad 
        u \in \calS(\bbR^{n},\Alt{k}), 
        \quad 
        v \in \calS(\bbR^{n},\Alt{k-1}). 
    \end{align*}
    By duality
    (see also~\cite[Theorem~18.1.7]{hoermander},~\cite[Lemma~7.1.1]{taylor2}),
    {}~\eqref{math:improved-operator:schwartz} embeds into a continuous operator
    \begin{align}\label{math:improved-operator:tempered-distribution}
        {\euQ}^{k} : \calS'(\bbR^{n},\Alt{k}) \to \calS'(\bbR^{n},\Alt{k-1})
        .
    \end{align}

    We show that $({\euQ}^{k})^{t}$ maps $\calD(\bbR^{n},\Alt{k-1})$ into $\calD(\bbR^{n},\Alt{k})$.
    Suppose that $v \in \calD(\bbR^{n},\Alt{k-1})$ and let $u \in \calD(\bbR^{n},\Alt{k})$ such that $\supp(u) \cap H(\supp(v)) = \emptyset$.
    Then ${\euQ}^{k} u$ vanishes on $\supp(v)$. From 
    \begin{align*}
        \left\langle 
        u, ({\euQ}^{k})^{t} v
        \right\rangle 
        =
        \left\langle 
        {\euQ}^{k} u, v
        \right\rangle 
        =
        0
    \end{align*}
    we conclude that $\supp( ({\euQ}^{k})^{t} v ) \subseteq H(\supp v)$. 
    As $({\euQ}^{k})^{t}$ is already a continuous mapping on the Schwartz space, 
    for any compact set $K \subseteq \bbR^{n}$ we therefore have a continuous operator 
    $({\euQ}^{k})^{t} : \calD_{K}(\bbR^{n},\Alt{k-1}) \to \calD_{H(K)}(\bbR^{n},\Alt{k})$. 
    In particular, we have a continuous operator 
    \begin{align}
        ({\euQ}^{k})^{t} : \calD(\bbR^{n},\Alt{k-1}) \to \calD(\bbR^{n},\Alt{k})
        .
    \end{align}
    By duality, we have a continuous operator 
    \begin{align}\label{math:improved-operator:distributions}
        {\euQ}^{k} : \calD'(\bbR^{n},\Alt{k}) \to \calD'(\bbR^{n},\Alt{k-1})
        .
    \end{align}
    Note that~\eqref{math:improved-operator:schwartz} embeds into~\eqref{math:improved-operator:distributions}.

    If $u \in \calS'(\bbR^{n},\Alt{k})$ with $u \restriction_{\Domain} = 0$, 
    then for any $v \in C^{\infty}_{c}(\Domain,\Alt{k-1})$ we observe 
    \begin{align*}
        \left\langle 
        {\euQ}^{k} u, v
        \right\rangle 
        =
        \left\langle 
        u, ({\euQ}^{k})^{t} v
        \right\rangle 
        =
        0
    \end{align*}
    because $\supp( ({\euQ}^{k})^{t} v ) \subseteq \Domain$. 
    Since~\eqref{math:improved-operator:tempered-distribution} is continuous and preserves the kernel of the restriction to $\Domain$,
    we conclude that we have a well-defined operator 
    \begin{align}\label{math:improved-operator-transpose:restricted-tempered-distribution}
        {\euQ}^{k} : \calS'(\Domain,\Alt{k}) \to \calS'(\Domain,\Alt{k-1})
    \end{align}
    that is continuous with respect to the quotient topology.

    Recall that for any compact set $K \subseteq \Domain$ we have $H(K) \subseteq \Domain$.
    We have seen that we have a continuous operator 
    $({\euQ}^{k})^{t} : \calD_{K}(\bbR^{n},\Alt{k-1}) \to \calD_{H(K)}(\bbR^{n},\Alt{k})$.
    Thus, we have a well-defined continuous operator 
    \begin{align}\label{math:improved-operator-transpose:domain-distributions}
        ({\euQ}^{k})^{t} : \calD(\Domain,\Alt{k-1}) \to \calD(\Domain,\Alt{k})
        .
    \end{align}
    We define a linear operator 
    \begin{align}\label{math:improved-operator:domain-distributions}
        {\euQ}^{k} : \calD'(\Domain,\Alt{k}) \to \calD'(\Domain,\Alt{k-1}), \qquad u \mapsto \langle u, ({\euQ}^{k})^{t} \cdot \rangle
        .
    \end{align}
    The continuity of~\eqref{math:improved-operator-transpose:domain-distributions} implies that~\eqref{math:improved-operator:domain-distributions}
    is continuous.
    Furthermore, the restriction to $\calS'(\Domain,\Alt{k})$ coincides with the original operator~\eqref{math:improved-operator-transpose:restricted-tempered-distribution}.

    We have already established the continuity of the operator $\Poincare^{k}_{\weight} : \calD(\bbR^{n},\Alt{k}) \to \calD'(\bbR^{n},\Alt{k-1})$. 
    By restriction, we have a continuous operator 
    \begin{align}\label{math:poincare-auxiliary}
        \Poincare^{k}_{\weight} : \calD(\Domain,\Alt{k}) \to \calD'(\Domain,\Alt{k-1})
        .
    \end{align}
    The latter agrees with ${\euQ}^{k} : \calD'(\Domain,\Alt{k}) \to \calD'(\Domain,\Alt{k-1})$ on $\calD(\Domain,\Alt{k})$, which is dense in $\calD'(\Domain,\Alt{k})$. Therefore, 
    the operator~\eqref{math:poincare-auxiliary} embeds into the continuous operator~\eqref{math:improved-operator:domain-distributions}.
    We conclude that we have continuous linear operators~\eqref{math:averaged-poincare-operator-on-domains:distributions} and~\eqref{math:averaged-poincare-operator-on-domains:tempered-distributions} that extend~\eqref{math:averaged-poincare-operator-on-domains:smooth}.

    We now draw on continuity properties of pseudodifferential operators with symbols in the global Hörmander class $S^{-1}_{1,0}$.
    Since ${\euQ^{k}}$ has a symbol in the global Hörmander class, we have continuous linear operators (see~\cite[Theorem~4.1]{kalton2007interpolation},~\cite[Theorem~7, Theorem~9]{marschall1996nonregular}):
    \begin{align*}
        \euQ^{k} &: B^{s}_{p,q}(\bbR^{n},\Alt{k}) \to B^{s+1}_{p,q}(\bbR^{n},\Alt{k-1}), \qquad s \in \bbR, \quad p, q \in (0,\infty]
        ,
        \\
        \euQ^{k} &: F^{s}_{p,q}(\bbR^{n},\Alt{k}) \to F^{s+1}_{p,q}(\bbR^{n},\Alt{k-1}), \qquad s \in \bbR, \quad p \in (0,\infty), \quad q \in (0,\infty]
        .
    \end{align*}
    These operators map distributions that vanish on $\Domain$ to distributions that vanish on $\Domain$,
    and therefore we have well-defined continuous operators
    \begin{align*}
        \euQ^{k} &: B^{s}_{p,q}(\Domain,\Alt{k}) \to B^{s+1}_{p,q}(\Domain,\Alt{k-1}), \qquad s \in \bbR, \quad p, q \in (0,\infty]
        ,
        \\
        \euQ^{k} &: F^{s}_{p,q}(\Domain,\Alt{k}) \to F^{s+1}_{p,q}(\Domain,\Alt{k-1}), \qquad s \in \bbR, \quad p \in (0,\infty), \quad q \in (0,\infty]
        ,
    \end{align*}
    with respect to the quotient topologies. 
    The domains of the latter two mappings embed into $\calS'(\Domain,\Alt{k})$, where $\euQ^{k}$ and $\Poincare^{k}_{\weight}$ agree. 
    As a consequence, we obtain the continuity of the operators~\eqref{math:averaged-poincare-operator-on-domains:besov} and~\eqref{math:averaged-poincare-operator-on-domains:triebellizorkin}.

    The continuity condition~\eqref{math:averaged-poincare-operator-on-domains:sobolev} between Sobolev--Slobodeckij spaces is now evident due to the usual identification with Besov and Triebel--Lizorkin spaces, depending on whether both $s$ and $s+1$ are either integral or non-integral; see Remark~\ref{remark:bessel-sobolev-besov-triebel-lizorkin}.

    The homotopy identities~\eqref{math:averaged-poincare-operator-on-domains:homotopy} are known to hold for $u \in \calD(\Domain,\Alt{k})$. 
    Smooth compactly supported differential forms on $\Domain$ are dense in the distributional differential forms on $\Domain$, 
    and the involved operations are continuous on these distribution spaces. 
    By density,~\eqref{math:averaged-poincare-operator-on-domains:homotopy} holds for any $u \in \calD'(\Domain,\Alt{k})$.

    It remains to show that $\Poincare^{k}_{\weight}$ preserves local smoothness.
    Here, we only need to consider $k \in \{1,\dots,n\}$. 
    Let $u \in \calD'(\Domain,\Alt{k})$ and $\subdomain \subseteq \closure\Domain$ be relatively open 
    such that $u \restriction_{\subdomain \cap \Domain} \in C^{\infty}_{\bounded}(\subdomain \cap \Domain,\Alt{k})$.
    Let $x_{0} \in \subdomain$. 
    If $x_{0} \in \Domain$, then we restrict $\subdomain$ to a Lipschitz domain $\subdomain \Subset \Domain$, 
    and we write $\widetilde\subdomain := \subdomain$. 
    If $x_{0} \in \partial\Domain$, then we may restrict $\subdomain$ such that there exists a Lipschitz-graph chart domain $\widetilde\subdomain$ at $x_{0}$ with $\subdomain = \closure\Domain \cap \widetilde\subdomain$;
    in this case, both $\widetilde\subdomain$ and $\subdomain \cap \Domain$ are Lipschitz domains.
    
    Since $\subdomain \cap \Domain$ is now assumed to be a Lipschitz domain, 
    there exists $\widetilde{u} \in C^{\infty}(\bbR^{n},\Alt{k})$ with $\widetilde{u} \restriction_{\subdomain \cap \Domain} = u \restriction_{\subdomain \cap \Domain}$. 
    
    There exists $\eta \in C^{\infty}_{c}(\bbR^{n},\bbR)$ with $\eta \equiv 1$ in an open neighborhood of $x_{0}$ and $\supp(\eta) \Subset \widetilde\subdomain$.
    Therefore, $\eta \widetilde{u} \in C^{\infty}_{c}(\widetilde\subdomain,\Alt{k})$, from which we conclude 
    \begin{align*}
        \euQ^{k}( \eta \widetilde{u} ) \in C^{\infty}(\bbR^{n},\Alt{k-1}).
    \end{align*}
    This implies 
    \begin{align*}
        \Poincare^{k}_{\weight}( \eta u ) = \euQ^{k}( \eta \widetilde{u} ) \restriction_{\Domain} \in C^{\infty}_{\bounded}(\Domain,\Alt{k-1})
        .
    \end{align*}
    The expression $\Poincare^{k}_{\weight}( (1-\eta) u )$ is analyzed via its Schwartz kernel. 
    Abbreviate $R := 2\diameter{\Domain}$. 
    There exists $\epsilon > 0$ small enough such that $\eta \equiv 1$ on $\Ball{x_0}{2\epsilon}$.
    For all $x \in \Ball{x_0}{\epsilon}$, 
    the properties of $\weight$ and the definition of the Poincar\'e operator yield 
    \begin{align*}
        \Poincare^{k}_{\weight}\left( (1-\eta) u \right)
        &
        =
        \int_{\Domain}
        \int_{1}^{\infty} 
        (\tau-1)^{k-1} \tau^{n-k} \weight\left( x + \tau( y - x ) \right) 
        ( 1 - \eta(y) ) ( x - y ) \contract u(y) 
        \,d\tau
        \,dy
        \\&
        =
        \iint\limits_{\substack{ y \in \Ball{x_0}{2R} \setminus \Ball{x_0}{2\epsilon} \\ 1 < \tau < \diameter{\Domain} \epsilon^{-1} }}
        (\tau-1)^{k-1} \tau^{n-k} \weight\left( x + \tau( y - x ) \right) 
        ( 1 - \eta(y) )
        \,d\tau
        \ ( x - y ) \contract u(y) 
        \,dy
        .
    \end{align*}
    In particular, for $x$ near $x_{0}$, the Schwartz kernel is smooth in both variables and its support in the second variable is contained in a compact subset of $\Domain$.
    It follows that $\Poincare^{k}_{\weight}\left( (1-\eta) u \right)$ is smooth near $x_{0}$ for any $u \in \calD'(\Domain,\Alt{k})$.
    We conclude that 
    \begin{align*}
        \Poincare^{k}_{\weight}( u ) = \Poincare^{k}_{\weight}( \eta u ) + \Poincare^{k}_{\weight}( (1-\eta) u )
    \end{align*}
    is smooth near $x_{0}$ with bounded derivatives of all orders. 
    
    The proof is complete.
\end{proof}

\begin{remark}
    The averaged Poincar\'e operators have symbols only in the local Hörmander class, which is why the multiplication with a cut-off function is used.
    Costabel and McIntosh~\cite{costabel2010} focus on the Bessel potential scale with $p=q=2$ for the sake of exposition. 
    The proof in this manuscript, albeit perhaps verbose, collects elementary facts on pseudodifferential operators, explicitly available in the literature, 
    to make the key results reproducible for a wider audience.
\end{remark}

\subsection{Global potential operators}

We follow the construction of a global potential as described by Costabel and McIntosh~\cite{costabel2010}. 
We develop two geometric auxiliary results. 
The first one is that within sufficiently small local chart domains, the Lipschitz domain is star-shaped with respect to an interior ball. 
The second one is that each bounded Lipschitz domain admits a family of open covers with subordinate partitions of unity satisfying additional properties.

\begin{lemma}\label{lemma:lipschitz-graphs-have-star-domains}
    Let $\Domain \subseteq \bbR^{n}$ be a Lipschitz domain and $p \in \partial\Domain$. 
    Then every open neighborhood of $p$ contains a convex Lipschitz chart domain $U$ with $p \in U$
    such that 
    \begin{itemize}
        \item 
        $U \cap \Domain$ is star-shaped with respect to a ball $\smallball \subseteq U \cap \Domain$.
        \item
        $U \cap \bbR^{n} \setminus \closure\Domain$ is star-shaped with respect to a ball $\smallball' \subseteq U \cap \bbR^{n} \setminus \closure\Domain$.
    \end{itemize}
\end{lemma}
\begin{proof}
    The statement is trivial if $n=1$, so we consider the case $n \geq 2$. 
    By assumption, there exists an open neighborhood $U \subseteq \bbR^{n}$ of $p$, numbers $r,h > 0$, an isometry 
    \begin{gather*}
        \varphi : U \to (-r,r)^{n-1} \times (-h,h),
    \end{gather*}
    and a Lipschitz-continuous function $\gamma : (-r,r)^{n-1} \to \bbR$ such that~\eqref{math:stronglylipschitz} holds.
    We write $L > 0$ for some positive upper bound on the Lipschitz constant of $\gamma$ and abbreviate
    \begin{gather*}
        \subdomain := (-r,r)^{n-1} \times (-h,h)
        ,
        \qquad 
        \Gamma := \varphi( U \cap \partial\Domain )
        ,
        \qquad 
        \calC := \left\{ (y',y_n) \in \bbR^{n} \suchthat y_n < - L \|y'\| \right\}
        .
    \end{gather*}
    Note that $\Gamma$ is the graph of $\gamma$.
    Then $\calC$ is a cone whose apex is the origin. 
    Abbreviate $\delta := 2\sqrt{n-1} r$. 
    The graph of $\gamma$ is contained within $(-r,r)^{n-1} \times (-L\delta,L\delta)$. 
    Hence, 
    \begin{gather*}
        ( -r, r )^{n-1} \times ( -h, -3L\delta ) \subseteq \Gamma + \calC
        .
    \end{gather*}
    We restrict $r > 0$ so that $h/2 > 3L\delta$. 
    Then $\varphi( U \cap \Domain )$ is star-shaped with respect to 
    $(-r,r)^{n-1} \times ( -h, -6L\sqrt{n-1} r )$,
    which includes the desired ball $\varphi(\smallball)$. 
    In this situation, $\varphi( U \cap ( \bbR^{n} \setminus \closure\Domain ) )$ is analogously star-shaped with respect to $(-r,r)^{n-1} \times ( 6L\sqrt{n-1} r, h )$, which includes the desired ball $\varphi(\smallball')$.
    We can require both $r$ and $h$ to be as small as needed. 
\end{proof}

\begin{theorem}\label{theorem:cover-and-partition-of-unity}
    Let $\Domain \subseteq \bbR^{n}$ be a bounded Lipschitz domain. 
    There exist finitely many finite families
    \begin{align*}
        \Cover_{1}, \quad \Cover_{2}, \quad \dots, \quad \Cover_{N}
    \end{align*}
    of sets open in $\bbR^{n}$ having the form 
    \begin{align*}
        \Cover_{i} := \left\{ \subdomain_{i,j} \suchthat* 1 \leq j \leq m_i \right\}, \qquad 1 \leq i \leq N,
    \end{align*}
    with the following properties:
    \begin{itemize}
        \item 
        For each $1 \leq i \leq N$, the family $\Cover_{i}$ is an open cover of $\closure\Domain$ by convex sets.
        \item 
        For each $1 \leq i \leq N$ and each $1 \leq j \leq m_{i}$, 
        there exists an open ball $\smallball_{i,j} \subseteq \subdomain_{i,j} \cap \Domain$ such that $\subdomain_{i,j} \cap \Domain$ is star-shaped with respect to $\smallball_{i,j}$.
        \item 
        For each $1 \leq i \leq N$, there exists a smooth partition of unity $(\chi_{i,j})_{1 \leq j \leq m_i}$ satisfying
        \[
            \chi_{i,j} \in C^{\infty}_{c}(\subdomain_{i,j},\bbR),
            \qquad
            \chi_{i,j} \geq 0,
            \qquad
            \sum_{j=1}^{m_i} \chi_{i,j} = 1 \text{ on an open neighborhood of } \closure{\Domain}.
        \]
        \item 
        For every $x \in \closure\Domain$ there exists $1 \leq i \leq N$ such that
        \[
            \chi_{i,1}, \quad \chi_{i,2}, \quad \dots, \quad \chi_{i,m_{i}} \text{ are constant near } x.
        \]
    \end{itemize}
\end{theorem}

\begin{proof}
    We first show that each point $p \in \closure{\Domain}$ admits an arbitrarily small convex open neighborhood $U \subseteq \bbR^{n}$ 
    such that $U \cap \Domain$ is star-shaped with respect to a ball $\smallball \subset U \cap \Domain$.

    This is clear if $p \in \Domain$, because then we can take a sufficiently small open ball $U$ with $\closure{U} \subset \Domain$ and a smaller concentric open ball $\smallball \subseteq U$. 
    If $p \in \partial\Domain$, then we use Lemma~\ref{lemma:lipschitz-graphs-have-star-domains}.
    
    With that auxiliary result in place, we construct the open covers and the partitions of unity. 
    For each $z \in \closure\Domain$, we choose an open neighborhood $\subdomain_{z} \subseteq \bbR^{n}$ 
    such that both $\subdomain_{z}$ and $\subdomain_{z} \cap \Domain$ are star-shaped with respect to an open ball $\smallball_{z} \subseteq \subdomain_{z} \cap \Domain$.
    In addition, we choose $V_{z}$ open and convex such that $z \in V_{z} \Subset \subdomain_{z}$.
    
    Because $\closure\Domain$ is compact, 
    there exists a finite set $\calZ := \{ z_1, \dots, z_{N} \} \subseteq \closure\Domain$ 
    such that $\{ V_{z} \suchthat z \in \calZ \}$ is a finite open cover of $\closure\Domain$.
    
    Consider any $1 \leq i \leq N$. 
    Using the above auxiliary result and a compactness argument, 
    we obtain a finite family of open sets $\Cover_{i}$ such that
    \begin{itemize}
        \item $\closure{\Domain} \subset \bigcup\Cover_{i}$.
        \item For each $U \in \Cover_{i}$, both $U$ and $U \cap \Domain$ are star-shaped with respect to a ball $\smallball \subseteq U \cap \Domain$.
        \item $\subdomain_{z_{i}} \in \Cover_{i}$.
        \item $V_{z_{i}} \cap \bigcup \left( \Cover_{i} \setminus \{\subdomain_{z_i}\} \right) = \emptyset$.
    \end{itemize}
    Let $\left( \chi_{U} \right)_{U \in \Cover_{i}}$ be a smooth partition of unity subordinate to $\Cover_{i}$, 
    which means that each $\chi_{U} \in C^{\infty}_{c}(U,\bbR)$ is non-negative
    and there exists an open neighborhood $\widetilde{\Domain}$ of $\closure\Domain$ such that 
    \begin{align*}
        1 = \tsum_{U \in \Cover_{i}} \chi_{U} \restriction_{\widetilde{\Domain}}.
    \end{align*}
    In addition, only one of the functions is non-zero on $V_{z_{i}}$, and that function must equal $1$ on $V_{z_{i}} \cap \widetilde{\Domain}$ because we have a partition of unity.
    Since the sets $V_{z}$, $z \in \calZ$, cover $\closure{\Domain}$, this shows that for every $x \in \closure{\Domain}$ there exists $1 \leq i \leq N$ such that the functions $\chi_{i,1}, \dots, \chi_{i,m_{i}}$ are constant near $x$.
    This defines the families $\Cover_{1}, \dots, \Cover_{N}$ of open covers and the partitions of unity, 
    as in the statement of the theorem.
\end{proof}

\begin{remark}
    By taking the union of the open covers $\Cover_{1}, \dots, \Cover_{N}$ in Theorem~\ref{theorem:cover-and-partition-of-unity},
    we may assume that the different partitions of unity are all subordinate to a single open cover. 
\end{remark}

We proceed with the construction of the global potential operator.
Let $\Domain$ be a bounded Lipschitz domain. 
By Theorem~\ref{theorem:cover-and-partition-of-unity}, there exists a finite family of open covers 
\begin{align*}
    \Cover_{1} := (\subdomain_{1,j})_{1 \leq j \leq m_{1}}, 
    \quad
    \Cover_{2} := (\subdomain_{2,j})_{1 \leq j \leq m_{2}}, 
    \quad
    \dots,
    \quad
    \Cover_{N} := (\subdomain_{N,j})_{1 \leq j \leq m_{N}},
\end{align*}
and associated partitions of unity 
\begin{align*}
    ( \chi_{1,j})_{1 \leq j \leq m_{1}}, 
    \quad
    ( \chi_{2,j})_{1 \leq j \leq m_{2}}, 
    \quad
    \dots, 
    \quad
    ( \chi_{N,j})_{1 \leq j \leq m_{N}}.
\end{align*}
For each $1 \leq i \leq N$ and $1 \leq j \leq m_{i}$, the set $\subdomain_{i,j}$ is star-shaped with respect to a ball $\smallball_{i,j} \subseteq \subdomain_{i,j} \cap \Domain$. 
We assume that $\weight_{i,j} \in C^{\infty}_{c}(\bbR^{n},\bbR)$ is a non-negative smooth function with unit integral and support in $\smallball_{i,j}$.
This allows us to define the local Poincar\'e potential operators 
\begin{align*}
    \Poincare^{k}_{i,j} : \calD'(\subdomain_{i,j} \cap \Domain,\Alt{k}) \to \calD'(\subdomain_{i,j} \cap \Domain,\Alt{k-1})
    ,
    \qquad 
    k \in \bbZ
    .
\end{align*}
We define the intermediate potential operators 
\begin{align*}
    \Poincare^{k}_{i} : \calD'(\Domain,\Alt{k}) \to \calD'(\Domain,\Alt{k-1})
    ,
    \qquad 
    k \in \bbZ
\end{align*}
by the sum 
\begin{align*}
    \Poincare^{k}_{i}(u) := \sum_{1 \leq j \leq m_{i}} \chi_{i,j} \Poincare^{k}_{i,j}( u \restriction_{ \subdomain_{i,j}  \cap \Domain} ).
\end{align*}
We comment on why the operator $\Poincare^{k}_{i}$ is well-defined: 
the support properties of $\Poincare^{k}_{i,j}$ imply that $\Poincare^{k}_{i,j}( u \restriction_{ \subdomain_{i,j} \cap \Domain } )$ on the open set $\Domain \cap \subdomain_{i,j}$ depends only on $u \restriction_{ \subdomain_{i,j} \cap \Domain }$. 

Suppose that $u \in \calD'(\Domain,\Alt{k})$. 
In the case $k = 0$, we observe 
\begin{align*}
    \Poincare^{1}_{i} \cartan^{k} u
    &
    =
    \sum_{1 \leq j \leq m_{i}} \chi_{i,j} \Poincare^{1}_{i,j}( \cartan^{0} u \restriction_{ \subdomain_{i,j} \cap \Domain } )
    =
    \sum_{1 \leq j \leq m_{i}} \chi_{i,j} \left( u - \langle u, \weight_{i,j} \rangle \right)
    =
    u - \sum_{1 \leq j \leq m_{i}} \langle u, \weight_{i,j} \rangle \chi_{i,j} 
    .
\end{align*}
In the case $k \geq 1$, we observe the homotopy property
\begin{align*}
    \cartan^{k-1} \Poincare^{k}_{i} u + \Poincare^{k+1}_{i} \cartan^{k} u 
    &
    =
    \cartan^{k-1} \sum_{1 \leq j \leq m_{i}} ( \chi_{i,j} \Poincare^{k}_{i,j}(u) ) 
    + 
    \sum_{1 \leq j \leq m_{i}} \chi_{i,j} \Poincare^{k+1}_{i,j} \cartan^{k} u 
    \\&
    =
    \sum_{1 \leq j \leq m_{i}}
    \cartan^{k-1} ( \chi_{i,j} \Poincare^{k}_{i,j}(u) ) 
    - 
    \chi_{i,j} \cartan^{k-1} \Poincare^{k}_{i,j} u 
    + 
    \chi_{i,j} \cartan^{k-1} \Poincare^{k}_{i,j} u 
    + 
    \chi_{i,j} \Poincare^{k+1}_{i,j} \cartan^{k} u 
    \\&
    =
    u
    +
    \sum_{1 \leq j \leq m_{i}}
    \cartan^{k-1} ( \chi_{i,j} \Poincare^{k}_{i,j}(u) )
    - 
    \chi_{i,j} \cartan^{k-1} \Poincare^{k}_{i,j} u 
    .
\end{align*}
Note that 
\begin{align*}
    \sum_{1 \leq j \leq m_{i}}
    \cartan^{k-1} ( \chi_{i,j} \Poincare^{k}_{i,j}(u) )
    - 
    \chi_{i,j} \cartan^{k-1} \Poincare^{k}_{i,j} u 
    =
    \sum_{1 \leq j \leq m_{i}}
    ( \cartan^{0} \chi_{i,j} ) \wedge \Poincare^{k}_{i,j}(u) 
    .
\end{align*}
This motivates the following definitions of residual operators. 
We define 
\begin{align*}
    \Defect^{k}_{i} : \calD'(\Domain,\Alt{k}) \to \calD'(\Domain,\Alt{k})
    ,
    \qquad 
    k \in \bbZ
    ,
\end{align*}
by setting 
\begin{align*}
    \Defect^{0}_{i}(u)
    &
    :=
    \sum_{1 \leq j \leq m_{i}} 
    \langle u, \weight_{i,j} \rangle \chi_{i,j},
    \\
    \Defect^{k}_{i}(u)
    &
    :=
    -\sum_{1 \leq j \leq m_{i}}
    ( \cartan^{0} \chi_{i,j} ) \wedge \Poincare^{k}_{i,j}(u), \quad k \geq 1
    .
\end{align*}
The homotopy property takes the form 
\begin{align*}
    &
    u = \Poincare^{1}_{i} \cartan^{0} u + \Defect^{0}_{i}(u),
    \\
    &
    u = \cartan^{k-1} \Poincare^{k}_{i} u + \Poincare^{k+1}_{i} \cartan^{k} u + \Defect^{k}_{i}(u), \quad k \geq 1
    .
\end{align*}
This residual operator commutes with the exterior derivative:
in the case $k = 0$, we find
\begin{align*}
    \cartan^{0} \Defect^{0}_{i}(u)
    &
    =
    \sum_{1 \leq j \leq m_{i}} 
    \langle u, \weight_{i,j} \rangle \cartan^{0} \chi_{i,j} 
    =
    \sum_{1 \leq j \leq m_{i}} 
    ( u - \Poincare^{1}_{i,j} \cartan^{0} u ) \cartan^{0} \chi_{i,j} 
    \\&
    =
    -
    \sum_{1 \leq j \leq m_{i}} 
    ( \Poincare^{1}_{i,j} \cartan^{0} u ) \cartan^{0} \chi_{i,j}
    =
    -
    \sum_{1 \leq j \leq m_{i}} 
    \cartan^{0} \chi_{i,j} \wedge ( \Poincare^{1}_{i,j} \cartan^{0} u )
    =
    \Defect^{1}_{i}( \cartan^{0} u) 
    ,
\end{align*}
and in the case $k \geq 1$, we find that 
\begin{align*}
    \cartan^{k} \Defect^{k}_{i}( u )
    &
    =
    -\sum_{1 \leq j \leq m_{i}}
    \cartan^{k}\left( ( \cartan^{0} \chi_{i,j} ) \wedge \Poincare^{k}_{i,j}(u) \right)
    =
    \sum_{1 \leq j \leq m_{i}}
    ( \cartan^{0} \chi_{i,j} ) \wedge \cartan^{k-1} \Poincare^{k}_{i,j} u
    \\&
    =
    \sum_{1 \leq j \leq m_{i}}
    ( \cartan^{0} \chi_{i,j} ) \wedge u
    -
    \sum_{1 \leq j \leq m_{i}}
    ( \cartan^{0} \chi_{i,j} ) \wedge \Poincare^{k+1}_{i,j}( \cartan^{k} u)
    =
    \Defect^{k+1}_{i}( \cartan^{k} u)
    ,
\end{align*}
where we have used the local homotopy formula and $\sum_{j} \cartan^{0} \chi_{i,j} = 0$. 
We define the global Poincar\'e potential operator as the following sum of compositions:
\begin{align*}
    \Poincare^{k} 
    := 
    \sum_{i=1}^{N} \Defect^{k-1}_{1} \Defect^{k-1}_{2} \cdots \Defect^{k-1}_{i-1} \Poincare^{k}_{i}
    ,
    \qquad 
    k \geq 1
    . 
\end{align*}
Using that the residual operators commute with the exterior derivative, we again observe a homotopy property:
\begin{align*}
    \cartan^{k-1} \Poincare^{k} u + \Poincare^{k+1} \cartan^{k} u 
    &
    =
    \cartan^{k-1} \sum_{i=1}^{N} \Defect^{k-1}_{1} \Defect^{k-1}_{2} \cdots \Defect^{k-1}_{i-1} \Poincare^{k}_{i} u
    +
    \sum_{i=1}^{N} \Defect^{k}_{1} \Defect^{k}_{2} \cdots \Defect^{k}_{i-1} \Poincare^{k+1}_{i} \cartan^{k} u
    \\&
    =
    \sum_{i=1}^{N} \Defect^{k}_{1} \Defect^{k}_{2} \cdots \Defect^{k}_{i-1} \cartan^{k-1} \Poincare^{k}_{i} u
    +
    \sum_{i=1}^{N} \Defect^{k}_{1} \Defect^{k}_{2} \cdots \Defect^{k}_{i-1} \Poincare^{k+1}_{i} \cartan^{k} u
    \\&
    =
    \sum_{i=1}^{N} \Defect^{k}_{1} \Defect^{k}_{2} \cdots \Defect^{k}_{i-1} ( \cartan^{k-1} \Poincare^{k}_{i} + \Poincare^{k+1}_{i} \cartan^{k} ) u
    \\&
    =
    \sum_{i=1}^{N} \Defect^{k}_{1} \Defect^{k}_{2} \cdots \Defect^{k}_{i-1} ( u - \Defect^{k}_{i} u )
    \\&
    =
    u - \Defect^{k}_{1} \Defect^{k}_{2} \cdots \Defect^{k}_{N} u
    .
\end{align*}
We consider the composition 
\begin{align*}
    \Defect^{k} := \Defect^{k}_{1} \Defect^{k}_{2} \cdots \Defect^{k}_{N}
    ,
    \qquad 
    k \in \bbZ
    .
\end{align*}
We have shown that 
\begin{align*}
    u = \cartan^{k-1} \Poincare^{k} u + \Poincare^{k+1} \cartan^{k} u + \Defect^{k} u
    .
\end{align*}
It is clear that 
\begin{align*}
    \cartan^{k} \Defect^{k} = \Defect^{k+1} \cartan^{k}.
\end{align*}
As a consequence for closed forms, if $\cartan^{k} u = 0$, then
\begin{align*}
    u = \cartan^{k-1} \Poincare^{k} u + \Defect^{k} u
    .
\end{align*}
As a consequence for exact forms, if $u = \cartan^{k-1} v$, then
\begin{align*}
    u 
    = 
    \cartan^{k-1} \Poincare^{k} \cartan^{k-1} v + \Defect^{k} \cartan^{k-1} v
    = 
    \cartan^{k-1} \left( \Poincare^{k} u + \Defect^{k-1} v \right)
    .
\end{align*}

We review a few known mapping properties of these operators; see~\cite[Corollary~3.4, Theorem~4.6]{costabel2010}. 
It may come as a surprise that $\Defect^{k}$ is, in fact, infinitely smoothing.

\begin{theorem}\label{theorem:remainder-operator-is-smoothing}
    We have continuous linear operators 
    \begin{align}\label{math:defect-operator:distributions}
        \Defect^{k} : \calD'(\Domain,\Alt{k}) \to C^{\infty}(\closure\Domain,\Alt{k})
        ,
        \qquad 
        k \in \bbZ
        .
    \end{align}
\end{theorem}

\begin{proof} 
    We show that each $\Defect^{k}_{i}$ preserves local smoothness. 
    Suppose that $u \in \calD'(\Domain,\Alt{k})$ such that for some relatively open set $\subdomain \subseteq \closure\Domain$ 
    we have $u \restriction_{\subdomain \cap \Domain} \in C^{\infty}_{\bounded}(\subdomain \cap \Domain,\Alt{k})$.
    Given any $x_{0} \in \subdomain$, 
    we use Theorem~\ref{theorem:averaged-poincare-operator-on-domains} and finiteness of the open cover $\Cover_{i}$ 
    to conclude that there exists a set $\subdomain' \subseteq \closure\Domain$ relatively open with $x_{0} \in \subdomain'$
    such that 
    $( \Defect^{k}_{i} u ) \restriction_{\subdomain' \cap \Domain} \in C^{\infty}_{\bounded}(\subdomain' \cap \Domain,\Alt{k})$.

    For every $x_{0} \in \closure\Domain$ there exists $1 \leq i \leq N$ and $\subdomain \subseteq \closure\Domain$ relatively open 
    such that $x_{0} \in \subdomain$ and $\chi_{i,1}, \chi_{i,2}, \dots, \chi_{i,m_{i}}$ are constant on $\subdomain$. 
    It follows that $(\Defect^{k}_{i} u) \restriction_{\subdomain \cap \Domain}$ is smooth on $\subdomain \cap \Domain$ with bounded derivatives of all orders.

    Taking into account that $\closure\Domain$ is compact, 
    we conclude that the composition $\Defect^{k} = \Defect^{k}_{1} \Defect^{k}_{2} \cdots \Defect^{k}_{N}$ defines a linear map 
    \begin{align}\label{math:defect-operator-as-linear-map}
        \Defect^{k} : \calD'(\Domain,\Alt{k}) \to C^{\infty}_{\bounded}(\Domain,\Alt{k}).
    \end{align}

    It remains to prove that~\eqref{math:defect-operator-as-linear-map} is continuous. 
    As an auxiliary result, 
    we first notice that $\Defect^{k}_{i} : \calD'(\Domain,\Alt{k}) \to \calD'(\Domain,\Alt{k})$ are continuous operators,
    and so their compositions are continuous:
    \begin{align}\label{math:defect-between-distributions}
        \Defect^{k} : \calD'(\Domain,\Alt{k}) \to \calD'(\Domain,\Alt{k})
        .
    \end{align}
    Suppose now that $u_{i} \in \calD'(\Domain,\Alt{k})$ is a sequence that converges to $u \in \calD'(\Domain,\Alt{k})$
    such that the sequence $\Defect^{k} u_{i}$ converges to $v \in C^{\infty}(\Domain,\Alt{k})$,
    each with respect to the locally convex topology.
    Then $\Defect^{k} u_{i}$ converges to $v$ in the topology of $\calD'(\Domain,\Alt{k})$.
    By continuity of $\Defect^{k} : \calD'(\Domain,\Alt{k}) \to \calD'(\Domain,\Alt{k})$, we infer that $v = \Defect^{k} u$. 
    We conclude that $\Defect^{k} : \calD'(\Domain,\Alt{k}) \to C^{\infty}_{\bounded}(\Domain,\Alt{k})$ has a sequentially closed graph.
    
    As $\calD'(\Domain,\Alt{k})$ is ultrabornological and $C^{\infty}_{\bounded}(\Domain,\Alt{k})$ is a Fr\'echet space, hence a webbed space,
    we can apply a classical generalization of the closed graph theorem~\cite[35.2]{koethe2} for operators between such spaces with sequentially closed graphs 
    to infer the desired continuity of $\Defect^{k}$.
    On a Lipschitz domain, $C^{\infty}_{\bounded}(\Domain,\Alt{k}) = C^{\infty}(\closure\Domain,\Alt{k})$, and thus the proof is complete. %
\end{proof}

\begin{theorem}\label{theorem:continuity-of-poincare-and-remainder-operator}
    We have linear maps 
    \begin{align*}
        \Poincare^{k} : \calD'(\Domain,\Alt{k}) \to \calD'(\Domain,\Alt{k-1}),
        \qquad 
        \Defect^{k} : \calD'(\Domain,\Alt{k}) \to C^{\infty}(\closure\Domain,\Alt{k}). 
    \end{align*}
    We have continuous linear operators 
    \begin{subequations}\label{math:continuity-of-poincare-and-remainder-operator:poincare}
    \begin{align}
        &\Poincare^{k} : \calS'(\Domain,\Alt{k}) \to \calS'(\Domain,\Alt{k-1}), 
        \\
        &\Poincare^{k} : C^{\infty}(\Domain,\Alt{k}) \to C^{\infty}(\Domain,\Alt{k-1}), 
        \\
        &\Poincare^{k} : \Wres^{s,p}(\Domain,\Alt{k}) \to \Wres^{s+1,p}(\Domain,\Alt{k-1}), \qquad s \in \bbR, \quad p \in (1,\infty),
        \\
        &\Poincare^{k} : B^{s}_{p,q}(\Domain,\Alt{k}) \to B^{s+1}_{p,q}(\Domain,\Alt{k-1}), \qquad s \in \bbR, \quad p,q \in (0,\infty],
        \\
        &\Poincare^{k} : F^{s}_{p,q}(\Domain,\Alt{k}) \to F^{s+1}_{p,q}(\Domain,\Alt{k-1}), \qquad s \in \bbR, \quad p \in (0,\infty), \quad q \in (0,\infty].
    \end{align}
    \end{subequations}
    We have compact linear operators 
    \begin{align}\label{math:continuity-of-poincare-and-remainder-operator:defect-operator}
        \Defect^{k} : \Wres^{s,p}(\Domain,\Alt{k}) \to \Wres^{t,p}(\Domain,\Alt{k}), \qquad s,t \in \bbR, \quad p \in (1,\infty).
    \end{align}
    In addition, if $u \in \calD'(\Domain,\Alt{k})$, then 
    \begin{align}\label{math:continuity-of-poincare-and-remainder-operator:homotopy}
        u = \cartan^{k-1} \Poincare^{k} u + \Poincare^{k+1} \cartan^{k} u + \Defect^{k} u
        .
    \end{align}
\end{theorem}
\begin{proof}
    The properties~\eqref{math:continuity-of-poincare-and-remainder-operator:poincare} and~\eqref{math:continuity-of-poincare-and-remainder-operator:homotopy} follow from the previous results.
    The compactness of~\eqref{math:continuity-of-poincare-and-remainder-operator:defect-operator} follows with a Rellich--Kondrachov compactness theorem on bounded Lipschitz domains. 
\end{proof}

The potential operators $\Poincare^{k}$ and the residual operators $\Defect^{k}$ are universal 
in the sense that they are defined and continuous on a range of function spaces simultaneously.
We now introduce additional universal constructions: 
a universal harmonic complement and a universal regular decomposition.
These generalize the homotopy formulas.

\begin{theorem}\label{theorem:universal-regular-decomposition}
    For each $k \in \bbZ$, there exists a finite-dimensional subspace 
    \[
        \Harmonic^{k} \subset C^{\infty}(\closure\Domain,\Alt{k}), \qquad \cartan^{k} \Harmonic^{k} = 0,
    \]
    and continuous operators 
    \begin{align*}
        \widehat\Poincare^{k} : \calD'(\Domain,\Alt{k}) \to \calD'(\Domain,\Alt{k-1}),
        \qquad 
        \Proj^{k} : \calD'(\Domain,\Alt{k}) \to \Harmonic^{k}
    \end{align*}
    such that the following properties hold:
    \begin{itemize}
        \item 
        $\Proj^{k}$ is a projection. 
        \item 
        For any $u \in \calD'(\Domain,\Alt{k})$, we have 
        \begin{gather*}
            \cartan^{k} \widehat\Poincare^{k+1} \cartan^{k} u
            =
            \cartan^{k} u,
            \qquad 
            \cartan^{k} \Proj^{k} u
            =
            0
            ,
            \qquad
            \Proj^{k+1} \cartan^{k} u = 0
            ,
            \\
            u
            =
            \widehat\Poincare^{k+1} \cartan^{k} u
            + 
            \Proj^{k} u
            + 
            \cartan^{k-1} \widehat\Poincare^{k} u
            .
        \end{gather*}
        \item 
        If $u \in \calD'(\Domain,\Alt{k})$ with $\cartan^{k} u = 0$, then
        \begin{align*}
            u = \Proj^{k} u + \cartan^{k-1} \widehat\Poincare^{k} u.
        \end{align*}
        \item 
        We have 
        \begin{align*}
            \dim \Harmonic^{k}
            =
            \dim 
            \frac{
                \ker\left( \cartan^{k  } : C^{\infty}(\closure\Domain,\Alt{k}) \to C^{\infty}(\closure\Domain,\Alt{k+1}) \right)
            }{
                \rng\left( \cartan^{k-1} : C^{\infty}(\closure\Domain,\Alt{k-1}) \to C^{\infty}(\closure\Domain,\Alt{k}) \right)
            }
            .
        \end{align*}
        \item 
        The operators $\widehat\Poincare^{k}$ have the following continuous restrictions:
        \begin{subequations}\label{math:universal-regular-decomposition:continuity}
        \begin{align}
            &\widehat\Poincare^{k} : \calS'(\Domain,\Alt{k}) \to \calS'(\Domain,\Alt{k-1}), \label{math:universal-regular-decomposition:tempered-distributions}
            \\
            &\widehat\Poincare^{k} : C^{\infty}(\Domain,\Alt{k}) \to C^{\infty}(\Domain,\Alt{k-1}), \label{math:universal-regular-decomposition:smooth}
            \\
            &\widehat\Poincare^{k} : C^{\infty}_{\bounded}(\Domain,\Alt{k}) \to C^{\infty}_{\bounded}(\Domain,\Alt{k-1}), \label{math:universal-regular-decomposition:bounded-smooth}
            \\
            &\widehat\Poincare^{k} : \Wres^{s,p}(\Domain,\Alt{k}) \to \Wres^{s+1,p}(\Domain,\Alt{k-1}), \qquad s \in \bbR, \quad p \in (1,\infty), \label{math:universal-regular-decomposition:sobolev}
            \\
            &\widehat\Poincare^{k} : B^{s}_{p,q}(\Domain,\Alt{k}) \to B^{s+1}_{p,q}(\Domain,\Alt{k-1}), \qquad s \in \bbR, \quad p,q \in (0,\infty], \label{math:universal-regular-decomposition:besov}
            \\
            &\widehat\Poincare^{k} : F^{s}_{p,q}(\Domain,\Alt{k}) \to F^{s+1}_{p,q}(\Domain,\Alt{k-1}), \qquad s \in \bbR, \quad p \in (0,\infty), \quad q \in (0,\infty]. \label{math:universal-regular-decomposition:triebellizorkin}
        \end{align}
        \end{subequations}
        \item 
        We have direct sums 
        \begin{subequations}\label{math:universal-regular-decomposition:kernel}
            \begin{align}
                 \left\{ u \in \calD'(\Domain,\Alt{k}) \suchthat \cartan^{k} u = 0 \right\}
                 &
                 =
                 \Harmonic^{k} \oplus \cartan^{k-1} \calD'(\Domain,\Alt{k-1}),
                 \\
                 \left\{ u \in \Wres^{s,p}(\Domain,\Alt{k}) \suchthat \cartan^{k} u = 0 \right\}
                 &
                 =
                 \Harmonic^{k} \oplus \cartan^{k-1} \Wres^{s+1,p}(\Domain,\Alt{k-1}),
                 \quad 
                 p \in (1,\infty), \quad s \in \bbR,
                 \\
                 \left\{ u \in B^{s}_{p,q}(\Domain,\Alt{k}) \suchthat \cartan^{k} u = 0 \right\}
                 &
                 =
                 \Harmonic^{k} \oplus \cartan^{k-1} B^{s+1}_{p,q}(\Domain,\Alt{k-1}),
                 \quad 
                 s \in \bbR, \quad p,q \in (0,\infty],
                 \\
                 \left\{ u \in F^{s}_{p,q}(\Domain,\Alt{k}) \suchthat \cartan^{k} u = 0 \right\}
                 &
                 =
                 \Harmonic^{k} \oplus \cartan^{k-1} F^{s+1}_{p,q}(\Domain,\Alt{k-1}),
                 \quad 
                 s \in \bbR, \quad p \in (0,\infty), \quad q \in (0,\infty].
            \end{align}
        \end{subequations}
    \end{itemize}
\end{theorem}

\begin{proof}
    \newcommand{\Hilbertcomplement}{\euY}
    We use a differential complex of Hilbert spaces for auxiliary computations. Define 
    \begin{align*}
        X^{k} := \Wres^{n-k,2}(\Domain,\Alt{k}).
    \end{align*}
    Write $X^{k}_{\perp} \subseteq X^{k}$ for the orthogonal complement of the nullspace of $\cartan^{k} : X^{k} \to X^{k+1}$,
    and write $\Proj_{\perp}^{k} : X^{k} \to X^{k}_{\perp}$ for the orthogonal projection. 
    We have bounded operators $\cartan^{k} : X^{k} \to X^{k+1}$. We first prove that these differentials have closed range.
    Suppose to the contrary that $(v_{m})_{m \in \bbN}$ is a sequence of unit vectors in $X^{k}_{\perp}$ 
    such that $(\cartan^{k} v_{m})_{m \in \bbN}$ converges to zero. 
    Recall that 
    \begin{align*}
        v_{m} 
        = 
        \Proj_{\perp}^{k} v_{m} 
        = 
        \Proj_{\perp}^{k} \Defect^{k} v_{m} + \Proj_{\perp}^{k} \cartan^{k-1} \Poincare^{k} v_{m} + \Proj_{\perp}^{k} \Poincare^{k+1} \cartan^{k} v_{m}
        = 
        \Proj_{\perp}^{k} \Defect^{k} v_{m} + \Proj_{\perp}^{k} \Poincare^{k+1} \cartan^{k} v_{m}
        .
    \end{align*}
    By continuity of $\Poincare^{k+1}$, the summand $\Proj_{\perp}^{k} \Poincare^{k+1} \cartan^{k} v_{m}$ converges to zero.
    Recall that $\Defect^{k} : X^{k} \to X^{k}$ is compact. 
    By that compactness, we can restrict $(\Proj_{\perp}^{k} \Defect^{k} v_{m})_{m \in \bbN}$ to a subsequence that is convergent to a vector $v \in X^{k}_{\perp}$.
    That is a unit vector. 
    But since $\cartan^{k} : X^{k} \to X^{k+1}$ is continuous,
    \begin{align*}
        \cartan^{k} v 
        = 
        \lim_{m \rightarrow \infty} \cartan^{k} \Proj_{\perp}^{k} \Defect^{k} v_{m} 
        = 
        \lim_{m \rightarrow \infty} \cartan^{k} \Defect^{k} v_{m} 
        = 
        \lim_{m \rightarrow \infty} \Defect^{k+1} \cartan^{k} v_{m} 
        = 
        0
        .
    \end{align*}
    We conclude that $\cartan^{k} : X^{k} \to X^{k+1}$ has closed range.

    There exists a closed subspace $\Hilbertcomplement^{k} \subseteq X^{k}$ such that we have an orthogonal sum 
    \begin{align*}
        \ker\left( \cartan^{k} : X^{k} \to X^{k+1} \right) = \cartan^{k-1} X^{k-1} \oplus \Hilbertcomplement^{k}
        .
    \end{align*}
    Write $\Projection_{\Hilbertcomplement}^{k} : X^{k} \to \Hilbertcomplement^{k}$ for the orthogonal projection. 
    In addition, the closed range property implies the existence of a bounded operator 
    $\Poincare_{X}^{k+1} : X^{k+1} \to X^{k}$ that is the Moore--Penrose pseudoinverse of $\cartan^{k} : X^{k} \to X^{k+1}$.
    
    Suppose that $z \in \Hilbertcomplement^{k}$. Then $\cartan^{k} z = 0$, and the homotopy formula in Theorem~\ref{theorem:continuity-of-poincare-and-remainder-operator} shows 
    \begin{align*}
        z
        = 
        \Defect^{k} z + \cartan^{k-1} \Poincare^{k} z + \Poincare^{k+1} \cartan^{k} z
        = 
        \Defect^{k} z + \cartan^{k-1} \Poincare^{k} z
        .
    \end{align*}
    Since $\Poincare^{k} z \in X^{k-1}$, we get
    \begin{align*}
         z 
         = 
         \Projection_{\Hilbertcomplement}^{k} z
         = 
         \Projection_{\Hilbertcomplement}^{k} \left( \Defect^{k} z + \cartan^{k-1} \Poincare^{k} z \right)
         = 
         \Projection_{\Hilbertcomplement}^{k} \Defect^{k} z
         .
    \end{align*}
    As $z \in \Hilbertcomplement^{k}$ and $\Defect^{k} : X^{k} \to X^{k}$ is compact, the space $\Hilbertcomplement^{k}$ is finite-dimensional. 
    
    We define the projected subspace and the operator 
    \begin{align*}
        \Harmonic^{k} := \Defect^{k} \Hilbertcomplement^{k}, \qquad \Proj^{k}_{0} := \Defect^{k} \Projection_{\Hilbertcomplement}^{k} \Defect^{k}.
    \end{align*}
    We have $\Harmonic^{k} \subseteq C^{\infty}(\closure\Domain,\Alt{k})$ by Theorem~\ref{theorem:continuity-of-poincare-and-remainder-operator}.
    If $u \in \Hilbertcomplement^{k}$ with $\Defect^{k} u = 0$, then $u = \cartan^{k-1} \Poincare^{k} u$, which implies $u = 0$.
    This means that $\Defect^{k} : \Hilbertcomplement^{k} \to \Harmonic^{k}$ is injective, hence bijective. 
    In particular, 
    \begin{align*}
        \dim \Hilbertcomplement^{k} = \dim \Harmonic^{k}.
    \end{align*}
    
    We review the properties of $\Proj^{k}_{0}$.
    If $u \in \calD'(\Domain,\Alt{k-1})$, then 
    \begin{align*}
        \Proj^{k}_{0} \cartan^{k-1} u = \Defect^{k} \Projection_{\Hilbertcomplement}^{k} \Defect^{k} \cartan^{k-1} u = \Defect^{k} \Projection_{\Hilbertcomplement}^{k} \cartan^{k-1} \Defect^{k-1} u = 0
        .
    \end{align*}
    If $u \in \calD'(\Domain,\Alt{k})$, then 
    \begin{align*}
        \cartan^{k} \Proj^{k}_{0} u = \cartan^{k} \Defect^{k} \Projection_{\Hilbertcomplement}^{k} \Defect^{k} u = \Defect^{k+1} \cartan^{k} \Projection_{\Hilbertcomplement}^{k} \Defect^{k} u = 0
        .
    \end{align*}
    We verify that $\Proj^{k}_{0} : \calD'(\Domain,\Alt{k}) \to \Harmonic^{k}$ is idempotent. 
    To see that, suppose $v \in \Hilbertcomplement^{k}$ and let $z := \Defect^{k} v$. 
    Recall that $\cartan^{k} v = \cartan^{k} z = 0$.
    Using $z \in C^{\infty}(\closure\Domain,\Alt{k})$, we get 
    \begin{align*}
        \Proj^{k}_{0} z
        &
        =
        \Defect^{k} \Projection_{\Hilbertcomplement}^{k} \Defect^{k} z
        =
        \Defect^{k} \Projection_{\Hilbertcomplement}^{k} ( z - \cartan^{k-1} \Poincare^{k} z )
        =
        \Defect^{k} \Projection_{\Hilbertcomplement}^{k} z 
        \\&
        =
        \Defect^{k} \Projection_{\Hilbertcomplement}^{k} \Defect^{k} v 
        =
        \Defect^{k} \Projection_{\Hilbertcomplement}^{k} ( v - \cartan^{k-1} \Poincare^{k} v )
        =
        \Defect^{k} \Projection_{\Hilbertcomplement}^{k} v 
        =
        \Defect^{k} v 
        =
        z
        .
    \end{align*}
    Hence, $\Proj^{k}_{0} : \calD'(\Domain,\Alt{k}) \to \Harmonic^{k}$ is a projection. 
    We also conclude 
    \begin{align}\label{math:smooth-forms:disjointness}
        \cartan^{k-1} \calD'(\Domain,\Alt{k-1}) \cap \Harmonic^{k} = \{0\}
    \end{align}
    because $\Proj^{k}_{0}$ annihilates the first space in the intersection and acts as the identity on the second.
    
    We prove the universal regular decomposition. 
    Let $u \in \calD'(\Domain,\Alt{k})$ with $\cartan^{k} u = 0$. The homotopy formula of Theorem~\ref{theorem:continuity-of-poincare-and-remainder-operator} yields 
    \begin{align*}
        u = \cartan^{k-1} \Poincare^{k} u + \Defect^{k} u
        .
    \end{align*}
    The orthogonal decomposition gives that $w := \Defect^{k} u - \Projection_{\Hilbertcomplement}^{k} \Defect^{k} u \in \cartan^{k-1} X^{k-1}$. 
    Consequently, $\Defect^{k} u = \Projection_{\Hilbertcomplement}^{k} \Defect^{k} u + w$ with $w = \cartan^{k-1} \Poincare_{X}^{k} \Defect^{k} u$, and we have
    \begin{align*}
        \Defect^{k} u = \Projection_{\Hilbertcomplement}^{k} \Defect^{k} u + \cartan^{k-1} \Poincare_{X}^{k} \Defect^{k} u.
    \end{align*}
    The last decomposition is orthogonal in $X^{k}$. Now,
    \begin{align*}
        \Defect^{k} u
        &
        = 
        \Defect^{k} \cartan^{k-1} \Poincare^{k} u + \Defect^{k} \Defect^{k} u
        \\&
        = 
        \cartan^{k-1} \Defect^{k-1} \Poincare^{k} u + \Defect^{k} \Projection_{\Hilbertcomplement}^{k} \Defect^{k} u + \Defect^{k} \cartan^{k-1} \Poincare_{X}^{k} \Defect^{k} u
        = 
        \cartan^{k-1} \left( \Defect^{k-1} \Poincare^{k} u + \Defect^{k-1} \Poincare_{X}^{k} \Defect^{k} u \right) + \Defect^{k} \Projection_{\Hilbertcomplement}^{k} \Defect^{k} u
        .
    \end{align*}
    It follows that 
    \begin{align*}
        u 
        = 
        \cartan^{k-1} \left( \Poincare^{k} u + \Defect^{k-1} \Poincare^{k} u + \Defect^{k-1} \Poincare_{X}^{k} \Defect^{k} u \right) 
        + 
        \Defect^{k} \Projection_{\Hilbertcomplement}^{k} \Defect^{k} u
        .
    \end{align*}
    Define the operator 
    \begin{align*}
        \widehat\Poincare^{k}_{0} := \Poincare^{k} + \Defect^{k-1} \Poincare^{k} + \Defect^{k-1} \Poincare_{X}^{k} \Defect^{k}
        .
    \end{align*}
    
    We have shown that if $u \in \calD'(\Domain,\Alt{k})$ with $\cartan^{k} u = 0$, then 
    \begin{align}\label{math:kernel-decomposition}
        u = \cartan^{k-1} \widehat\Poincare^{k}_{0} u + \Proj_{0}^{k} u
        .
    \end{align}
    As an immediate consequence, 
    \begin{align*}
        \cartan^{k} \widehat\Poincare^{k+1}_{0} \cartan^{k} u = \cartan^{k} u.
    \end{align*}
    For any $u \in \calD'(\Domain,\Alt{k})$, write $v := u - \widehat\Poincare^{k+1}_{0} \cartan^{k} u$. Thus, $\cartan^{k} v = 0$, and 
    \begin{align*}
        u 
        = 
        \widehat\Poincare^{k+1}_{0} \cartan^{k} u + v
        = 
        \widehat\Poincare^{k+1}_{0} \cartan^{k} u + \cartan^{k-1} \widehat\Poincare^{k}_{0} v + \Proj_{0}^{k} v
        .
    \end{align*}

    Define the operators 
    \begin{align*}
        \widehat\Poincare^{k} := \widehat\Poincare^{k}_{0} \left( \Id - \widehat\Poincare^{k+1}_{0} \cartan^{k} \right),
        \qquad 
        \Proj^{k} := \Proj^{k}_{0} \left( \Id - \widehat\Poincare^{k+1}_{0} \cartan^{k} \right)
        .
    \end{align*}
    If $u \in \calD'(\Domain,\Alt{k})$, then 
    \begin{gather*}
        u 
        = 
        \widehat\Poincare^{k+1} \cartan^{k} u + \cartan^{k-1} \widehat\Poincare^{k} u + \Proj^{k} u
        ,
        \qquad 
        \cartan^{k} u 
        = 
        \cartan^{k} \widehat\Poincare^{k+1} \cartan^{k} u
        ,
        \\
        \cartan^{k} u = 0 \implies u = \cartan^{k-1} \widehat\Poincare^{k}_{0} u + \Proj^{k} u,
        \\
        \cartan^{k} \Proj^{k} u = 0, \qquad \Proj^{k+1} \cartan^{k} u = 0, \qquad \Proj^{k} \Proj^{k} u = \Proj^{k} u
        .
    \end{gather*}
    The continuity conditions~\eqref{math:universal-regular-decomposition:continuity} hold after composition.

    Now,~\eqref{math:smooth-forms:disjointness} together with the decomposition~\eqref{math:kernel-decomposition}
    via Theorem~\ref{theorem:continuity-of-poincare-and-remainder-operator} 
    shows the topological direct sums~\eqref{math:universal-regular-decomposition:kernel}. 
    The proof is complete. 
\end{proof}

\begin{remark}
    The auxiliary differential complex in the proof of Theorem~\ref{theorem:universal-regular-decomposition} is an arbitrary choice; different choices need not lead to the same universal complement space. In particular, the universal regular decomposition is not identical to the orthogonal Hodge decomposition of those Hilbert complexes. 
    
    We also remark that the last harmonic complement space $\Harmonic^{n}(\Domain)$ of an $n$-dimensional bounded Lipschitz domain is always trivial,
    i.e., $\dim \Harmonic^{n}(\Domain) = 0$, as is known from algebraic topology. 
    In particular, the homotopy operator for the divergence always finds a potential vector field whose divergence equals the original source field. 
\end{remark}

\section{Universal extension of differential forms}\label{section:extension-of-forms}

The combination of the universal extension operators on scalar fields and the universal regular decomposition on de~Rham complexes leads to the universal extension operators of differential forms.
We first present an operator that extends differential forms on a bounded Lipschitz domain onto the whole ambient space and preserves the regularity of both the differential form and its exterior derivative; this preservation of regularity is what we refer to as \emph{compatible}.
For example, if the differential form is in a Sobolev--Slobodeckij space and its exterior derivative has additional regularity, then this will be reflected in the extension onto the enclosing ambient space.

The second extension operator modifies that simple construction: 
it is not onto the whole space but only onto a collar neighborhood of the original Lipschitz domain. 
However, it commutes with the exterior derivative and induces isomorphisms
between the cohomology groups on the Lipschitz domain and its extension.

An exploration of the smooth de~Rham cohomology on a bounded Lipschitz domain subject to a partial flatness condition 
is a critical auxiliary result and also a contribution of independent interest.

We now develop the universal extension operators for differential forms on bounded Lipschitz domains. 

\subsection{Compatible extension of differential forms}

We write $\Extension^{k}$ for the componentwise universal extension operator of Theorem~\ref{theorem:extension-componentwise} acting on $k$-forms.

\begin{theorem}\label{theorem:universal-extension}
    Let $\Domain \subseteq \bbR^{n}$ be a bounded Lipschitz domain and $k \in \bbZ$. 
    Let $\widehat\Poincare^{k}$ and $\Proj^{k}$ be the corrected potential and projection from Theorem~\ref{theorem:universal-regular-decomposition}.
    We have a continuous operator 
    \begin{align}\label{math:universal-extension}
        \Erweiterung^{k} : 
        \calS'(\Domain,\Alt{k})
        \to 
        \calS'(\bbR^{n},\Alt{k})
        ,
        \qquad 
        u \mapsto \Extension^{k} \widehat\Poincare^{k+1} \cartan^{k} u + \cartan^{k-1} \Extension^{k-1} \widehat\Poincare^{k} u + \Extension^{k} \Proj^{k} u
        ,
    \end{align}
    such that for all $u \in \calS'(\Domain,\Alt{k})$ we have 
    \begin{align}\label{math:universal-extension:restriction}
        ( \Erweiterung^{k} u) \restriction_{\Domain} = u.
    \end{align}
    Moreover, we have a continuous operator 
    \begin{align}\label{math:universal-extension:smooth}
        \Erweiterung^{k} : C^{\infty}(\closure\Domain,\Alt{k}) \to C^{\infty}(\bbR^{n},\Alt{k})
        .
    \end{align}
    Each of the following inequalities holds if its right-hand side is finite:
    \begin{subequations}
    \begin{itemize}
        \item 
        If $s \in \bbR$ and $p \in (1,\infty)$, then 
        \begin{align}
            \| \Erweiterung^{k} u \|_{\Wres^{s,p}(\bbR^{n},\Alt{k})} &\leq C \| u \|_{\Wres^{s,p}(\Domain,\Alt{k})}
            , \label{math:universal-extension:sobolev}
            \\
            \| \cartan^{k} \Erweiterung^{k} u \|_{\Wres^{s,p}(\bbR^{n},\Alt{k+1})} &\leq C \left( \| \cartan^{k} u \|_{\Wres^{s,p}(\Domain,\Alt{k+1})} + \| \Proj^{k} u \|_{\Wres^{s+1,p}(\Domain,\Alt{k})} \right)
            . \label{math:universal-extension:sobolev-cartan}
        \end{align}
        \item 
        If $s \in \bbR$ and $p,q \in (1,\infty]$, then 
        \begin{align}
            \| \Erweiterung^{k} u \|_{B^{s}_{p,q}(\bbR^{n},\Alt{k})} &\leq C \| u \|_{B^{s}_{p,q}(\Domain,\Alt{k})}
            ,  \label{math:universal-extension:besov}
            \\
            \| \cartan^{k} \Erweiterung^{k} u \|_{B^{s}_{p,q}(\bbR^{n},\Alt{k+1})} &\leq C \left( \| \cartan^{k} u \|_{B^{s}_{p,q}(\Domain,\Alt{k+1})} + \| \Proj^{k} u \|_{B^{s+1}_{p,q}(\Domain,\Alt{k})} \right)
            . \label{math:universal-extension:besov-cartan}
        \end{align}
        \item 
        If $s \in \bbR$ and $p,q \in (1,\infty]$ with $p < \infty$, then 
        \begin{align}
            \| \Erweiterung^{k} u \|_{F^{s}_{p,q}(\bbR^{n},\Alt{k})} &\leq C \| u \|_{F^{s}_{p,q}(\Domain,\Alt{k})}
            , \label{math:universal-extension:triebel-lizorkin}
            \\
            \| \cartan^{k} \Erweiterung^{k} u \|_{F^{s}_{p,q}(\bbR^{n},\Alt{k+1})} &\leq C \left( \| \cartan^{k} u \|_{F^{s}_{p,q}(\Domain,\Alt{k+1})} + \| \Proj^{k} u \|_{F^{s+1}_{p,q}(\Domain,\Alt{k})} \right)
            . \label{math:universal-extension:triebel-lizorkin-cartan}
        \end{align}
    \end{itemize}
    \end{subequations}
    Here, the constant $C$ only depends on $\Domain$ and the relevant parameters $s,p,q$.
\end{theorem}
\begin{proof}
    The continuity conditions~\eqref{math:universal-extension} and~\eqref{math:universal-extension:smooth} follow from Theorems~\ref{theorem:extension-componentwise} and~\ref{theorem:universal-regular-decomposition}.
    Throughout the proof, let $u \in \calS'(\Domain,\Alt{k})$. 
    To show the extension property, we observe 
    \begin{align*}
        ( \Erweiterung^{k} u ) \restriction_{\Domain}
        &
        =
        ( \cartan^{k-1} \Extension^{k-1} \widehat\Poincare^{k} u ) \restriction_{\Domain}
        + 
        ( \Extension^{k} \Proj^{k} u ) \restriction_{\Domain}
        + 
        ( \Extension^{k} \widehat\Poincare^{k+1} \cartan^{k} u ) \restriction_{\Domain}
        \\&
        =
        \cartan^{k-1} ( \Extension^{k-1} \widehat\Poincare^{k} u ) \restriction_{\Domain}
        + 
        ( \Extension^{k} \Proj^{k} u ) \restriction_{\Domain}
        + 
        ( \Extension^{k} \widehat\Poincare^{k+1} \cartan^{k} u ) \restriction_{\Domain}
        \\&
        =
        \cartan^{k-1} \widehat\Poincare^{k} u
        + 
        \Proj^{k} u
        + 
        \widehat\Poincare^{k+1} \cartan^{k} u
        =
        u
        .
    \end{align*}
    The last equality is the homotopy identity in Theorem~\ref{theorem:universal-regular-decomposition}.
    It remains to prove the other continuity properties. 
    We begin with the Sobolev--Slobodeckij scale. 
    Let $s,t \in \bbR$ and $p \in (1,\infty)$ with $t \geq s-1$. 
    In what follows, we write $A \lesssim B$ if $A \leq C B$ holds with a constant $C \geq 0$
    that only depends on $\Domain$, $p$, $s$, and $t$. 
    
    From the definition of the extension operator, 
    \begin{align*}
        \| \Erweiterung^{k} u \|_{\Wres^{s,p}(\bbR^{n},\Alt{k})}
        &
        \lesssim 
        \| 
            \Extension^{k} \widehat\Poincare^{k+1} \cartan^{k} u
            + 
            \cartan^{k-1} \Extension^{k-1} \widehat\Poincare^{k} u
            + 
            \Extension^{k} \Proj^{k} u 
        \|_{\Wres^{s,p}(\bbR^{n},\Alt{k})}
        \\&
        \lesssim
        \| \Extension^{k} \widehat\Poincare^{k+1} \cartan^{k} u \|_{\Wres^{s,p}(\bbR^{n},\Alt{k})}
        + 
        \| \cartan^{k-1} \Extension^{k-1} \widehat\Poincare^{k} u \|_{\Wres^{s,p}(\bbR^{n},\Alt{k})}
        +
        \| \Extension^{k} \Proj^{k} u \|_{\Wres^{s,p}(\bbR^{n},\Alt{k})}
        ,
        \\
        \| \cartan^{k} \Erweiterung^{k} u \|_{\Wres^{t,p}(\bbR^{n},\Alt{k+1})}
        &
        \lesssim 
        \| 
            \cartan^{k} \Extension^{k} \widehat\Poincare^{k+1} \cartan^{k} u
            + 
            \cartan^{k} \Extension^{k} \Proj^{k} u 
        \|_{\Wres^{t,p}(\bbR^{n},\Alt{k+1})}
        \\&
        \lesssim
        \| \cartan^{k} \Extension^{k} \widehat\Poincare^{k+1} \cartan^{k} u \|_{\Wres^{t,p}(\bbR^{n},\Alt{k+1})}
        +
        \| \cartan^{k} \Extension^{k} \Proj^{k} u \|_{\Wres^{t,p}(\bbR^{n},\Alt{k+1})}
        .
    \end{align*}
    We address the harmonic terms, using Theorem~\ref{theorem:extension-componentwise} and the fact that $\Proj^{k}$ has finite-dimensional range:
    \begin{align*}
        &
        \| \Extension^{k} \Proj^{k} u \|_{\Wres^{s,p}(\bbR^{n},\Alt{k})}
        \lesssim
        \| \Proj^{k} u \|_{\Wres^{s,p}(\Domain,\Alt{k})}
        \lesssim
        \| u \|_{\Wres^{s,p}(\Domain,\Alt{k})}
        ,
        \\
        &
        \| \cartan^{k} \Extension^{k} \Proj^{k} u \|_{\Wres^{t,p}(\bbR^{n},\Alt{k+1})}
        \lesssim
        \| \Extension^{k} \Proj^{k} u \|_{\Wres^{t+1,p}(\bbR^{n},\Alt{k})}
        \lesssim
        \| \Proj^{k} u \|_{\Wres^{t+1,p}(\Domain,\Alt{k})}
        \lesssim
        \| u \|_{\Wres^{s,p}(\Domain,\Alt{k})}
        .
    \end{align*}
    Next, we notice 
    \begin{align*}
        \| \cartan^{k-1} \Extension^{k-1} \widehat\Poincare^{k} u \|_{\Wres^{s,p}(\bbR^{n},\Alt{k})}
        \lesssim
        \| \Extension^{k-1} \widehat\Poincare^{k} u \|_{\Wres^{s+1,p}(\bbR^{n},\Alt{k-1})}
        \lesssim
        \| \widehat\Poincare^{k} u \|_{\Wres^{s+1,p}(\Domain,\Alt{k-1})}
        \lesssim
        \| u \|_{\Wres^{s,p}(\Domain,\Alt{k})}
        .
    \end{align*}
    A similar combination of techniques gives 
    \begin{align*}
        \| \Extension^{k} \widehat\Poincare^{k+1} \cartan^{k} u \|_{\Wres^{s,p}(\bbR^{n},\Alt{k})}
        \lesssim
        \| \widehat\Poincare^{k+1} \cartan^{k} u \|_{\Wres^{s,p}(\Domain,\Alt{k})}
        \lesssim
        \| \cartan^{k} u \|_{\Wres^{s-1,p}(\Domain,\Alt{k+1})}
        \lesssim
        \| u \|_{\Wres^{s,p}(\Domain,\Alt{k})}
        ,
    \end{align*}
    and 
    \begin{align*}
        \| \cartan^{k} \Extension^{k} \widehat\Poincare^{k+1} \cartan^{k} u \|_{\Wres^{t,p}(\bbR^{n},\Alt{k+1})}
        &
        \lesssim
        \| \Extension^{k} \widehat\Poincare^{k+1} \cartan^{k} u \|_{\Wres^{t+1,p}(\bbR^{n},\Alt{k})}
        \\&
        \lesssim
        \| \widehat\Poincare^{k+1} \cartan^{k} u \|_{\Wres^{t+1,p}(\Domain,\Alt{k})}
        \lesssim
        \| \cartan^{k} u \|_{\Wres^{t,p}(\Domain,\Alt{k+1})}
        .
    \end{align*}
    This shows the desired continuity estimates for the Sobolev--Slobodeckij spaces. 
    
    The proof of the continuity estimates for Besov spaces and Triebel--Lizorkin spaces of differential forms is fully analogous.
\end{proof}

\begin{remark}
    The exterior derivative of the universal extension operator depends only on the exterior derivative of the original differential form 
    and its projection onto the universal harmonic complement. 
    The exact part in the universal regular decomposition does not enter the extension of the exterior derivative.
    We point out that on any contractible domain, the harmonic complement vanishes in all degrees except degree zero, where it corresponds to the constant functions.
\end{remark}

\subsection{Aspects of sheaf theory}

The main result of this manuscript will build on sheaf theory, for which we now give a short overview.
Our exposition is inspired by work in global analysis on manifolds~\cite{goldshtein2011hodge,mitrea2008poisson},
and the reader is referred to the literature for additional background on sheaf~theory~\cite{stocker_zieschang}. %

Suppose that $\calX$ is a topological space whose collection of open sets we denote by $\calO$. 
A \emph{sheaf} on $\calX$ is a collection of vector spaces $\{ \calF(U) \}_{U \in \calO}$
together with a family of linear mappings $\{ \rho_{V,U} \}_{ V, U \in \calO, V \subseteq U}$ such that the following conditions hold:
\begin{itemize}
    \item If $V, U \in \calO$ and $V \subseteq U$, then we have a linear map $\rho_{V,U} : \calF(U) \to \calF(V)$.
    \item For each $U \in \calO$, we have $\rho_{U,U} = \Id_{\calF(U)}$.
    \item For all $U,V,W \in \calO$ with $W \subseteq V \subseteq U$, we have 
    \begin{align*}
        \rho_{W,V} \rho_{V,U} = \rho_{W,U}.
    \end{align*}
    \item 
    If $U \in \calO$ and $\Cover = \{ U_{i} \}_{i \in I} \subseteq \calO$ is an open cover of $U$, 
    and the family $\{ f_{i} \}_{i \in I}$, with $f_{i} \in \calF(U_i)$, satisfies the compatibility condition
    \begin{align*}
        \rho_{U_{i} \cap U_{j}, U_{i}} f_{i} = \rho_{U_{i} \cap U_{j}, U_{j}} f_{j}, \qquad i,j \in I
        ,
    \end{align*}
    then there exists a unique $f \in \calF(U)$ such that $\rho_{U_{i},U}(f) = f_{i}$ for all $i \in I$. 
\end{itemize}
The mappings $\rho_{V,U}$ are called \emph{restriction operators}. 
We follow the convention of abbreviating $f \restriction_{V} := \rho_{V,U}(f)$ for any $V, U \in \calO$ with $V \subseteq U$ and $f \in \calF(U)$.

A \emph{sheaf homomorphism} $\theta : \calF \to \calG$ from a sheaf $\calF$ to another sheaf $\calG$ on $\calX$ is a collection of linear maps 
\begin{align*}
    \{ \theta_{U} : \calF(U) \to \calG(U) \}_{U \in \calO}
\end{align*}
that commute with the restriction operators, i.e., $\rho_{V,U} \theta_{U} = \theta_{V} \rho_{V,U}$ for all $U, V \in \calO$ with $V \subseteq U$. 
The support $\supp(\theta)$ of the sheaf homomorphism is the smallest closed subset outside of which $\theta$ restricts to the trivial homomorphism.

A sheaf is called \emph{fine} if for each locally finite open cover $\Cover = \{ U_{i} \}_{i \in I}$ of $\calX$
there exists a family of sheaf homomorphisms $\{ \theta_{i} : \calF \to \calF \}_{i \in I}$ such that
\begin{align*}
    \sum_{i \in I} \theta_{i} = \Id
    ,
    \qquad 
    \supp(\theta_{i}) \subseteq U_{i}, \qquad i \in I
    .
\end{align*}
We emphasize that the last sum is locally finite and hence well-defined. In practice, these families of sheaf homomorphisms correspond to multiplication by partition-of-unity functions subordinate to the given open cover. 

Suppose that $\calF^{0}, \calF^{1}, \calF^{2}, \dots$ is a sequence of sheaves over the space $\calX$
and that we have a sequence of sheaf homomorphisms $\theta^{k} : \calF^{k} \to \calF^{k+1}$ such that
\begin{align*}
    \theta^{k+1} \theta^{k} = 0.
\end{align*}
We have a well-defined diagram 
\begin{align}\label{math:differential-complex-of-sheaves:rewrite}
    \begin{CD}
        0 @>>> \calF^{0}  @>\theta^{0}>> \calF^{1} @>\theta^{1}>> \calF^{2} @>\theta^{2}>> \dots 
    \end{CD}
\end{align}
that we call a differential complex of sheaves. 

The differential complex of sheaves~\eqref{math:differential-complex-of-sheaves:rewrite} is called \emph{exact in positive degrees} if the following condition is satisfied for all $k \in \bbN$:
if $U \in \calO$ and $f \in \calF^{k}(U)$ with $\theta^{k} f = 0$, 
then for each $x \in U$ there exists $V \in \calO$ with $x \in V \subseteq U$ and $F \in \calF^{k-1}(V)$ such that $\theta^{k-1}(F) = \rho_{V,U}(f)$.

We consider the sheaf of locally constant functions on the space $\calX$ that vanish on a fixed subset $\Gamma \subseteq \calX$. 
Here, we define 
\begin{align*}
    \bbR_{\calX,\Gamma}(U) 
    := 
    \left\{ 
        u : U \to \bbR \suchthat* u \restriction_{U \cap \Gamma} = 0, \ u \textrm{ is locally constant}
    \right\}
\end{align*}
for each $U \in \calO$ together with the obvious restriction operators. Note that the restrictions are not necessarily surjective.

\begin{theorem}\label{theorem:resolution-of-locally-constant-functions}
    Let $\calX$ be a Hausdorff topological space, and let $\Gamma \subseteq \calX$ be a closed subspace.
    Suppose that every subspace of $\calX$ is paracompact,
    and that the topology of $\calX$ and the relative topology of $\Gamma$ have bases of contractible open sets.
    If we have a differential complex of sheaves 
    \begin{align} 
        \begin{CD}
            0
            @>>>
            \calF^{0}
            @>{\theta^{0}}>>
            \calF^{1}
            @>{\theta^{1}}>>
            \calF^{2}
            @>{\theta^{2}}>>
            \dots 
        \end{CD}
    \end{align}
    that is exact in positive degrees, 
    where the sheaves $\calF^{0}, \calF^{1}, \calF^{2}, \dots$ are fine, 
    and where 
    \[
        \ker\left( \theta^0 : \calF^{0} \to \calF^{1} \right)
        =
        \bbR_{\calX,\Gamma}
        ,
    \]
    then we have vector space isomorphisms 
    \begin{align*}
        H^0_{\mathrm{sing}}(\calX,\Gamma;\bbR)
        &
        \simeq 
        \ker\left( \theta^0 : \calF^{0}(\calX) \to \calF^{1}(\calX) \right)
        ,
        \\
        H^{k}_{\mathrm{sing}}(\calX,\Gamma;\bbR)
        &
        \simeq 
        \frac{ 
            \ker\left( \theta^{k} : \calF^{k}(\calX) \to \calF^{k+1}(\calX) \right)
        }{
            \rng\left( \theta^{k-1} : \calF^{k-1}(\calX) \to \calF^{k}(\calX) \right)
        } 
        ,
        \qquad 
        k \in \bbN
        .
    \end{align*}
\end{theorem}
\begin{proof}
    Under the assumptions of the theorem, the diagram
    \begin{align}\label{math:resolution-of-locally-constant-functions:rewrite}
        \begin{CD}
            0
            @>>>
            \bbR_{\calX,\Gamma}
            @>>>
            \calF^{0}
            @>{\theta^{0}}>>
            \calF^{1}
            @>{\theta^{1}}>>
            \calF^{2}
            @>{\theta^{2}}>>
            \dots 
        \end{CD}
    \end{align}
    is a so-called \emph{acyclic resolution} of the sheaf $\bbR_{\calX,\Gamma}$;
    see~\cite[Lemma~II.3.3]{bartocci1991geometry}, also~\cite[Satz 4.5.12]{deitmar2025kohomologie},~\cite[Satz 2.17, Satz 3.2]{cap1998garben}.
    Therefore, the cohomology groups of the differential complex of global sections 
    \begin{align}\label{math:complex-of-global-sections:rewrite}
        \begin{CD}
            0
            @>>>
            \calF^{0}(\calX)
            @>{\theta^{0}}>>
            \calF^{1}(\calX)
            @>{\theta^{1}}>>
            \calF^{2}(\calX)
            @>{\theta^{2}}>>
            \dots 
        \end{CD}
    \end{align}
    are isomorphic to the so-called \emph{sheaf cohomology groups} of the sheaf $\bbR_{\calX,\Gamma}$;
    see~\cite[Proposition~II.2.4]{bartocci1991geometry}, also~\cite[Satz 3.5]{cap1998garben}.
    
    The latter cohomology groups can be shown to be isomorphic to the singular cohomology groups of $\calX$ relative to $\Gamma$ with coefficients in the real numbers.\footnote{While this is folklore in sheaf theory, it is surprisingly difficult to find an explicit reference in the literature.}
    An explicit statement
    \footnote{Theorem~4.3.15 in~\cite{theeigenbook} seems to be stated best with an additional hypothesis, e.g., that the subspace be locally contractible, too. As a counterexample that violates that hypothesis, consider the case $\calX := [0,1]$ and $\Gamma := \{0,1,1/2,1/3,\dots\}$ and the cohomologies of $\calX$ relative to $\Gamma$: the first sheaf cohomology is isomorphic to $\bigoplus_{i=0}^\infty \bbR$ and the first singular cohomology is isomorphic to $\prod_{i=0}^\infty \bbR$.}
    is found in the textbook by Bellaïche~\cite{theeigenbook},
    namely Theorem~4.3.15 together with the definitions in Sections~4.3.1--6 and~4.3.11. 
    In combination, the claimed theorem follows.
\end{proof}

\subsection{Aspects of smooth de~Rham cohomology}

Let $\Domain$ be a bounded Lipschitz domain.
Let $\Gamma \subseteq \partial\Domain$ be a relatively open subset of the boundary. 
We recall that there are topological invariants 
\begin{align*}
    b_{0}(\closure\Domain,\Gamma), \quad b_{1}(\closure\Domain,\Gamma), \quad \dots, \quad b_{n}(\closure\Domain,\Gamma),
\end{align*}
which are known as the Betti numbers of $\closure\Domain$ relative to $\Gamma$. 
The $k$-th Betti number is the dimension of the $k$-th singular cohomology group over the real numbers of $\closure\Domain$ relative to the subspace $\Gamma$: 
\begin{align}\label{math:definition-of-betti-numbers}
    b_{k}(\closure\Domain,\Gamma) := \dim H^{k}_{\mathrm{sing}}(\closure\Domain,\Gamma,\bbR).
\end{align}
The Betti numbers are a standard concept in algebraic topology and reflect topological features of the domain~\cite{stocker_zieschang}.
The absolute Betti numbers $b_{k}(\closure\Domain) := b_{k}( \closure\Domain, \emptyset )$ are more widely discussed in the literature.
For example, $b_{0}(\closure\Domain) = b_{0}(\closure\Domain,\emptyset)$ is the number of connected components of $\closure\Domain$. 
The Betti numbers are closely related to the cohomology groups of de~Rham complexes,
which is well-known for the absolute Betti numbers but deserves some exposition for the relative Betti numbers.

We introduce the vector space $C^{\infty}(\Domain,\Gamma,\VEC) \subseteq C^{\infty}(\Domain,\VEC)$
by requiring that each $u \in C^{\infty}(\Domain,\Gamma,\VEC)$ satisfies the following condition:
every $x \in \Gamma$ has an open neighborhood $\subdomain_x \subseteq \bbR^{n}$
such that the zero extension of $u$ to $\subdomain_x$ (i.e., the function that coincides with $u$ on $\subdomain_x \cap \Domain$ and vanishes on $\subdomain_x \setminus \Domain$) belongs to $C^{\infty}(\subdomain_x,\VEC)$. 
Note that this encodes a boundary condition that is more restrictive than the common Dirichlet condition: 
the function and its derivatives of all orders are required to vanish towards $\Gamma$.

It follows from the definition that the exterior derivative preserves the boundary condition, so that we have well-defined mappings 
\begin{align*}
    \cartan^{k} : C^{\infty}(\Domain,\Gamma,\Alt{k}) \to C^{\infty}(\Domain,\Gamma,\Alt{k+1})
    .
\end{align*}
In particular, we have a differential complex 
\begin{align}\label{math:relative-smooth-homology:complex-on-domain}
    \begin{CD}
        0 \to 
        C^{\infty}(\Domain,\Gamma,\Alt{0}) 
        @>>> 
        C^{\infty}(\Domain,\Gamma,\Alt{1}) 
        @>>> 
        \dots
        @>>> 
        C^{\infty}(\Domain,\Gamma,\Alt{n})
        \to 0.
    \end{CD}
\end{align}
We study the cohomology of the smooth de~Rham complex subject to these boundary conditions.

\begin{proposition}\label{proposition:relative-smooth-homology}
    Consider a bounded Lipschitz domain $\Domain$.
    Let $\Gamma$ be a union of connected components of $\partial\Domain$. 
    Then the differential complex~\eqref{math:relative-smooth-homology:complex-on-domain}
    realizes the Betti numbers of $\closure\Domain$ relative to $\Gamma$ on cohomology. In particular, 
    \begin{align*}
        \dim 
        \frac{ 
            \ker\left( \cartan^{k} : C^{\infty}(\Domain,\Gamma,\Alt{k}) \to C^{\infty}(\Domain,\Gamma,\Alt{k+1}) \right)
        }{
            \rng\left( \cartan^{k-1} : C^{\infty}(\Domain,\Gamma,\Alt{k-1}) \to C^{\infty}(\Domain,\Gamma,\Alt{k}) \right)
        } 
        = 
        b_{k}(\closure\Domain,\Gamma)
        .
    \end{align*}
\end{proposition}
\begin{proof}
    Recall that $\Domain$ is connected.
    If $\Gamma = \emptyset$, then $b_{0}(\closure\Domain,\emptyset) = 1$ and $\bboldone_{\Domain} \in C^{\infty}(\Domain,\Gamma,\Alt{0})$.
    If $\Gamma \neq \emptyset$, then $b_{0}(\closure\Domain,\Gamma) = 0$ and $\bboldone_{\Domain} \notin C^{\infty}(\Domain,\Gamma,\Alt{0})$.
    In order to address the degrees $k \geq 1$,
    we use the elements of sheaf theory introduced above so that Theorem~\ref{theorem:resolution-of-locally-constant-functions} can be applied. 
    
    Write $\calX := \Domain \cup \Gamma$. 
    If $U \subseteq \calX$ is a relatively open set, 
    then $U \cap \Gamma$ is relatively open in $\partial U$, and we can define 
    \begin{align*}
        \calF^{k}(U) :=  C^{\infty}( U \cap \Domain, U \cap \Gamma, \Alt{k})
        .
    \end{align*}
    If $U,V \subseteq \calX$ are relatively open sets with $V \subseteq U$, 
    then the restriction yields operators 
    \begin{gather*}
        \rho_{V,U} : \calF^{k}(U) \to \calF^{k}(V)
        .
    \end{gather*}
    The families $\{ \calF^{k}(U) \}_{U \subseteq \calX}$ together with the restriction operators $\{ \rho_{V,U} \}_{V \subseteq U}$ form a sheaf on $\calX$.
    
    We show that the sheaves $\calF^{k}$ are fine.
    Let $\Cover = \{ U_{i} \}_{i \in I}$ be a locally finite relatively open cover of $\calX$.
    We can find a family of smooth functions $( \phi_{i} )_{i \in I}$ such that 
    \begin{gather*}
        \forall i \in I : \supp(\phi_{i}) \cap \calX \subseteq U_{i}
        ,
        \\
        \sum_{i \in I} \phi_{i} \equiv 1 \text{ on an open neighborhood of } \calX
        .
    \end{gather*}
    We define sheaf homomorphisms
    \begin{align*}
        \theta_{i} : \calF^{k} \to \calF^{k}, \qquad i \in I,
    \end{align*}
    by multiplication with the partition-of-unity functions $\phi_{i}$.
    It is clear that $\supp(\theta_{i}) \subseteq U_{i}$ and that $\sum_{i \in I} \theta_{i} = \Id$.
    Hence, the sheaves $\calF^{k}$ are fine.

    Recall the sheaf of locally constant real-valued functions on $\calX$ that vanish on $\Gamma$,
    which we write as $\bbR_{\calX,\Gamma}$.
    We have a well-defined sheaf homomorphism
    \begin{align*}
        \imath : \bbR_{\calX,\Gamma} \to \calF^{0}
    \end{align*}
    that gives an injection on each open set. 
    
    The exterior derivatives provide sheaf homomorphisms
    \begin{align*}
        \cartan^{k} : \calF^{k} \to \calF^{k+1}, \qquad k \in \bbZ 
        .
    \end{align*}
    We have a differential complex of sheaves 
    \begin{align}\label{math:relative-smooth-homology:complex-of-sheaves}
        \begin{CD}
            0 \to 
            \bbR_{\calX,\Gamma} 
            @>\imath>>
            \calF^{0} 
            @>\cartan^{0}>>
            \calF^{1} 
            @>\cartan^{1}>>
            \dots
            @>\cartan^{n-1}>>
            \calF^{n}
            \to 0
        \end{CD}
        .
    \end{align}
    Clearly, $\calX$ is Hausdorff and locally contractible.
    Since it is a metric space, all its topological subspaces are paracompact. 
    The subspace $\Gamma$ is locally contractible. 
    To satisfy the conditions of Theorem~\ref{theorem:resolution-of-locally-constant-functions},
    it remains to show that~\eqref{math:relative-smooth-homology:complex-of-sheaves} is an exact complex of sheaves.

    Consider now any relatively open set $U \subseteq \calX$ and let $w \in \calF^{k}(U)$ with $\cartan^{k} w = 0$. 
    We want to show that the closed form $w$ is locally exact, i.e., each $p \in U$ admits an open neighborhood on which $w$ is exact. 
    
    Recall that any interior point $p \in \Domain$ is contained within an open ball $V \subseteq \Domain$ with positive distance from $\partial\Domain$. 
    It then follows that the differential complex of sheaves~\eqref{math:relative-smooth-homology:complex-of-sheaves} gives rise to a local complex on $V$ that is just the smooth de~Rham complex 
    \begin{align*}
        \begin{CD}
            0 \to \bbR 
            @>\imath>>
            C^{\infty}(V,\Alt{0}) 
            @>\cartan^{  0}>> 
            C^{\infty}(V,\Alt{1}) 
            @>\cartan^{  1}>> 
            \dots
            @>\cartan^{n-1}>> 
            C^{\infty}(V,\Alt{n})
            \to 0
            .
        \end{CD}
    \end{align*}
    Here, $\imath : \bbR \to C^{\infty}(V,\Alt{0})$ denotes the inclusion of constant functions.
    This differential complex is exact. 
    Hence, any $p \in U \cap \Domain$ admits an open neighborhood $V \subseteq U$ such that $w \restriction_{V}$ is exact.

    Consider any point $p \in U \cap \Gamma$. Here, we use the Lipschitz property of the domain. 
    Let $\subdomain_{p} \subseteq \bbR^{n}$ be an open neighborhood of $p$ such that the extension-by-zero $\widetilde{w} \in C^{\infty}(\subdomain_{p},\Alt{k})$ is a smooth extension of $w$ on $\subdomain_{p}$.
    Lemma~\ref{lemma:lipschitz-graphs-have-star-domains} establishes that every open neighborhood of $p$ 
    contains a convex Lipschitz chart domain $V$ with 
    \begin{align*}
        p \in V \subseteq \subdomain_{p}
        ,
        \qquad 
        V \cap \calX \subseteq U
        ,
    \end{align*}
    and $V \setminus \calX$ is star-shaped with respect to a ball $\smallball \subseteq V \setminus \calX$.
    Let $b_0 \in \bbR^{n}$ denote the center of that ball. 
    We recall the classical homotopy operator $\Poincare^{k}_{b_{0}} : C^{\infty}(V,\Alt{k}) \to C^{\infty}(V,\Alt{k-1})$ with pivot at $b_0$. 
    Write $\widetilde{\xi} := \Poincare^{k}_{b_{0}} \widetilde{w} \in C^{\infty}(V,\Alt{k-1})$. 
    It is clear that $\cartan^{k-1} \widetilde{\xi} = \widetilde{w}$. 
    In addition, $\widetilde{\xi} \restriction_{ V \setminus \Domain } = 0$ by construction of the homotopy operator,
    which implies $\xi := \widetilde{\xi} \restriction_{V \cap \Domain} \in C^{\infty}(V \cap \Domain,V \cap \Gamma,\Alt{k-1})$ with $\cartan^{k-1} \xi = w$.
    
    We have verified the conditions of Theorem~\ref{theorem:resolution-of-locally-constant-functions}.
    Thus,
    \begin{align*}
        \dim 
        \frac{ 
            \ker\left( \cartan^{k} : C^{\infty}(\Domain,\Gamma,\Alt{k}) \to C^{\infty}(\Domain,\Gamma,\Alt{k+1}) \right)
        }{
            \rng\left( \cartan^{k-1} : C^{\infty}(\Domain,\Gamma,\Alt{k-1}) \to C^{\infty}(\Domain,\Gamma,\Alt{k}) \right)
        } 
        = 
        b_{k}(\Domain\cup\Gamma,\Gamma)
        = 
        b_{k}(\closure\Domain,\Gamma)
        .
    \end{align*}
    The last identity holds because the pairs $(\calX,\Gamma)$ and $(\closure\Domain,\Gamma)$ are homotopy-equivalent~\cite[2.4.2, 2.4.3]{stocker_zieschang}.
    This completes the proof. 
\end{proof}

We recall the existence of a collar neighborhood along the boundary of a Lipschitz domain. 

\begin{lemma}\label{lemma:outer-collar}
    Let $\Domain \subseteq \bbR^{n}$ be a bounded Lipschitz domain. 
    Then there exists an embedding 
    \begin{align*}
        \Phi : \partial\Domain \times [0,1] \to \bbR^{n} \setminus \Domain 
    \end{align*}
    that is bi-Lipschitz onto its image and satisfies 
    \begin{align*}
        \forall x \in \partial\Domain : \Phi(x, 0) = x.
    \end{align*}
    Moreover, for each $s \in (0,1]$, the open set $\Domain \cup \Phi( \partial\Domain \times [0,s) )$ is a domain with smooth boundary $\Phi( \partial\Domain, s )$.
\end{lemma}
\begin{proof}
    We use a result by Ball and Zarnescu~\cite[Theorem~5.1, Remark~5.3]{ball2017partial};
    see also~\cite[Theorem~8.3, Theorem~8.8]{chen2019cauchy} for further discussion. 
    For every $\epsilon > 0$ small enough, 
    there exists a bi-Lipschitz map $F_{\epsilon} : \bbR^{n} \to \bbR^{n}$
    such that $F_{\epsilon}$ is a diffeomorphism on $\bbR^{n} \setminus \partial\Domain$,
    such that $\Domain_{\epsilon} := F_{\epsilon}( \Domain )$ is a bounded domain with smooth boundary,
    and $F_{\epsilon}(x) = x$ for all $x \in \bbR^{n}$ with $\dist(x,\partial\Domain) > \epsilon$. 
    Notice that $\partial\Domain_{\epsilon} = F_{\epsilon}(\partial\Domain)$ is therefore a smooth hypersurface. 
    
    Fix $\epsilon > 0$ small enough.     
    By standard results in differential geometry, there exists an embedding of manifolds with boundary 
    $\Phi_{\epsilon} : \partial\Domain_{\epsilon} \times [0,1] \to \bbR^{n} \setminus \Domain_{\epsilon}$
    that satisfies $\Phi_{\epsilon}(x, 0) = x$ for all $x \in \partial\Domain_{\epsilon}$.
    One can show that $\Phi_{\epsilon}$ is bi-Lipschitz onto its image. 
    
    We define $\Phi : \partial\Domain \times [0,1] \to \bbR^{n} \setminus \Domain$ by setting 
    \begin{align*}
        \Phi( x, t ) := F_{\epsilon}^{-1}\left( \Phi_{\epsilon}( F_{\epsilon}(x) , t) \right), \qquad x \in \partial\Domain, \quad t \in [0,1].
    \end{align*}
    This is the required transformation. 
\end{proof}

For the bounded Lipschitz domain $\Domain \subseteq \bbR^{n}$,
we fix some bi-Lipschitz embedding $\Phi$ as in Lemma~\ref{lemma:outer-collar} and define
\begin{align*}
    \collaredDomain := \Domain \cup \Phi( \partial\Domain, [0,1) ),
    \qquad 
    \Collar := \Phi( \partial\Domain, (0,1) ).
\end{align*}
We recall that $\collaredDomain$, which we henceforth refer to as the \emph{extended domain}, 
is a Lipschitz domain with boundary $\Phi( \partial\Domain,1)$,
as in the statement of the preceding lemma, 
and that $\Collar$ is a finite union of pairwise disjoint Lipschitz domains with pairwise positive distance.
Furthermore, $\partial\Collar$ consists of two disjoint compact sets, 
namely $\Phi( \partial\Domain, 0 ) = \partial\Domain$ and $\Phi( \partial\Domain, 1 )$.
We study the smooth de~Rham complex on $\Collar$ subject to boundary conditions along $\partial\Domain$. 

\begin{proposition}\label{proposition:collar-homology}
    The de~Rham complex 
    \begin{align*}
        \begin{CD}
            0 \to C^{\infty}(\Collar,\partial\Domain,\Alt{0}) 
            @>\cartan^{  0}>> 
            C^{\infty}(\Collar,\partial\Domain,\Alt{1})
            @>\cartan^{  1}>> 
            \dots 
            @>\cartan^{n-1}>> 
            C^{\infty}(\Collar,\partial\Domain,\Alt{n}) 
            \to 0
        \end{CD}
    \end{align*}
    is exact. 
    In particular, the Betti numbers of $\closure\Collar$ relative to $\partial\Domain$ vanish:
    \begin{align*}
        b_{k}(\closure\Collar,\partial\Domain) 
        = 
        \dim 
        \frac{ 
            \ker\left( \cartan^{k} : C^{\infty}(\Collar,\partial\Domain,\Alt{k}) \to C^{\infty}(\Collar,\partial\Domain,\Alt{k+1}) \right)
        }{
            \rng\left( \cartan^{k-1} : C^{\infty}(\Collar,\partial\Domain,\Alt{k-1}) \to C^{\infty}(\Collar,\partial\Domain,\Alt{k}) \right)
        } 
        = 
        0.
    \end{align*}
\end{proposition}
\begin{proof}
    Taking into account Proposition~\ref{proposition:relative-smooth-homology} for each connected component of $\Collar$, 
    it only remains for us to show $b_{k}(\closure\Collar,\partial\Domain) = 0$. 
    To that end, we define a mapping $\Theta : \closure\Collar \times [0,1] \to \closure\Collar$ by setting 
    \begin{align*}
        \Theta( \Phi(x,t), s ) := \Phi(x, s t), \qquad x \in \partial\Domain, \quad s, t \in [0,1].
    \end{align*}
    Then $\Theta$ is continuous, and for all $x \in \partial\Domain$ and $s, t \in [0,1]$ we have 
    \begin{align*}
        \Theta( \Phi(x,t), 1 ) := \Phi(x,t), 
        \qquad 
        \Theta( \Phi(x,t), 0 ) := x, 
        \qquad 
        \Theta( x, s ) := x.
    \end{align*}
    This means that $\Theta$ is a \emph{strong deformation retraction} of $\closure\Collar$ onto $\partial\Domain$,
    which means that these two topological spaces are homotopy-equivalent~\cite[2.4.1, 2.4.3]{stocker_zieschang}. 
    It then follows that the embedding of $\partial\Domain$ into $\closure\Collar$ induces an isomorphism on the singular cohomology groups, 
    and hence the Betti numbers of $\closure\Collar$ relative to $\partial\Domain$ must be trivial:
    $b_{k}( \closure\Collar, \partial\Domain ) = 0$ for all $k \in \bbZ$. 
    This shows the desired result. 
\end{proof}

This has prepared the following auxiliary lemma. It concerns the closed smooth extension of closed smooth differential forms.  

\begin{lemma}\label{lemma:closed-extension}
    Let $\Domain$ be a bounded Lipschitz domain and let $\collaredDomain$ be a collar extension with respect to a boundary collar.
    If $u \in C^{\infty}(\closure\Domain,\Alt{k})$ with $\cartan^{k} u = 0$,
    then there exists $u^{e} \in C^{\infty}(\collaredDomain,\Alt{k})$ with $\cartan^{k} u^{e} = 0$ and $u^{e} \restriction_{\Domain} = u$.
\end{lemma}
\begin{proof}
    Let $u \in C^{\infty}(\closure\Domain,\Alt{k})$ with $\cartan^{k} u = 0$.
    By Theorem~\ref{theorem:extension-componentwise}, %
    there exists $w \in C^{\infty}(\collaredDomain,\Alt{k})$ with $w \restriction_{\Domain} = u$. 
    This implies $\cartan^{k} w \restriction_{\Domain} = 0$.
    In particular, all derivatives of any order of $\cartan^{k} w$ vanish at $\partial\Domain$, which means that
    \begin{align*}
        \cartan^{k} w \in C^{\infty}( \Collar, \partial\Domain, \Alt{k+1} ).
    \end{align*}
    By Proposition~\ref{proposition:collar-homology}, there exists $v \in C^{\infty}( \Collar, \partial\Domain, \Alt{k} )$
    such that $\cartan^{k} v = \cartan^{k} w$ on $\Collar$.
    We extend $v$ trivially onto $\Domain$ and call that extension $\widetilde{v}$.
    Because the derivatives of $v$ vanish at $\partial\Domain$, 
    this gives $\widetilde{v} \in C^{\infty}( \collaredDomain, \Alt{k} )$ with $\cartan^{k} \widetilde{v} = \cartan^{k} w$.
    Define $u^{e} := w - \widetilde{v}$. By construction, 
    \begin{align*}
        u^{e} \restriction_{\Domain} = u, \qquad \cartan^{k} u^{e} = 0.
    \end{align*}
    This proves the desired statement. 
\end{proof}

\begin{corollary}\label{corollary:closed-and-bounded-extension}
    Let $\Domain \subseteq \bbR^{n}$ be a bounded Lipschitz domain.
    Then there exists a Lipschitz domain $\extDomain$ with $\Domain \Subset \extDomain$ such that the following holds.
    If $u \in C^{\infty}(\closure\Domain,\Alt{k})$ with $\cartan^{k} u = 0$,
    then there exists $u^{e} \in C^{\infty}(\closure\extDomain,\Alt{k})$ with $\cartan^{k} u^{e} = 0$ and $u^{e} \restriction_{\Domain} = u$.
\end{corollary}
\begin{proof}
    \newcommand{\extextDomain}{{\sfE_{\star}\Domain}}
    Since $\Domain$ is a Lipschitz domain, Lemma~\ref{lemma:outer-collar} gives a boundary collar $\Phi : \partial\Domain \times [0,1] \to \bbR^{n}$
    by which we extend $\Domain$ to another Lipschitz domain $\collaredDomain$. 
    Moreover, $\extDomain := \Domain \cup \Phi(\partial\Domain,[0,1/2))$ is a bounded domain with smooth boundary $\Phi( \partial\Domain, 1/2 )$ and $\Domain \Subset \extDomain \Subset \collaredDomain$. 
    If $u \in C^{\infty}(\closure\Domain,\Alt{k})$ with $\cartan^{k} u = 0$,
    then there exists $u^{e} \in C^{\infty}(\collaredDomain,\Alt{k})$ with $\cartan^{k} u^{e} = 0$ and $u^{e} \restriction_{\Domain} = u$. 
    Since $\extDomain \Subset \collaredDomain$, we have $u^{e} \restriction_{\extDomain} \in C^{\infty}(\closure{\extDomain},\Alt{k})$.
\end{proof}

\subsection{Universal commuting extension of differential forms}

After this excursion, we return to the construction of extension operators for de~Rham complexes. 
We address the second objective of the manuscript, which is the construction of a universal extension operator onto the extended domain that commutes with the exterior derivative. 
The main difference to Theorem~\ref{theorem:universal-extension} is how the extension operator acts on the harmonic complement.

\begin{theorem}\label{theorem:commuting-universal-extension}
    Let $\Domain \subseteq \bbR^{n}$ be a bounded Lipschitz domain. Then there exists a bounded Lipschitz domain $\extDomain$ with $\Domain \Subset \extDomain$ such that the following holds.
    We have linear continuous operators 
    \begin{align}\label{math:commuting-universal-extension}
        \kommErweiterung^{k} : 
        \calS'(\Domain,\Alt{k})
        \to 
        \calS'(\extDomain,\Alt{k})
        ,
        \qquad 
        k \in \bbZ
        ,
    \end{align}
    such that for all $u \in \calS'(\Domain,\Alt{k})$ we have 
    \begin{align*}
        ( \kommErweiterung^{k} u) \restriction_{\Domain} = u,
        \qquad 
        \kommErweiterung^{k+1} \cartan^{k} u = \cartan^{k} \kommErweiterung^{k} u
        .
    \end{align*}
    Moreover, we have continuous operators 
    \begin{align}\label{math:commuting-universal-extension:smooth}
        &
        \kommErweiterung^{k} : C^{\infty}(\closure\Domain,\Alt{k}) \to C^{\infty}(\closure\extDomain,\Alt{k})
        ,
        \\
        &
        \kommErweiterung^{k} : \Wres^{s,p}(\Domain,\Alt{k}) \to \Wres^{s,p}(\extDomain,\Alt{k}), 
        \qquad 
        s \in \bbR, \quad 1 < p < \infty \label{math:commuting-universal-extension:sobolev}
        ,
        \\
        &
        \kommErweiterung^{k} : B^{s}_{p,q}(\Domain,\Alt{k}) \to B^{s}_{p,q}(\extDomain,\Alt{k}),
        \qquad
        s \in \bbR, \quad 1 < p < \infty, \quad 1 < q < \infty \label{math:commuting-universal-extension:besov}
        ,
        \\
        &
        \kommErweiterung^{k} : F^{s}_{p,q}(\Domain,\Alt{k}) \to F^{s}_{p,q}(\extDomain,\Alt{k}),
        \qquad
        s \in \bbR, \quad 1 < p < \infty, \quad 1 < q < \infty \label{math:commuting-universal-extension:triebellizorkin}
        .
    \end{align}
\end{theorem}

\begin{proof}
    Suppose that $\eta_1, \dots, \eta_L \in C^{\infty}(\closure\Domain,\Alt{k})$ form a basis of the finite-dimensional harmonic complement $\Harmonic^{k}(\Domain)$.
    Using Corollary~\ref{corollary:closed-and-bounded-extension}, we can find closed extensions $\eta^{e}_{1}, \dots, \eta^{e}_{L} \in C^{\infty}(\closure\extDomain,\Alt{k})$ of the basis elements of the harmonic complement.
    By linear extension, this defines a continuous linear operator
    \begin{align*}
        \Extension_{H}^{k} : \Harmonic^{k}(\Domain) \to C^{\infty}(\closure\extDomain,\Alt{k}),
    \end{align*}
    such that 
    \begin{align*}
        \forall \eta \in \Harmonic^{k}(\Domain)
        :
        \Extension_{H}^{k} \eta \restriction_{\Domain} = \eta, 
        \qquad 
        \forall \eta \in \Harmonic^{k}(\Domain)
        :
        \cartan^{k} \Extension_{H}^{k} \eta = 0
        .
    \end{align*}
    We define the universal commuting extension operator by
    \begin{align*}
        \kommErweiterung^{k} u := 
        ( \cartan^{k-1} \Extension^{k-1} \widehat\Poincare^{k} u ) \restriction_{\extDomain}
        + 
        \Extension_{H}^{k} \Proj^{k} u
        + 
        ( \Extension^{k} \widehat\Poincare^{k+1} \cartan^{k} u ) \restriction_{\extDomain}
        .
    \end{align*}
    The continuity conditions~\eqref{math:commuting-universal-extension} and~\eqref{math:commuting-universal-extension:smooth} follow from Theorems~\ref{theorem:extension-componentwise} and~\ref{theorem:universal-regular-decomposition}.
    Throughout the proof, let $u \in \calS'(\Domain,\Alt{k})$. 
    To show the extension property, we observe 
    \begin{align*}
        ( \kommErweiterung^{k} u ) \restriction_{\Domain}
        &
        =
        ( \cartan^{k-1} \Extension^{k-1} \widehat\Poincare^{k} u ) \restriction_{\Domain}
        + 
        ( \Extension_{H}^{k} \Proj^{k} u ) \restriction_{\Domain}
        + 
        ( \Extension^{k} \widehat\Poincare^{k+1} \cartan^{k} u ) \restriction_{\Domain}
        \\&
        =
        \cartan^{k-1} ( \Extension^{k-1} \widehat\Poincare^{k} u ) \restriction_{\Domain}
        + 
        ( \Extension_{H}^{k} \Proj^{k} u ) \restriction_{\Domain}
        + 
        ( \Extension^{k} \widehat\Poincare^{k+1} \cartan^{k} u ) \restriction_{\Domain}
        \\&
        =
        \cartan^{k-1} \widehat\Poincare^{k} u
        + 
        \Proj^{k} u
        + 
        \widehat\Poincare^{k+1} \cartan^{k} u
        =
        u
        .
    \end{align*}
    The last equality is the homotopy identity in Theorem~\ref{theorem:universal-regular-decomposition}.

    To prove the commuting property, we compute
    \begin{align*}
        \cartan^{k} \kommErweiterung^{k} u
        &
        =
        \cartan^{k} \cartan^{k-1} \Extension^{k-1} \widehat\Poincare^{k} u 
        +   
        \cartan^{k} \Extension_{H}^{k} \Proj^{k} u
        +
        \cartan^{k} \Extension^{k} \widehat\Poincare^{k+1} \cartan^{k} u
        \\&
        =
        \cartan^{k} \Extension^{k} \widehat\Poincare^{k+1} \cartan^{k} u
        \\&
        =
        \cartan^{k} \Extension^{k} \widehat\Poincare^{k+1} \cartan^{k} u
        +   
        \Extension_{H}^{k+1} \Proj^{k+1} \cartan^{k} u
        +
        \Extension^{k+1} \widehat\Poincare^{k+2} \cartan^{k+1} \cartan^{k} u
        =
        \kommErweiterung^{k+1} \cartan^{k} u
        .
    \end{align*}

    It remains to prove the other continuity properties. 
    We begin with the Sobolev--Slobodeckij scale. 
    Let $s \in \bbR$ and $p \in (1,\infty)$.
    In what follows, we write $A \lesssim B$ if $A \leq C B$ holds with a constant $C \geq 0$
    that only depends on $\Domain$, $p$, and $s$. 
    
    From the definition of the extension operator, 
    \begin{align*}
        \| \kommErweiterung^{k} u \|_{\Wres^{s,p}(\extDomain,\Alt{k})}
        &
        \lesssim 
        \| 
            \Extension^{k} \widehat\Poincare^{k+1} \cartan^{k} u
            + 
            \cartan^{k-1} \Extension^{k-1} \widehat\Poincare^{k} u
            + 
            \Extension_{H}^{k} \Proj^{k} u 
        \|_{\Wres^{s,p}(\extDomain,\Alt{k})}
        \\&
        \lesssim
        \| \Extension^{k} \widehat\Poincare^{k+1} \cartan^{k} u \|_{\Wres^{s,p}(\extDomain,\Alt{k})}
        + 
        \| \cartan^{k-1} \Extension^{k-1} \widehat\Poincare^{k} u \|_{\Wres^{s,p}(\extDomain,\Alt{k})}
        +
        \| \Extension_{H}^{k} \Proj^{k} u \|_{\Wres^{s,p}(\extDomain,\Alt{k})}
        .
    \end{align*}
    We address the harmonic terms, using Theorem~\ref{theorem:extension-componentwise} and the fact that $\Proj^{k}$ has finite-dimensional range:
    \begin{align*}
        &
        \| \Extension_{H}^{k} \Proj^{k} u \|_{\Wres^{s,p}(\extDomain,\Alt{k})}
        \lesssim
        \| \Proj^{k} u \|_{\Wres^{s,p}(\Domain,\Alt{k})}
        \lesssim
        \| u \|_{\Wres^{s,p}(\Domain,\Alt{k})}
        .
    \end{align*}
    Next, we notice 
    \begin{align*}
        \| \cartan^{k-1} \Extension^{k-1} \widehat\Poincare^{k} u \|_{\Wres^{s,p}(\extDomain,\Alt{k})}
        \lesssim
        \| \Extension^{k-1} \widehat\Poincare^{k} u \|_{\Wres^{s+1,p}(\extDomain,\Alt{k-1})}
        \lesssim
        \| \widehat\Poincare^{k} u \|_{\Wres^{s+1,p}(\Domain,\Alt{k-1})}
        \lesssim
        \| u \|_{\Wres^{s,p}(\Domain,\Alt{k})}
        .
    \end{align*}
    A similar combination of techniques gives 
    \begin{align*}
        \| \Extension^{k} \widehat\Poincare^{k+1} \cartan^{k} u \|_{\Wres^{s,p}(\extDomain,\Alt{k})}
        \lesssim
        \| \widehat\Poincare^{k+1} \cartan^{k} u \|_{\Wres^{s,p}(\Domain,\Alt{k})}
        \lesssim
        \| \cartan^{k} u \|_{\Wres^{s-1,p}(\Domain,\Alt{k+1})}
        \lesssim
        \| u \|_{\Wres^{s,p}(\Domain,\Alt{k})}
        .
    \end{align*}
    This shows the desired continuity estimates for the Sobolev--Slobodeckij spaces. 
    
    The proof of the continuity estimates for Besov spaces and Triebel--Lizorkin spaces of differential forms is fully analogous.
\end{proof}

\begin{remark}
    The closed extensions of the harmonic complement via Lemma~\ref{lemma:closed-extension} rely on smooth correction terms that not only ensure that the extension is closed but also are localized in the boundary collar so as not to disturb the extension property.
\end{remark}

\begin{theorem}\label{theorem:commuting-mollification}
    Let $\Domain \subseteq \bbR^{n}$ be a bounded Lipschitz domain. 
    There exists a family of linear maps 
    \begin{align}\label{math:commuting-mollifier}
        \Mollifier^{k}_{\epsilon} : \calS'(\Domain,\Alt{k}) \to C^{\infty}(\closure\Domain,\Alt{k}), 
        \qquad 
        \epsilon > 0,
        \qquad 
        k \in \bbZ,
    \end{align}
    such that 
    \begin{align*}
        \cartan^{k} \Mollifier^{k}_{\epsilon} u = \Mollifier^{k+1}_{\epsilon} \cartan^{k} u, 
        \qquad 
        u \in \calS'(\Domain,\Alt{k}), \quad k \in \bbZ,
    \end{align*}
    and the following holds. 
    
    We have the following continuity and convergence properties:
    \begin{subequations}
        \begin{align*}
            \Mollifier^{k}_{\epsilon} : \calS'(\Domain,\Alt{k}) \to \calS'(\Domain,\Alt{k}),
            \qquad 
            u \in \calS'(\Domain,\Alt{k})
            &
            \implies
            \lim_{\epsilon \to 0} \Mollifier^{k}_{\epsilon} u = u \textrm{ in } \calS'(\Domain,\Alt{k})
            ,
            \\
            \Mollifier^{k}_{\epsilon} : C^{\infty}(\closure\Domain,\Alt{k}) \to C^{\infty}(\closure\Domain,\Alt{k}),
            \qquad 
            u \in C^{\infty}(\closure\Domain,\Alt{k})
            &
            \implies
            \lim_{\epsilon \to 0} \Mollifier^{k}_{\epsilon} u = u \textrm{ in } C^{\infty}(\closure\Domain,\Alt{k})
            .
        \end{align*}
    \end{subequations}
    
    The following boundedness and convergence conditions hold:
    \begin{subequations}
        \begin{itemize}
            \item 
            If $s \in \bbR$ and $p \in (1,\infty)$, then
            \begin{align*}
                u \in \Wres^{s,p}(\Domain,\Alt{k})
                &
                \implies
                \left\{
                \begin{array}{l} %
                \| \Mollifier^{k}_{\epsilon} u \|_{\Wres^{s,p}(\Domain,\Alt{k})} \lesssim \| u \|_{\Wres^{s,p}(\Domain,\Alt{k})}
                \\
                \lim_{\epsilon \to 0} \| \Mollifier^{k}_{\epsilon} u - u \|_{\Wres^{s,p}(\Domain,\Alt{k})} = 0
                \end{array}
                \right.
                .
            \end{align*}
            \item 
            If $s \in \bbR$ and $p,q \in (1,\infty)$, then
            \begin{align*}
                u \in B^{s}_{p,q}(\Domain,\Alt{k})
                &
                \implies
                \left\{
                \begin{array}{l} %
                \| \Mollifier^{k}_{\epsilon} u \|_{B^{s}_{p,q}(\Domain,\Alt{k})} \lesssim \| u \|_{B^{s}_{p,q}(\Domain,\Alt{k})}
                \\
                \lim_{\epsilon \to 0} \| \Mollifier^{k}_{\epsilon} u - u \|_{B^{s}_{p,q}(\Domain,\Alt{k})} = 0
                \end{array}
                \right.
                .
            \end{align*}
            \item 
            If $s \in \bbR$ and $p,q \in (1,\infty)$, then
            \begin{align*}
                u \in F^{s}_{p,q}(\Domain,\Alt{k})
                &
                \implies
                \left\{
                \begin{array}{l} %
                \| \Mollifier^{k}_{\epsilon} u \|_{F^{s}_{p,q}(\Domain,\Alt{k})} \lesssim \| u \|_{F^{s}_{p,q}(\Domain,\Alt{k})}
                \\
                \lim_{\epsilon \to 0} \| \Mollifier^{k}_{\epsilon} u - u \|_{F^{s}_{p,q}(\Domain,\Alt{k})} = 0
                \end{array}
                \right.
                .
            \end{align*}
        \end{itemize}
    \end{subequations}
    The constants implicit in $\lesssim$ are independent of $\epsilon$ and only depend on $\Domain$ and the relevant parameters $s, p, q$.
\end{theorem}

\begin{proof}
    The theorem follows from combining the standard mollification operator with the commuting universal extension operator of Theorem~\ref{theorem:commuting-universal-extension}.
    Let $\kommErweiterung^{k}$ be the universal extension operator of Theorem~\ref{theorem:commuting-universal-extension}.
    Let $\weight \in C^{\infty}_{c}(\bbR^{n},\bbR)$ be a standard mollifier, that is, $\weight \geq 0$, $\int_{\bbR^{n}} \weight = 1$, and $\supp(\weight) \subseteq \Ball{0}{1}$.
    For any $\epsilon > 0$ and $x \in \bbR^{n}$, we define the scaled mollifier $\weight_{\epsilon}(x) := \epsilon^{-n} \weight(x/\epsilon)$.
    We define the global smoothing operator by 
    \begin{align*}
        \euS^{k}_{\epsilon} u := \weight_{\epsilon} \star u,
        \qquad 
        u \in \calS'(\bbR^{n},\Alt{k})
        .
    \end{align*}
    We notice that $\cartan^{k} \euS^{k}_{\epsilon} u = \euS^{k+1}_{\epsilon} \cartan^{k} u$ for all $u \in \calS'(\bbR^{n},\Alt{k})$.
    
    Write $\delta := \operatorname{dist}( \closure\Domain, \bbR^{n} \setminus \extDomain)$ and assume that $\epsilon \in (0,\delta)$. 
    Then the locality of the global smoothing operator shows that we have a well-defined operator
    \begin{align*}
        \euS^{k}_{\Domain,\epsilon} : \calS'(\extDomain,\Alt{k}) \to \calS'(\Domain,\Alt{k})
        .
    \end{align*}
    We define the mollification operator by applying the global smoothing to the universal extension and then restricting back to the domain:
    \begin{align*}
        \Mollifier^{k}_{\epsilon} u := ( \euS^{k}_{\Domain,\epsilon} \circ \kommErweiterung^{k} u ) %
        ,
        \qquad 
        u \in \calS'(\Domain,\Alt{k})
        .
    \end{align*}
    All claimed properties of the mollification operator now follow from Theorem~\ref{theorem:commuting-universal-extension} and the standard properties of convolution with a mollifier.
    This constructs the operator for $\epsilon > 0$ small enough, which gives the statement of the theorem.
\end{proof}

\begin{remark}
    Theorem~\ref{theorem:commuting-mollification} is restricted to $p, q < \infty$ for simplicity of exposition. 
    We refer to the literature~\cite[Section~2.1.3]{runstsickel} for background on the general case.
\end{remark}

The main application is the following. On any open set $\Domain \subseteq \bbR^{n}$, we define 
\begin{align*}
    \Sobdiff^{s,t,p}(\Domain,\Alt{k})
    := 
    \left\{ u \in \Wres^{s,p}(\Domain,\Alt{k}) \suchthat* \cartan u \in \Wres^{t,p}(\Domain,\Alt{k+1}) \right\}
    ,
    \qquad 
    p \in (1,\infty), 
    \quad 
    s \in \bbR,
    \quad 
    t \in [s-1,\infty)
    .
\end{align*}
This is a Banach space with the norm 
\begin{align*}
    \| u \|_{\Sobdiff^{s,t,p}(\Domain,\Alt{k})}
    :=
    \| u \|_{\Wres^{s,p}(\Domain,\Alt{k})}
    +
    \| \cartan u \|_{\Wres^{t,p}(\Domain,\Alt{k+1})}
    .
\end{align*}
It describes Sobolev differential forms whose exterior derivative shows additional smoothness. 
If $\Domain \subseteq \bbR^{3}$ is a three-dimensional domain, then we set 
\begin{align*}
    \Sobdiff^{s,t,p}(\Domain,\curl)
    &
    := 
    \left\{ u \in \Wres^{s,p}(\Domain,\bbR^{3}) \suchthat* \curl u \in \Wres^{t,p}(\Domain,\bbR^{3}) \right\}
    ,
    \\
    \Sobdiff^{s,t,p}(\Domain,\divergence)
    &
    := 
    \left\{ u \in \Wres^{s,p}(\Domain,\bbR^{3}) \suchthat* \divergence u \in \Wres^{t,p}(\Domain,\bbR) \right\}
    ,
\end{align*}
where $p \in (1,\infty)$ and $s, t \in \bbR$ with $t \geq s-1$.
These spaces of three-dimensional vector fields are equipped with the obvious norms
and are special cases of $\Sobdiff^{s,t,p}(\Domain,\Alt{k})$. 

The results in this manuscript allow the following smooth approximation.
We use Theorem~\ref{theorem:commuting-mollification} but the result can also be derived similarly by Theorem~\ref{theorem:universal-extension}. 

\begin{proposition}\label{proposition:smooth-forms-are-dense}
    Let $\Domain \subseteq \bbR^{n}$ be a bounded Lipschitz domain.
    Let $p \in (1,\infty)$ and $s, t \in \bbR$ with $t \geq s-1$. 
    Then we have a dense inclusion 
    \begin{align*}
        C^{\infty}(\closure\Domain,\Alt{k}) \subseteq \Sobdiff^{s,t,p}(\Domain,\Alt{k})
        .
    \end{align*}
    If $n = 3$, then we have dense inclusions 
    \begin{align*}
        C^{\infty}(\closure\Domain,\bbR^{3}) \subseteq \Sobdiff^{s,t,p}(\Domain,\curl)
        ,
        \qquad 
        C^{\infty}(\closure\Domain,\bbR^{3}) \subseteq \Sobdiff^{s,t,p}(\Domain,\divergence)
        .
    \end{align*}
\end{proposition}
\begin{proof}
    This immediately follows from Theorem~\ref{theorem:commuting-mollification}.
\end{proof}

\section{Examples in vector calculus}\label{section:vector-calculus}

We present the main results of this manuscript in the language of vector calculus on domains, and for the special case of Sobolev--Slobodeckij spaces, defined in terms of weak derivatives and Gagliardo seminorms. 

Throughout this section, $\Domain \subseteq \bbR^{3}$ denotes a bounded Lipschitz domain. 
Recall the Lebesgue spaces of scalar fields $L_{p}(\Domain) := L^{p}(\Domain,\bbR)$ and vector fields $\bfL_{p}(\Domain) := L^{p}(\Domain,\bbR^{3})$.
For $\sigma \geq 0$, we define the corresponding Sobolev--Slobodeckij spaces 
\begin{align*}
    W^{\sigma}_{p}(\Domain) := \Wint^{\sigma,p}(\Domain,\bbR),
    \qquad 
    \bfW^{\sigma}_{p}(\Domain) := \Wint^{\sigma,p}(\Domain,\bbR^{3})
    .
\end{align*}
Spaces of vector fields with non-standard smoothness are an intrinsic feature of vector calculus. 
Given $1 < p < \infty$ and $\sigma, \tau \in \bbR_{0}^{+}$ with $\tau \geq \sigma-1$, we define the spaces 
\begin{align*}
    \bfW^{\sigma,\tau}_{p}(\Domain,\curl)
    &
    :=
    \left\{ u \in \bfW^{\sigma}_{p}(\Domain) \suchthat* \curl u \in \bfW^{\tau}_{p}(\Domain) \right\},
    \\
    \bfW^{\sigma,\tau}_{p}(\Domain,\divergence)
    &
    :=
    \left\{ u \in \bfW^{\sigma}_{p}(\Domain) \suchthat* \divergence u \in W^{\tau}_{p}(\Domain) \right\}
    .
\end{align*}
These spaces are equipped with the obvious norms
\begin{align*}
    \| u \|_{\bfW^{\sigma,\tau}_{p}(\Domain,\curl)}
    :=
    \| u \|_{\bfW^{\sigma}_{p}(\Domain)}
    +
    \| \curl u \|_{\bfW^{\tau}_{p}(\Domain)}
    ,
    \qquad
    \| u \|_{\bfW^{\sigma,\tau}_{p}(\Domain,\divergence)}
    :=
    \| u \|_{\bfW^{\sigma}_{p}(\Domain)}
    +
    \| \divergence u \|_{W^{\tau}_{p}(\Domain)}
    .   
\end{align*}
Moreover, these spaces constitute a differential complex. If $1 < p < \infty$ and $\rho, \sigma, \tau \in \bbR_{0}^{+}$ with $\sigma \geq \rho-1$ and $\tau \geq \sigma-1$, then this differential complex takes the form 
\begin{align*}
    \begin{CD}
        W^{\rho+1}_{p}(\Domain) @> \grad >> \bfW^{\rho,\sigma}_{p}(\Domain,\curl) @> \curl >> \bfW^{\sigma,\tau}_{p}(\Domain,\divergence) @> \divergence >> W^{\tau}_{p}(\Domain)
        .
    \end{CD}
\end{align*}
As a common base space that contains all the preceding classes of scalar and vector fields, we define the restriction spaces of tempered distributions with values in scalars and vectors. This ultimately simplifies the exposition. We define 
\begin{align*}
    \calS'(\Domain) := \calS'(\Domain,\bbR),
    \qquad 
    \calS'(\Domain)^{3} := \calS'(\Domain,\bbR^{3})
    .
\end{align*}
Accordingly, we have a differential complex
\begin{align*}
    \begin{CD}
        \calS'(\Domain) @> \grad >> \calS'(\Domain)^{3} @> \curl >> \calS'(\Domain)^{3} @> \divergence >> \calS'(\Domain)
        .
    \end{CD}
\end{align*}

Earlier in this manuscript, we have reviewed the potential operators of Costabel and McIntosh. 
For completeness of exposition, we summarize the properties of these operators explicitly, including their behavior on spaces of vector fields.

\begin{theorem}
    Let $\Domain \subseteq \bbR^{3}$ be a bounded Lipschitz domain.
    There exist linear continuous operators
    \begin{align*}
        \Poincare^{\grad} : \calS'(\Domain)^{3} \to \calS'(\Domain)
        ,
        \qquad
        \Poincare^{\curl} : \calS'(\Domain)^{3} \to \calS'(\Domain)^{3}
        ,
        \qquad
        \Poincare^{\divergence} : \calS'(\Domain) \to \calS'(\Domain)^{3}
        ,
    \end{align*}
    and continuous operators with finite-dimensional range 
    \begin{gather*}
        \Defect^{\grad} : \calS'(\Domain) \to C^{\infty}(\Domain)
        ,
        \qquad
        \Defect^{\curl} : \calS'(\Domain)^{3} \to C^{\infty}(\Domain)^{3}
        ,
        \qquad 
        \Defect^{\divergence} : \calS'(\Domain)^{3} \to C^{\infty}(\Domain)^{3}
        .
    \end{gather*}
    In addition, the following homotopy identities hold:
    \begin{align*}
        \forall u \in \calS'(\Domain) 
        &:
        u = \Poincare^{\grad} \grad u + \Defect^{\grad} u,
        \\
        \forall u \in \calS'(\Domain)^{3} 
        &:
        u = \Poincare^{\curl} \curl u + \grad \Poincare^{\grad} u + \Defect^{\curl} u,
        \\
        \forall u \in \calS'(\Domain)^{3} 
        &:
        u = \Poincare^{\divergence} \divergence u + \curl \Poincare^{\curl} u + \Defect^{\divergence} u,
        \\
        \forall u \in \calS'(\Domain) 
        &:
        u = \divergence \Poincare^{\divergence} u
        .
    \end{align*}
    For any $1 < p,q < \infty$ and $\sigma, \tau \geq 0$, the following operators are compact:
    \begin{gather*}
        \Defect^{\grad} : W^{\sigma}_{p}(\Domain) \to W^{\tau}_{q}(\Domain)
        ,
        \qquad 
        \Defect^{\curl} : \bfW^{\sigma}_{p}(\Domain) \to \bfW^{\tau}_{q}(\Domain)
        ,
        \qquad  
        \Defect^{\divergence} : \bfW^{\sigma}_{p}(\Domain) \to \bfW^{\tau}_{q}(\Domain)
        .
    \end{gather*}
    In addition, if $1 < p < \infty$ and $\sigma, \tau \in \bbR_{0}^{+}$ with $\tau \geq \sigma-1$, then we have the following continuity properties:
    \begin{align*}
        \Poincare^{\grad} &: \bfW^{\sigma,\tau}_{p}(\Domain,\curl) \to W^{\sigma+1}_{p}(\Domain),
        \\
        \Poincare^{\curl} &: \bfW^{\sigma,\tau}_{p}(\Domain,\divergence) \to \bfW^{\sigma+1,\sigma}_{p}(\Domain,\curl),
        \\
        \Poincare^{\divergence} &: W^{\sigma}_{p}(\Domain) \to \bfW^{\sigma+1,\sigma}_{p}(\Domain,\divergence)
        .
    \end{align*}
\end{theorem}

We have introduced the commuting extension operators.
These have the following specialization to vector calculus. 

\begin{theorem}
    Let $\Domain \subseteq \bbR^{3}$ be a bounded Lipschitz domain. 
    Then there exists a bounded Lipschitz domain $\extDomain$ with $\Domain \Subset \extDomain$ such that the following holds.
    We have linear continuous operators 
    \begin{gather*}
        \kommErweiterung^{\grad} : \calS'(\Domain) \to \calS'(\extDomain),
        \qquad
        \kommErweiterung^{\curl} : \calS'(\Domain)^{3} \to \calS'(\extDomain)^{3},
        \\
        \kommErweiterung^{\divergence} : \calS'(\Domain)^{3} \to \calS'(\extDomain)^{3},
        \qquad 
        \kommErweiterung^{0} : \calS'(\Domain) \to \calS'(\extDomain)
        ,
    \end{gather*}
    such that for all $u$ in the respective spaces we have 
    \begin{align*}
        ( \kommErweiterung^{\grad} u) \restriction_{\Domain} = u,
        \qquad 
        ( \kommErweiterung^{\curl} u) \restriction_{\Domain} = u,
        \qquad 
        ( \kommErweiterung^{\divergence} u) \restriction_{\Domain} = u,
        \qquad 
        ( \kommErweiterung^{0} u) \restriction_{\Domain} = u,
    \end{align*}
    and we have 
    \begin{align*}
        \forall u \in \calS'(\Domain)
        &: \kommErweiterung^{\curl} (\grad u) = \grad (\kommErweiterung^{\grad} u),
        \\
        \forall u \in \calS'(\Domain)^{3}
        &: \kommErweiterung^{\divergence} (\curl u) = \curl (\kommErweiterung^{\curl} u),
        \\
        \forall u \in \calS'(\Domain)^{3}
        &: \kommErweiterung^{0} (\divergence u) = \divergence (\kommErweiterung^{\divergence} u).
    \end{align*}
    In addition, if $1 < p < \infty$ and $\sigma, \tau \in \bbR_{0}^{+}$ with $\tau \geq \sigma-1$, then we have continuous operators
    \begin{align*}
        \kommErweiterung^{\grad} &: W^{\sigma}_{p}(\Domain) \to W^{\sigma}_{p}(\extDomain),
        \\
        \kommErweiterung^{\curl} &: \bfW^{\sigma,\tau}_{p}(\Domain,\curl) \to \bfW^{\sigma,\tau}_{p}(\extDomain,\curl),
        \\
        \kommErweiterung^{\divergence} &: \bfW^{\sigma,\tau}_{p}(\Domain,\divergence) \to \bfW^{\sigma,\tau}_{p}(\extDomain,\divergence),
        \\
        \kommErweiterung^{0} &: W^{\sigma}_{p}(\Domain) \to W^{\sigma}_{p}(\extDomain)
        .
    \end{align*}
\end{theorem}

As a result, we have commuting mollification operators for these de~Rham complexes. 
These smoothing operators converge to the identity in the strong operator topology, which shows that smooth fields are dense in the spaces of non-standard regularity. 

\begin{theorem}
    Let $\Domain \subseteq \bbR^{3}$ be a bounded Lipschitz domain. 
    There exists a family of linear maps 
    \begin{gather*}
        \Mollifier^{\grad}_{\epsilon} : \calS'(\Domain) \to C^{\infty}(\closure\Domain),
        \qquad
        \Mollifier^{\curl}_{\epsilon} : \calS'(\Domain)^{3} \to C^{\infty}(\closure\Domain)^{3},
        \\
        \Mollifier^{\divergence}_{\epsilon} : \calS'(\Domain)^{3} \to C^{\infty}(\closure\Domain)^{3},
        \qquad
        \Mollifier^{0}_{\epsilon} : \calS'(\Domain) \to C^{\infty}(\closure\Domain),
    \end{gather*}
    indexed over $\epsilon > 0$, 
    such that the following commuting properties hold:
    \begin{align*}
        \forall u \in \calS'(\Domain)
        &: \grad (\Mollifier^{\grad}_{\epsilon} u) = \Mollifier^{\curl}_{\epsilon} (\grad u),
        \\
        \forall u \in \calS'(\Domain)^{3}
        &: \curl (\Mollifier^{\curl}_{\epsilon} u) = \Mollifier^{\divergence}_{\epsilon} (\curl u),
        \\
        \forall u \in \calS'(\Domain)^{3}
        &: \divergence (\Mollifier^{\divergence}_{\epsilon} u) = \Mollifier^{0}_{\epsilon} (\divergence u).
    \end{align*}
    If $1 < p < \infty$ and $\sigma, \tau \in \bbR_{0}^{+}$ with $\tau \geq \sigma-1$, then we have continuous operators
    \begin{align*}
        \Mollifier^{\grad}_{\epsilon} &: W^{\sigma}_{p}(\Domain) \to W^{\sigma}_{p}(\Domain),
        \\
        \Mollifier^{\curl}_{\epsilon} &: \bfW^{\sigma,\tau}_{p}(\Domain,\curl) \to \bfW^{\sigma,\tau}_{p}(\Domain,\curl),
        \\
        \Mollifier^{\divergence}_{\epsilon} &: \bfW^{\sigma,\tau}_{p}(\Domain,\divergence) \to \bfW^{\sigma,\tau}_{p}(\Domain,\divergence),
        \\
        \Mollifier^{0}_{\epsilon} &: W^{\sigma}_{p}(\Domain) \to W^{\sigma}_{p}(\Domain)
        .
    \end{align*}
    If $1 < p < \infty$ and $\sigma, \tau \in \bbR_{0}^{+}$ with $\tau \geq \sigma-1$, then we have dense inclusions
    \begin{align*}
        C^{\infty}(\closure\Domain) \subseteq W^{\sigma}_{p}(\Domain),
        \qquad
        C^{\infty}(\closure\Domain)^{3} \subseteq \bfW^{\sigma,\tau}_{p}(\Domain,\curl),
        \qquad
        C^{\infty}(\closure\Domain)^{3} \subseteq \bfW^{\sigma,\tau}_{p}(\Domain,\divergence)
        ,
    \end{align*}
    and the mollification operators converge strongly to the identity in the respective norms as $\epsilon$ goes to zero.
\end{theorem}

Lastly, we point out a reformulation in vector calculus of our extension result on smooth forms.

\begin{proposition}
    Let $\Domain \subseteq \bbR^{3}$ be a bounded Lipschitz domain.
    Then there exists a Lipschitz domain $\extDomain$ with $\Domain \Subset \extDomain$ such that the following holds.
    \begin{itemize}
        \item 
        If $u \in C^{\infty}(\closure\Domain)^{3}$ with $\curl u = 0$, 
        then there exists $u^{e} \in C^{\infty}(\closure\extDomain)^{3}$ with $\curl u^{e} = 0$ and $u^{e} \restriction_{\Domain} = u$.
        \item 
        If $u \in C^{\infty}(\closure\Domain)^{3}$ with $\divergence u = 0$, 
        then there exists $u^{e} \in C^{\infty}(\closure\extDomain)^{3}$ with $\divergence u^{e} = 0$ and $u^{e} \restriction_{\Domain} = u$.
    \end{itemize}
\end{proposition}

\bibliographystyle{plain}%

\begin{thebibliography}{10}

\bibitem{agricola2010vektoranalysis}
Ilka Agricola and Thomas Friedrich.
\newblock {\em Vektoranalysis. {Differentialformen} in {Analysis}, {Geometrie}
  und {Physik}}.
\newblock Wiesbaden: Vieweg+Teubner, 2nd revised and enlarged edition, 2010.

\bibitem{amann2009analysisiii}
Herbert Amann and Joachim Escher.
\newblock {\em Analysis {III}}.
\newblock Basel: Birkh{\"a}user, 2009.
\newblock Translated from the German by Silvio Levy and Matthew Cargo.

\bibitem{amrouche1998vector}
C.~Amrouche, C.~Bernardi, M.~Dauge, and V.~Girault.
\newblock Vector potentials in three-dimensional non-smooth domains.
\newblock {\em Math. Methods Appl. Sci.}, 21(9):823--864, 1998.

\bibitem{amrouche1994decomposition}
Ch{\'e}rif Amrouche and Vivette Girault.
\newblock Decomposition of vector spaces and application to the {Stokes}
  problem in arbitrary dimension.
\newblock {\em Czech. Math. J.}, 44(1):109--140, 1994.

\bibitem{AFW1}
Douglas~N. Arnold, Richard~S. Falk, and Ragnar Winther.
\newblock Finite element exterior calculus, homological techniques, and
  applications.
\newblock {\em Acta Numerica}, 15:1--155, 2006.

\bibitem{baldi2019l1}
Annalisa Baldi, Bruno Franchi, and Pierre Pansu.
\newblock {{\(L^1\)}}-poincar{\'e} and {Sobolev} inequalities for differential
  forms in {Euclidean} spaces.
\newblock {\em Sci. China, Math.}, 62(6):1029--1040, 2019.

\bibitem{ball2017partial}
John~M. Ball and Arghir Zarnescu.
\newblock Partial regularity and smooth topology-preserving approximations of
  rough domains.
\newblock {\em Calc. Var. Partial Differ. Equ.}, 56(1), 2017.
\newblock Article 13, 32 pp.

\bibitem{bartocci1991geometry}
Claudio Bartocci, Ugo Bruzzo, and Daniel Hern{\'a}ndez-Ruip{\'e}rez.
\newblock {\em The geometry of supermanifolds}, volume~71 of {\em Math. Appl.,
  Dordr.}
\newblock Dordrecht etc.: Kluwer Academic Publishers, 1991.

\bibitem{bauer2016maxwell}
Sebastian Bauer, Dirk Pauly, and Michael Schomburg.
\newblock The {Maxwell} compactness property in bounded weak {Lipschitz}
  domains with mixed boundary conditions.
\newblock {\em SIAM J. Math. Anal.}, 48(4):2912--2943, 2016.

\bibitem{theeigenbook}
Jo{\"e}l Bella{\"{\i}}che.
\newblock {\em The {Eigenbook}. {Eigenvarieties}, families of {Galois}
  representations, {{\(p\)}}-adic {{\(L\)}}-functions}.
\newblock Pathw. Math. Cham: Birkh{\"a}user, 2021.

\bibitem{birman1987theory}
M.~Sh. Birman and M.~Z. Solomyak.
\newblock {{\(L_ 2\)}}-theory of the {Maxwell} operator in arbitrary domains.
\newblock {\em Russ. Math. Surv.}, 42(6):75--96, 1987.

\bibitem{bossavit1992new}
A.~Bossavit.
\newblock A new viewpoint on mixed elements.
\newblock {\em Meccanica}, 27(1):3--11, 1992.

\bibitem{bourgain2007new}
Jean Bourgain and Ha{\"{\i}}m Brezis.
\newblock New estimates for elliptic equations and {Hodge} type systems.
\newblock {\em J. Eur. Math. Soc. (JEMS)}, 9(2):277--315, 2007.

\bibitem{calderon1961lebesgue}
Alberto~P. Calder{\'o}n.
\newblock Lebesgue spaces of differentiable functions and distributions.
\newblock Proc. {Sympos}. {Pure} {Math}. 4, 33--49, 1961.

\bibitem{cap1998garben}
Andreas {\v{C}}ap.
\newblock {Garben und Garbenkohomologie}, 1998.
\newblock Seminarskript, Sommersemester 1998.

\bibitem{cartan1967formes}
Henri Cartan.
\newblock {\em Formes diff{\'e}rentielles. {Applications} {\'e}l{\'e}mentaires
  au calcul des variations et {\`a} la th{\'e}orie des courbes et des
  surfaces}.
\newblock Collection m{\'e}thodes. Hermann \& Cie, Paris, 1967.
\newblock 185 pp.

\bibitem{chen2019cauchy}
Gui-Qiang~G. Chen, Giovanni~E. Comi, and Monica Torres.
\newblock Cauchy fluxes and {Gauss}-{Green} formulas for divergence-measure
  fields over general open sets.
\newblock {\em Arch. Ration. Mech. Anal.}, 233(1):87--166, 2019.

\bibitem{christ1984extension}
M.~Christ.
\newblock The extension problem for certain function spaces involving
  fractional orders of differentiability.
\newblock {\em Ark. Mat.}, 22:63--81, 1984.

\bibitem{costabel2010}
Martin Costabel and Alan McIntosh.
\newblock On {B}ogovski\u\i\ and regularized {P}oincar\'e integral operators
  for de {R}ham complexes on {L}ipschitz domains.
\newblock {\em Math. Z.}, 265(2):297--320, 2010.

\bibitem{costabel2013potential}
Martin Costabel, Alan McIntosh, and Robert~J. Taggart.
\newblock Potential maps, {Hardy} spaces, and tent spaces on special
  {Lipschitz} domains.
\newblock {\em Publ. Mat., Barc.}, 57(2):295--331, 2013.

\bibitem{derham1984differentiable}
Georges de~Rham.
\newblock {\em Differentiable manifolds. {Forms}, currents, harmonic forms},
  volume 266 of {\em Grundlehren Math. Wiss.}
\newblock Springer, Berlin, Heidelberg, 1984.
\newblock Translated from the French by F. R. Smith. Introduction to the
  English edition by S. S. Chern.

\bibitem{deitmar2025kohomologie}
Anton Deitmar.
\newblock {Kohomologie und Garben}.
\newblock Vorlesungsskript, Universit{\"a}t T{\"u}bingen, 2025.
\newblock Sommersemester 2025.

\bibitem{demlow2014posteriori}
Alan Demlow and Anil~N. Hirani.
\newblock A posteriori error estimates for finite element exterior calculus:
  the de {R}ham complex.
\newblock {\em Found. Comput. Math.}, 14(6):1337--1371, 2014.

\bibitem{edwards2014advanced}
Harold~M. Edwards.
\newblock {\em Advanced calculus. {A} differential forms approach}.
\newblock Mod. Birkh{\"a}user Class. New York, NY: Birkh{\"a}user/Springer,
  2014.
\newblock Reprint of the 1994 edition.

\bibitem{ern2021finite}
Alexandre Ern and Jean-Luc Guermond.
\newblock {\em Finite {E}lements {I}. {A}pproximation and {I}nterpolation},
  volume~72 of {\em Texts in Applied Mathematics}.
\newblock Springer International Publishing, Springer Nature Switzerland AG,
  2021.

\bibitem{federer}
Herbert Federer.
\newblock {\em Geometric measure theory.}
\newblock Class. Math. Berlin: Springer-Verlag, 1996.
\newblock Reprint of the 1969 edition.

\bibitem{goldshtein2026some}
V.~Gol'dshtein, Ya. Kopylov, and R.~Panenko.
\newblock Some {Poincar{\'e}}-{Sobolev} inequalities for differential forms.
\newblock {\em J. Math. Sci., New York}, 298(4):532--549, 2026.

\bibitem{goldshtein2011hodge}
V.~Gol'dshtein, I.~Mitrea, and M.~Mitrea.
\newblock Hodge decompositions with mixed boundary conditions and applications
  to partial differential equations on {Lipschitz} manifolds.
\newblock {\em J. Math. Sci., New York}, 172(3):347--400, 2011.

\bibitem{goldshtein1982differential}
V.~M. Gol'dshtejn, V.~I. Kuz'minov, and I.~A. Shvedov.
\newblock Differential forms on {Lipschitz} manifolds.
\newblock {\em Sib. Math. J.}, 23:151--161, 1982.

\bibitem{gopalakrishnan2018auxiliary}
J.~Gopalakrishnan, M.~Neum{\"u}ller, and P.~S. Vassilevski.
\newblock The auxiliary space preconditioner for the de {Rham} complex.
\newblock {\em SIAM J. Numer. Anal.}, 56(6):3196--3218, 2018.

\bibitem{gopalakrishnan2012partial}
Jay Gopalakrishnan and Weifeng Qiu.
\newblock Partial expansion of a {Lipschitz} domain and some applications.
\newblock {\em Front. Math. China}, 7(2):249--272, 2012.

\bibitem{grisvard2011elliptic}
Pierre Grisvard.
\newblock {\em Elliptic problems in nonsmooth domains}, volume~69 of {\em
  Class. Appl. Math.}
\newblock Philadelphia, PA: Society for Industrial {and} Applied Mathematics
  (SIAM), 2011.
\newblock Reprint of the 1985 hardback edition.

\bibitem{grubb}
Gerd Grubb.
\newblock {\em Distributions and operators}, volume 252 of {\em Grad. Texts
  Math.}
\newblock New York, NY: Springer, 2009.

\bibitem{hiptmair2002finite}
R.~Hiptmair.
\newblock Finite elements in computational electromagnetism.
\newblock {\em Acta Numer.}, 11:237--339, 2002.

\bibitem{hiptmair2012universal}
Ralf Hiptmair, Jingzhi Li, and Jun Zou.
\newblock Universal extension for {Sobolev} spaces of differential forms and
  applications.
\newblock {\em J. Funct. Anal.}, 263(2):364--382, 2012.

\bibitem{hiptmair2019discrete}
Ralf Hiptmair and Clemens Pechstein.
\newblock Discrete regular decompositions of tetrahedral discrete 1-forms.
\newblock In {\em Maxwell's equations. Analysis and numerics. Contributions
  from the workshop on analysis and numerics of acoustic and electromagnetic
  problems, RICAM, Linz, Austria, October 17--22, 2016}, pages 199--258.
  Berlin: De Gruyter, 2019.

\bibitem{hiptmair2020review}
Ralf Hiptmair and Clemens Pechstein.
\newblock A review of regular decompositions of vector fields: Continuous,
  discrete and structure-preserving.
\newblock In {\em Spectral and High Order Methods for Partial Differential
  Equations ICOSAHOM 2018: Selected Papers from the ICOSAHOM Conference,
  London, UK, July 9-13, 2018}, pages 45--60. Springer Nature, 2020.

\bibitem{hiptmair2007nodal}
Ralf Hiptmair and Jinchao Xu.
\newblock Nodal auxiliary space preconditioning in {H}(curl) and {H}(div)
  spaces.
\newblock {\em SIAM J. Numer. Anal.}, 45(6):2483--2509, 2007.

\bibitem{hoermander}
Lars H{\"o}rmander.
\newblock {\em The analysis of linear partial differential operators. {III}:
  {Pseudo}-differential operators}.
\newblock Class. Math. Berlin: Springer, 2007.
\newblock Reprint of the 1994 edition.

\bibitem{iwaniec1999nonlinear}
T.~Iwaniec, C.~Scott, and B.~Stroffolini.
\newblock Nonlinear {H}odge theory on manifolds with boundary.
\newblock {\em Annali di Matematica Pura ed Applicata}, 177(1):37--115, 1999.

\bibitem{iwaniec1993integral}
Tadeusz Iwaniec and Adam Lutoborski.
\newblock Integral estimates for null {Lagrangians}.
\newblock {\em Arch. Ration. Mech. Anal.}, 125(1):25--79, 1993.

\bibitem{jones1981quasiconformal}
Peter~W. Jones.
\newblock Quasiconformal mappings and extendability of functions in {Sobolev}
  spaces.
\newblock {\em Acta Math.}, 147:71--88, 1981.

\bibitem{kalton2007interpolation}
Nigel Kalton, Svitlana Mayboroda, and Marius Mitrea.
\newblock Interpolation of {Hardy}-{Sobolev}-{Besov}-{Triebel}-{Lizorkin}
  spaces and applications to problems in partial differential equations.
\newblock In {\em Interpolation theory and applications. A conference in honor
  of Michael Cwikel on the occasion of his 59th birthday, March 29--31, 2006
  and AMS special session on interpolation theory and applications, AMS
  sectional meeting, Miami, FL, USA, April 1--2, 2006}, pages 121--177.
  Providence, RI: American Mathematical Society (AMS), 2007.

\bibitem{khelemskij1989homology}
A.~Ya. Khelemskij.
\newblock {\em The homology of {Banach} and topological algebras. {Transl}.
  from the {Russian} by {Alan} {West}}, volume~41 of {\em Math. Appl., Sov.
  Ser.}
\newblock Dordrecht etc.: Kluwer Academic Publishers, 1989.

\bibitem{koethe2}
Gottfried K{\"o}the.
\newblock {\em Topological vector spaces. {II}}, volume 237 of {\em Grundlehren
  Math. Wiss.}
\newblock Springer, New York, 1979.

\bibitem{kuzminov1998homological}
V.~I. Kuz'minov and I.~A. Shvedov.
\newblock Homological aspects of the theory of {Banach} complexes.
\newblock {\em Sib. Math. J.}, 40(4):893--904, 1998.

\bibitem{kuzminov2003compactness}
V.~I. Kuz'minov and I.~A. Shvedov.
\newblock On the compactness theorem for differential forms.
\newblock {\em Sib. Mat. Zh.}, 44(1):132--142, 2003.

\bibitem{lee2012smooth}
John~M. Lee.
\newblock {\em Introduction to smooth manifolds}, volume 218 of {\em Graduate
  Texts in Mathematics}.
\newblock Springer, New York, second edition, 2013.

\bibitem{lewicka2023calculus}
Marta Lewicka.
\newblock {\em Calculus of variations on thin prestressed films. {Asymptotic}
  methods in elasticity}, volume 101 of {\em Prog. Nonlinear Differ. Equ.
  Appl.}
\newblock Cham: Birkh{\"a}user, 2023.

\bibitem{licht2019mixed}
Martin~W. Licht.
\newblock Smoothed projections and mixed boundary conditions.
\newblock {\em Math. Comput.}, 88(316):607--635, 2019.

\bibitem{licht2019smoothed}
Martin~W. Licht.
\newblock Smoothed projections over weakly {Lipschitz} domains.
\newblock {\em Math. Comput.}, 88(315):179--210, 2019.

\bibitem{licht2023averaging}
Martin~W. Licht.
\newblock Averaging-based local projections in finite element exterior
  calculus.
\newblock arXiv preprint arXiv:2301.03007, 2023.

\bibitem{lichtenstein1929elementare}
L.~Lichtenstein.
\newblock Eine elementare {Bemerkung} zur reellen {Analysis}.
\newblock {\em Math. Z.}, 30:794--795, 1929.

\bibitem{marschall1996nonregular}
J.~Marschall.
\newblock Nonregular pseudo-differential operators.
\newblock {\em Z. Anal. Anwend.}, 15(1):109--148, 1996.

\bibitem{melzak1959class}
Z.~A. Melzak.
\newblock A class of star-shaped bodies.
\newblock {\em Can. Math. Bull.}, 2:175--180, 1959.

\bibitem{mitrea2008poisson}
Dorina Mitrea, Marius Mitrea, and Sylvie Monniaux.
\newblock The {P}oisson problem for the exterior derivative operator with
  {D}irichlet boundary condition in nonsmooth domains.
\newblock {\em Commun. Pure Appl. Anal.}, 7(6):1295--1333, 2008.

\bibitem{mitrea2008traces}
Dorina Mitrea, Marius Mitrea, and Mei-Chi Shaw.
\newblock Traces of differential forms on {Lipschitz} domains, the boundary de
  {Rham} complex and {Hodge} decompositions.
\newblock {\em Indiana Univ. Math. J.}, 57(5):2061--2095, 2008.

\bibitem{mitrea2004sharp}
Marius Mitrea.
\newblock Sharp {Hodge} decompositions, {Maxwell}'s equations, and vector
  {Poisson} problems on nonsmooth, three-dimensional {Riemannian} manifolds.
\newblock {\em Duke Math. J.}, 125(3):467--547, 2004.

\bibitem{nilsson2026uniformly}
Erik Nilsson and Silvano Pitassi.
\newblock Uniformly bounded cochain extensions and uniform {Poincar\'e}
  inequalities.
\newblock {\em arXiv preprint arXiv:2604.04927}, 2026.

\bibitem{palamodov1972homological}
V.~P. Palamodov.
\newblock Homological methods in the theory of locally convex spaces.
\newblock {\em Russ. Math. Surv.}, 26(1):1--64, 1972.

\bibitem{pauly2022hilbert}
Dirk Pauly and Michael Schomburg.
\newblock Hilbert complexes with mixed boundary conditions. {I}: {De} {Rham}
  complex.
\newblock {\em Math. Methods Appl. Sci.}, 45(5):2465--2507, 2022.

\bibitem{zbMATH07859130}
Dirk Pauly and Walter Zulehner.
\newblock The elasticity complex: compact embeddings and regular
  decompositions.
\newblock {\em Appl. Anal.}, 102(16):4393--4421, 2023.

\bibitem{rapetti2009whitney}
Francesca Rapetti and Alain Bossavit.
\newblock Whitney forms of higher degree.
\newblock {\em SIAM J. Numer. Anal.}, 47(3):2369--2386, 2009.

\bibitem{rogers2006degree}
Luke~G. Rogers.
\newblock Degree-independent {Sobolev} extension on locally uniform domains.
\newblock {\em J. Funct. Anal.}, 235(2):619--665, 2006.

\bibitem{runstsickel}
Thomas Runst and Winfried Sickel.
\newblock {\em Sobolev spaces of fractional order, {Nemytskij} operators and
  nonlinear partial differential equations}, volume~3 of {\em De Gruyter Ser.
  Nonlinear Anal. Appl.}
\newblock Berlin: de Gruyter, 1996.

\bibitem{rychkov1999restrictions}
Vyacheslav~S. Rychkov.
\newblock On restrictions and extensions of the {Besov} and
  {Triebel}-{Lizorkin} spaces with respect to {Lipschitz} domains.
\newblock {\em J. Lond. Math. Soc., II. Ser.}, 60(1):237--257, 1999.

\bibitem{schoeberl2008posteriori}
Joachim Sch{\"o}berl.
\newblock A posteriori error estimates for {Maxwell} equations.
\newblock {\em Math. Comput.}, 77(262):633--649, 2008.

\bibitem{scott1995Lp}
Chad Scott.
\newblock {$L^p$} theory of differential forms on manifolds.
\newblock {\em Transactions of the American Mathematical Society},
  347(6):2075--2096, 1995.

\bibitem{seeley1964extension}
R.~Seeley.
\newblock Extension of {{\(C^ \infty\)}} functions defined in a half space.
\newblock {\em Proc. Am. Math. Soc.}, 15:625--626, 1964.

\bibitem{shi2024estimates}
Ziming Shi and Liding Yao.
\newblock New estimates of {Rychkov}'s universal extension operator for
  {Lipschitz} domains and some applications.
\newblock {\em Math. Nachr.}, 297(4):1407--1443, 2024.

\bibitem{stein1970singular}
Elias~M. Stein.
\newblock {\em Singular integrals and differentiability properties of
  functions}, volume~30 of {\em Princeton Math. Ser.}
\newblock Princeton University Press, Princeton, NJ, 1970.

\bibitem{stocker_zieschang}
Ralph St{\"o}cker and Heiner Zieschang.
\newblock {\em Algebraische {Topologie}: {Eine} {Einf{\"u}hrung}}.
\newblock Stuttgart: B. G. Teubner, 2nd expanded edition, 1994.

\bibitem{taylor2}
Michael~E. Taylor.
\newblock {\em Partial differential equations {II}. {Qualitative} studies of
  linear equations}, volume 116 of {\em Appl. Math. Sci.}
\newblock Cham: Springer, 3rd corrected and expanded edition, 2023.

\bibitem{treves1}
Fran{\c{c}}ois Tr{\`e}ves.
\newblock {\em Introduction to Pseudodifferential and {Fourier} Integral
  Operators}, volume~1 of {\em The University Series in Mathematics}.
\newblock Plenum Press, New York and London, 1980.
\newblock Pseudodifferential Operators.

\bibitem{whitney1934analytic}
Hassler Whitney.
\newblock Analytic extensions of differentiable functions defined in closed
  sets.
\newblock {\em Trans. Am. Math. Soc.}, 36:63--89, 1934.

\bibitem{wong}
M.~W. Wong.
\newblock {\em An introduction to pseudo-differential operators}.
\newblock Singapore: World Scientific, 1991.

\end{thebibliography}

\end{document}